\PassOptionsToPackage{sort&compress}{natbib}
\documentclass[3p,12pt]{elsarticle}

\usepackage{amssymb}
\usepackage{amsmath}

\usepackage[utf8]{inputenc}
\usepackage[hidelinks]{hyperref}
\usepackage[section]{placeins}
\usepackage{subcaption}
\usepackage{xcolor}

\newtheorem{remark}{Remark}

\journal{Computer Methods in Applied Mechanics and Engineering}

\begin{document}

\begin{frontmatter}



\title{A dual lumping procedure for static condensation in mixed NURBS-based isogeometric elements with optimal convergence rates for arbitrary open knot vectors} 

 \affiliation{organization={Chair of Applied Mechanics (LTM), RPTU~University~Kaiserslautern-Landau},
             addressline={Gottlieb-Daimler-Str.},
             city={Kaiserslautern},
             postcode={67663},
             country={Germany}}

\author{Lisa Stammen}
\author{Wolfgang Dornisch}

\begin{abstract}

Locking is a common effect in finite element and isogeometric analysis. In the case of plates, transverse shear locking is most prominent, for shells several other types of locking exist. A common cure are mixed methods that introduce additional fields of unknowns into the variational formulation. These fields reduce constraints and thus alleviate locking significantly. As a drawback, the discretized additional fields increase computational costs significantly. These fields are often eliminated by static condensation, which requires the inverse of a part of the stiffness matrix. In Lagrange-based finite elements, this inverse is computed on element level, due to a discontinuous interpolation of additional fields. Since isogeometric analysis features higher continuity, static condensation must be performed on patch level, which requires a costly matrix inversion on that level. In this contribution, the virtual shear parameters of a mixed isogeometric plate formulation are interpolated by enhanced approximate dual basis functions. This allows to conduct row-sum lumping of the relevant matrix part at a minimal loss of accuracy, since this part becomes diagonal dominant. For a properly chosen integration space, this lumped matrix becomes the identity matrix. Thus, the proposed condensation procedure does not require an inversion anymore. The crucial and novel point is the proposed treatment of knot vectors with limited internal continuity. With the help of several single- and multi-patch examples, both with full and with limited internal continuity, we show that the proposed procedure obtains optimal error convergence rates in all cases, while without these alterations, convergence rates are significantly deteriorated.

\end{abstract}

\begin{highlights}
	\item Mixed isogeometric Reissner-Mindlin plate formulation to alleviate shear locking
	\item Approximate/enhanced approximate dual basis functions for multi-patch NURBS surfaces
	\item Efficient static condensation based on enhanced approximate dual basis functions
	\item Treatment of single- and multi-patch discretizations with limited internal continuity
	\item Optimal error convergence rates for single- and multi-patch discretizations
\end{highlights}


\end{frontmatter}



\section{Introduction}
\label{sec:introduction}

Isogeometric analysis was introduced in \cite{HughesEtAl2005, CottrellEtAl2009} in order to unify computer-aided design and finite element analysis by employing the shape functions of the geometry representation as basis functions within the analysis as well.
Nevertheless, locking phenomena can occur within both finite element analysis and isogeometric analysis and deteriorate the accuracy of results by various specific effects \cite{BieberEtAl2018}.
Thus, \cite{BombardeEtAl2024} reviews different methods to tackle shear locking, membrane locking, volumetric locking and trapezoidal locking for beams, plates and shells in the context of finite element analysis; In particular, assumed natural strain and enhanced (assumed) strain methods as well as displacement-stress-mixed elements and (selectively) reduced integration procedures are studied and compared to conventional finite element methods.
A class of purely displacement-based methods aiming at circumventing locking theoretically for arbitrary discretizations based on non-adapted approximation spaces of equal order has been proposed in~\cite{EchterEtAl2013,OesterleEtAl2016,OesterleEtAl2017,ThiererEtAl2024}. The hierarchical concept for the analysis of shells of \cite{EchterEtAl2013} combines features of Kirchhoff-Love-, Reissner-Mindlin- and solid shells and eliminates shear locking based on separating the kinematics of bending and shear deformations. A slightly different approach is followed in \cite{BieberEtAl2018}, where a specific shear-displacement-type field originating from the strain construction is incorporated for several types of elements.
In order to counteract shear locking, other alternative Reissner-Mindlin-type discretizations, that employ shear strains instead of rotation parameters for isogeometric plates \cite{BeiraoDaVeigaEtAl2015} or involve a shear vector within the discretization of subdivision shells \cite{LongEtAl2012}, have been proposed as well.

A general overview on different mixed finite element formulations for different types of problems, including Reissner-Mindlin plate formulations for linear elasticity, can be found in~\cite{AuricchioEtAl2017}.
There are also numerous contributions investigating the ability of adapted approximation spaces to tackle locking phenomena for Reissner-Mindlin plates based on mixed formulations involving stress parameters in the context of finite element analysis~\cite{Arnold1991,ArnoldBrezzi1993} and isogeometric analysis \cite{BeiraoDaVeigaEtAl2012a,KikisKlinkel2022} as well as isogeometric Reissner-Mindlin shell formulations in mixed~\cite{KikisKlinkel2022} and primal~\cite{KikisEtAl2019} form.
A special focus of investigations is also the  exact fulfillment of the Kirchhoff limits \cite{LongEtAl2012,BeiraoDaVeigaEtAl2012a}.
Various formulations involving different types of assumed strain, force or stress fields that are derived based on the Hu-Washizu functional have also been carried out in order to counteract various types of locking like shear locking for Reissner-Mindlin plates~\cite{SimoHughes1986,WangHu2006}, solid shells \cite{KlinkelEtAl2006,CaseiroEtAl2014} or beams \cite{ChoiEtAl2023, ChoiEtAl2025}. Additionally, the approaches in~\cite{KlinkelEtAl2006, KikisKlinkel2022, ChoiEtAl2023, ChoiEtAl2025} involve a static condensation at element-level that is enabled by employing discontinuous discretizations for the additional field for~\cite{KikisKlinkel2022, ChoiEtAl2023, ChoiEtAl2025}, while a local projection method is employed in \cite{CaseiroEtAl2014}. Furthermore, the formulations in~\cite{KlinkelEtAl2006, CaseiroEtAl2014, KikisEtAl2019, KikisKlinkel2022, ChoiEtAl2023, ChoiEtAl2025} also employ adapted approximation spaces or selectively reduced approximation orders.
Continuous-assumed-strain elements that employ special discretizations of the respective strain parts have been derived for quadratic NURBS in order to tackle membrane locking for Kirchhoff rods \cite{CasqueroGolestanian2022}, linear \cite{CasqueroMathews2023} and non-linear \cite{MathewsCasquero2024} Kirchhoff-Love shells as well as to additionally counteract shear locking for Timoshenko rods~\cite{GolestanianCasquero2023} or volumetric locking for problems involving nearly-incompressibility \cite{CasqueroGolestanian2024}. In \cite{GolestanianCasquero2025}, these approaches have been generalized for cubic NURBS.
Furthermore, \cite{MathewsCasquero2024,GolestanianCasquero2025} also compare a symmetric to a non-symmetric variant of the respective formulation, where a higher accuracy of a non-symmetric variant was obtained compared to the symmetric variant of the formulation \cite{GolestanianCasquero2025}.
Lumped-assumed-strain elements that are equivalent to an assumed strain method based on a strain projection to a reduced-order space employing a lumped mass matrix were introduced in \cite{FaruqueCasquero2024} in order to alleviate membrane locking and shear locking for isogeometric Timoshenko rods.

The correspondence of displacement-based formulations involving reduced integration or selective integration to some mixed formulations was investigated in \cite{MalkusHughes1978} for various types of formulations.
Thus, \cite{ChinosiLovadina1995} studied mixed formulations based on a split of the shear energy that result in a partial selectively reduced integration procedure in the context of Reissner-Mindlin plates.
In order to counteract shear locking and membrane locking for isogeometric Reissner-Mindlin plates and shells, \cite{AdamEtAl2015a} proposes specific continuity- and approximation order-dependent reduced integration rules for B-splines of high continuity.
In \cite{AdamEtAl2015b}, optimized efficient element-based selective integration and reduced integration procedures that utilize the continuity of NURBS are derived by employing different approximations of the integration spaces in order to counteract numerical and volumetric locking in problems of nearly-incompressibility for isogeometric Reissner-Mindlin plates and shells.

In \cite{BouclierEtAl2012}, the existence of equal mixed forms and $\bar{B}$-projection methods that involve lower-order spaces was determined for isogeometric beam formulations.
Different projection methods that are also able to tackle shear locking are investigated in \cite{ElguedjEtAl2008} within the context of isogeometric analysis for problems of nearly-incompressibility.
Within the scope of projection-based isogeometric methods, that have been derived from corresponding mixed formulations in order to tackle different locking phenomena, several research also focused on reducing the computational effort of the required inversion by aiming at a preferable matrix structure by implementing various adapted interpolations. This involves local projection procedures for shell formulations \cite{BouclierEtAl2013, GrecoEtAl2018} as well as elasticity formulations \cite{AntolinEtAl2017}, that can be based on the use of discontinuous fields \cite{GrecoEtAl2018, AntolinEtAl2017} and may lead to purely displacement-based formulations \cite{BouclierEtAl2013, AntolinEtAl2017}. Additionally, the investigations in \cite{AntolinEtAl2017} also include a potential lumping and a comparison of symmetric and non-symmetric variants of the formulation, where the non-symmetric variants lead to improved accuracy compared to the symmetric variants as well.
A comparison of a symmetric and non-symmetric finite element methods regarding accuracy and computational efficiency is conducted in \cite{EisentraegerEtAl2026}, where a specific combination of basis function types proves advantageous for analyzing distorted meshes and more complex benchmark examples but still requires further locking treatment and goes along with increased computational costs.

While \cite{LamichhaneEtAl2006, DjokoEtAl2006} make use of an orthogonal decomposition of stresses within the analysis of the well-posedness and stability in connection with the ability to counteract locking of the proposed modified displacement-strain-stress mixed formulation with respect to incompressibility in plane elasticity,
\cite{KlinkelEtAl2008} construct a partially orthogonal and discontinuous interpolation of strains and stresses in order to achieve a decoupling on element level to enable an efficient static condensation procedure within a non-linear enhanced three-field finite element shell formulation.
Thus, the orthogonality of basis functions can be exploited for the efficient static condensation of additional unknowns, leading to various investigations employing different types of dual basis functions.
Dual functions for B-Splines have been defined explicitly in \cite{DeBoor1975, DeBoor1990} as well as in~\cite{Schumaker2007}.
Local dual bases for Lagrange shape functions that are computed involving the inverse Gram matrix are introduced in \cite{OswaldWohlmuth2002} and employed for the Lagrange multiplier space within a finite element mortar framework.
Approximate dual functions for B-splines that provide a minimal support have been proposed in \cite{ChuiEtAl2004} in the context of harmonic analysis.
\cite{DornischEtAl2017}~compared the aforementioned three types of dual functions and derived approximate duals for NURBS within the scope of a mortar method for patch-coupling, while dual functions for NURBS have been proposed in \cite{SeitzEtAl2016}  within an isogeometric mortar contact formulation.
Dual basis functions are also employed for the Lagrange multiplier spaces within Lagrange-based finite element mortar formulations \cite{OswaldWohlmuth2002, LamichhaneWohlmuth2002, LamichhaneWohlmuth2004, LamichhaneWohlmuth2007} and applied for the coupling of multiple discretizations \cite{Wohlmuth2000}, material behavior \cite{Lamichhane2009}, as well as for mortar-based contact formulations \cite{WohlmuthEtAl2012}. Especially in NURBS-based isogeometric analysis, patch coupling is of high relevance. Thus, the use of dual mortar formulations in this context is proposed in \cite{DornischEtAl2017,ZouEtAl2018,WunderlichEtAl2019,MiaoEtAl2020}.
Furthermore, investigations related to the size of the support of the different kinds of dual functions have been conducted in \cite{LamichhaneWohlmuth2002, LamichhaneWohlmuth2007, DornischEtAl2017, MiaoEtAl2020}.
Thus, such formulations can make use of the orthogonality in order to enable an efficient inversion within a static condensation procedure by using different types of dual basis functions \cite{Lamichhane2009}, due to a diagonalized \cite{Wohlmuth2000, LamichhaneWohlmuth2007} or sparser matrix structure~\cite{LamichhaneWohlmuth2004, MiaoEtAl2020} of the relevant matrix.
Employed accordingly, such function spaces can be used to enable an efficient static condensation for mixed isogeometric shell formulations~\cite{ZouEtAl2020,NguyenEtAl2024} that counteract different locking phenomena as well.
In particular, efficient static condensation procedures for mixed formulations, that are based on a sparser matrix structure or a diagonalized matrix part achieved by the employment of dual basis functions, can result in symmetric or non-symmetric formulations like compared in \cite{MiaoEtAl2018} in order to tackle shear locking for Timoshenko beams and volumetric locking for nearly-incompressible plate benchmark examples within isogeometric analysis as well as in~\cite{Lamichhane2008b, LamichhaneStephan2012} to counteract locking for problems of nearly-incompressibility within finite element analysis.
In \cite{LamichhaneStephan2012,Lamichhane2013}, a three-field functional is employed, in order to enable a symmetric mixed system matrix if dual basis functions are used, whereas the non-symmetric variant of the formulation in \cite{MiaoEtAl2018} yields an improved convergence in contrast to its symmetric variant also for these investigated plate benchmark examples.
Specifically,
quasi-orthogonal basis functions are involved in \cite{LamichhaneWohlmuth2004,Lamichhane2008a,Lamichhane2013} to increase the sparsity of the system matrix and simplify the computation of the inverse within the static condensation procedure. Furthermore, \cite{Wohlmuth2000, OswaldWohlmuth2002, LamichhaneWohlmuth2002, LamichhaneWohlmuth2004, LamichhaneWohlmuth2007, Lamichhane2014, MiaoEtAl2018, WunderlichEtAl2019} also focus on the locality of the approach.
Such procedures can also be applied for mass lumping within isogeometric formulations for dynamic problems~\cite{NguyenEtAl2023,HeldEtAl2024,HiemstraEtAl2025}.

In \cite{StammenDornisch2023}, we investigated that, for a displacement-shear-mixed isogeometric Reissner-Mindlin plate formulation, in order to counteract shear locking, the determined combinations of selectively reduced approximation orders, that would simultaneously lead to consistent degrees within the shear-strain equation and for the constitutive equation, should only be employed for the mixed parameters and not within the irreducible part, as this would, while effectively counteracting shear locking, deteriorate the expected convergence rate. Based on this formulation, the present contribution focuses on enabling an efficient static condensation procedure to eliminate the additional shear parameters. To reduce the computational effort for the matrix inversion that is required within the computation of the condensed system of equations, lumping the according matrix part to a diagonal matrix would be most efficient. However, as will be shown within the scope of the investigated numerical examples, performing row-sum lumping of this sub-matrix for a standard NURBS-based interpolation introduces a huge error and deteriorates the convergence rate. In order to minimize these effects, employing (approximate) dual basis functions of NURBS as trial functions is proposed within this contribution. Furthermore, the influence of the presence of points with limited internal continuities is studied for single-patch and multi-patch discretizations of a benchmark example with analytical solution. As these investigations will reveal the requirement for an adapted dual condensation procedure, employing enhanced approximate dual basis functions, that have been lately introduced in \cite{Stoeckler2025}, is shown to circumvent deteriorated convergence rates for the present investigations. In addition, it is shown that the weights associated with the NURBS basis functions employed for the interpolation of the additional shear parameters can be dropped at a minimal loss of accuracy, enabling a more efficient implementation of the proposed dual procedure.
Moreover, investigations on dual basis functions for T-Splines \cite{BeiraoDaVeigaEtAl2012b} indicate that the present concept could also be implemented for other types of spline basis functions, that would also enable local refinement.

\section{NURBS and their dual basis functions}
\label{sec:basis_functions}
In this section, we provide the basic formulas for B-splines and NURBS, in order to establish the notation for the introduction of their approximate dual basis functions.
\subsection{Construction of B-splines}
The construction of uni-variate and bi-variate B-splines is conducted according to \cite{PieglTiller1997}. B-spline curves of degree $p$ are constructed based on a non-decreasing knot vector
\begin{equation}
	\boldsymbol{\Xi}=\left(\xi_1,\dots,\xi_{n+p+1}\right)
	\text{,}\quad
	\xi_i\leq\xi_{i+1}
	\text{,}\quad
	i=1,\dots,n+p\,.
\end{equation}
Using the well-known Cox-de Boor recurrence algorithm, based on
\begin{equation}
	\begin{split}
	N_{i}^0\left(\xi\right)&=
	\left\{
	\begin{array}{ll}
		1\text{,} & \xi_i\leq\xi<\xi_{i+1}\\
		0\text{,} & \text{otherwise}
	\end{array}
	\right.
	\text{and}
	\\
	N_{i}^p\left(\xi\right)&=\frac{\xi-\xi_i}{\xi_{i+p}-\xi_i}\cdot N_{i}^{p-1}\left(\xi\right)+\frac{\xi_{i+p+1}-\xi}{\xi_{i+p+1}-\xi_{i+1}}\cdot N_{i+1}^{p-1}\left(\xi\right)\text{,}
	\end{split}
\end{equation}
the B-spline curves and their corresponding derivatives
\begin{equation}
	N_{i,\xi}^p\left(\xi\right)=\frac{p}{\xi_{i+p}-\xi_i}\cdot N_{i}^{p-1}\left(\xi\right)-\frac{p}{\xi_{i+p+1}-\xi_{i+1}}\cdot N_{i+1}^{p-1}\left(\xi\right)\text{\,}
\end{equation}
can be determined. Uni-variate NURBS basis functions are computed by
\begin{equation}
    R^p_i\left(\xi\right)=\frac{N^p_i\left(\xi\right)\cdot w_i}{\omega(\xi)}\qquad\text{with}\qquad\omega(\xi)=\sum_{k=1}^n{N^p_k\left(\xi\right)\cdot w_k}\,.
\end{equation}
 Within this work, we use the symbol $(\cdot)$ in the sense of a matrix multiplication.
Based on a second knot vector
$\boldsymbol{H}=\left(\eta_1,\dots,\eta_{m+q+1}\right)\text{,}~\eta_j\leq\eta_{j+1}\text{,}~ j=1,\dots,m+q$,
bi-variate B-spline basis functions, with respect to the parametric coordinates $\xi$ and $\eta$, can be constructed as follows:
\begin{equation}
	N_{ij}^{p,q}\left(\xi,\eta\right)
	=N_{i}^{p}\left(\xi\right)\cdot N_{j}^{q}\left(\eta\right)
\end{equation}
Therefore, independent approximation orders $p$ and $q$ can be selected for the interpolation of the relevant quantities in $\xi$- and in $\eta$-direction, respectively.
Constructing bi-variate NURBS basis functions additionally involves their corresponding weights as follows:
\begin{equation}
\label{eq:NURBS_basis}
	R_{ij}^{p,q}\left(\xi,\eta\right)
	=\frac{N_{ij}^{p,q}\left(\xi,\eta\right)\cdot w_{ij}}	{W\left(\xi,\eta\right)}
    \qquad\text{with}\qquad
    W\left(\xi,\eta\right)=\sum_{k=1}^{n}\sum_{l=1}^{m}N_{k}^{p}\left(\xi\right)\cdot N_{l}^{q}\left(\eta\right)\cdot w_{kl}
\end{equation}
In the remainder of the paper, we usually omit to display the dependency of $W$ and the basis functions on the parametric coordinates $(\xi,\eta)$ for the sake of readability.

\subsection{(Approximate) dual basis functions}
\label{sec:approximate_duals}
Approximate dual basis functions for B-Splines are constructed as introduced in \cite{ChuiEtAl2004}.
The eponymous property of dual basis functions ${\lambda}_j^p$ is the bi-orthogonality
\begin{equation}
	f_j\left(N_i^p\right)=\int_{I_j}N_i^p\cdot{\lambda}_j^p~\text{ds}=\delta_{ij}\,,
\end{equation}
where $I_j$ represents the support interval of the function with the index $j$.
For approximate dual basis functions $\tilde{\lambda}_j^p$, this criterion is fulfilled only approximately, i.e.
\begin{equation}
	f_j\left(N_i^p\right)=\int_{I_j}N_i^p\cdot\tilde{\lambda}_j^p~\text{ds}\approx\delta_{ij}\,,
\end{equation}
but their degree of reproduction $r$  of polynomials
\begin{equation}
	x^r=\sum_{j=1}^{n}\left(\int_{I_j}x^r\cdot\tilde{\lambda}_j^p~\text{ds}\right)\cdot N_j^p
	\text{,}\quad
	0\leq r\leq p
\end{equation}
can be selected as $r=p$, which enables full reproduction of B-spline basis functions of degree $p$ ~\cite{ChuiEtAl2004, DornischEtAl2017}. According to~\cite{ChuiEtAl2004},
 the vector $\tilde{\boldsymbol{\lambda}}^p_r\left(\xi\right)$ of uni-variate approximate dual basis functions of B-splines at the coordinate $\xi$ can be computed by
\begin{equation}
	\label{eq:duals1D}
	\tilde{\boldsymbol{\lambda}}^p_r\left(\xi\right)^T=\boldsymbol{S}_r^p\left(\boldsymbol{\Xi}\right)\cdot\boldsymbol{N}^p\left(\xi\right)^T\,,
\end{equation}
where $\tilde{\boldsymbol{\lambda}}^p_r\left(\xi\right)=\left(\tilde{\lambda}_1^p\left(\xi\right),\dots,\tilde{\lambda}_n^p\left(\xi\right)\right)$ and  $\boldsymbol{N}^p\left(\xi\right)=\left(N_1^p\left(\xi\right),\dots,N_n^p\left(\xi\right)\right)$ contain the particular values of the $n$ approximate dual and B-spline basis functions, respectively.
According to \cite{HeldEtAl2024}, the transformation matrix $\boldsymbol{S}^p_r\left(\boldsymbol{\Xi}\right)$ can be determined as follows:
\begin{equation}
	\label{eq:dual_transformation_matrix}	\boldsymbol{S}_r^p\left(\boldsymbol{\Xi}\right)=\boldsymbol{U}_{\boldsymbol{\Xi},0}+\sum_{v=1}^{r}\left[\left(\prod_{k=1}^{v}\boldsymbol{D}_{\boldsymbol{\Xi},p+k}\right)\cdot
	\boldsymbol{U}_{\boldsymbol{\Xi},v}\cdot
	\left(\prod_{k=1}^{v}\boldsymbol{D}_{\boldsymbol{\Xi},p+k}\right)^T\right]
\end{equation}
The formulas for the determination of the matrices involved in this equation are contained in \citep{ChuiEtAl2004}. It is to be noted, that the computation of $\boldsymbol{S}_r^p\left(\boldsymbol{\Xi}\right)$ needs to be performed only once for every knot vector, since it only depends on the knot vector $\boldsymbol{\Xi}$. Despite the rather complicated notation, the computation of $\boldsymbol{S}^p_r\left(\boldsymbol{\Xi}\right)$ is fast and robust, due to their analytical and explicit formulation. The symmetric matrix $\boldsymbol{S}^p_r\left(\boldsymbol{\Xi}\right)$ has the size $n\times n$ and is banded with a bandwidth $r$. Thus, up to $p+2r+1$ approximate dual basis functions $\tilde{\lambda}_i^p\left(\xi\right)$ have support in each element.
According to~\cite{DornischEtAl2017}, uni-variate approximate dual basis functions can also be constructed for NURBS by employing the adapted formula
\begin{equation}	{}^R\tilde{\boldsymbol{\lambda}}_r^p\left(\xi\right)^T=\boldsymbol{\Lambda}^{-1}\cdot\boldsymbol{S}_r^p\left(\boldsymbol{\Xi}\right)\cdot\boldsymbol{\Lambda}^{-1}\cdot\boldsymbol{R}^p\left(\xi\right)^T
	\text{,}
\end{equation}
where the diagonal matrix $\boldsymbol{\Lambda}=\operatornamewithlimits{diag}\left(w_j\right)$ for $j=1,\dots,n$ contains the weights and $\boldsymbol{R}^p\left(\xi\right)$ is the vector of the NURBS basis functions.
This also limits the reproduction of polynomials as follows:
\begin{equation}
\frac{x^r}{\omega}=\sum_{j=1}^{n}\left(\int_{I_j}\frac{x^r}{\omega}\cdot R_j^p\cdot \omega^2~\text{ds}\right)\cdot{}^R\tilde{\lambda}_j^p
\text{,}\quad
0\leq r\leq p
\end{equation}
For B-splines, employing Eq.~\eqref{eq:duals1D} and the corresponding definitions of the required matrices for each of the both parametric directions $\xi$ and $\eta$ separately, two uni-variate approximate dual basis functions $\tilde{\boldsymbol{\lambda}}_\xi^{p}\left(\xi\right)$ and $\tilde{\boldsymbol{\lambda}}_\eta^{q}\left(\eta\right)$ are determined by
\begin{equation}\tilde{\boldsymbol{\lambda}}_\xi^{p}\left(\xi\right)^T=\boldsymbol{S}_{r}^{p}\left(\boldsymbol{\Xi}\right)\cdot\boldsymbol{N}_\xi^{p}\left(\xi\right)^T
	\quad\text{and}\quad\tilde{\boldsymbol{\lambda}}_\eta^{q}\left(\eta\right)^T=\boldsymbol{S}_{s}^{q}\left(\boldsymbol{H}\right)\cdot\boldsymbol{N}_\eta^{q}\left(\eta\right)^T\,
\end{equation}
This also preserves the possibility to select both separate interpolation orders ($p$ and $q$) and degrees of reproduction ($r$ and $s$) for each direction. Based on these, bi-variate approximate  dual basis functions can be constructed by computing the tensor product
\begin{equation}	\tilde{\boldsymbol{\lambda}}^{p,q}_{r,s}\left(\xi,\eta\right)^T=\tilde{\boldsymbol{\lambda}}_\xi^{p}\left(\xi\right)^T\otimes\tilde{\boldsymbol{\lambda}}_\eta^{q}\left(\eta\right)^T	=\Bigl(\boldsymbol{S}_{r}^{p}\left(\boldsymbol{\Xi}\right)\cdot\boldsymbol{N}_\xi^{p}\left(\xi\right)^T\Bigr)\otimes\Bigl(\boldsymbol{S}_{s}^{q}\left(\boldsymbol{H}\right)\cdot\boldsymbol{N}_\eta^{q}\left(\eta\right)^T\Bigr)\,.
\end{equation}
Using the 2D-transformation matrix~$\boldsymbol{S}^{p,q}_{r,s}\left(\boldsymbol{\Xi},\boldsymbol{H}\right)$ and the vector of bi-variate B-spline basis functions~$\boldsymbol{N}^{p,q}\left(\xi,\eta\right)$, this expression can be rearranged to
\begin{eqnarray}
		\tilde{\boldsymbol{\lambda}}^{p,q}_{r,s}\left(\xi,\eta\right)^T
		=\Bigl(\boldsymbol{S}_{r}^{p}\left(\boldsymbol{\Xi}\right)\otimes\boldsymbol{S}_{s}^{q}\left(\boldsymbol{H}\right)\Bigr)\cdot\Bigl(\boldsymbol{N}_\xi^{p}\left(\xi\right)^T\otimes\boldsymbol{N}_\eta^{q}\left(\eta\right)^T\Bigr)
		=\boldsymbol{S}^{p,q}_{r,s}\left(\boldsymbol{\Xi},\boldsymbol{H}\right)\cdot\boldsymbol{N}^{p,q}\left(\xi,\eta\right)^T
\end{eqnarray}
in order to regain the same structure as for the uni-variate dual basis function in Eq.~\eqref{eq:duals1D} also for the bi-variate case.
This enables to maintain the advantage, pointed out by \cite{HeldEtAl2024} for the uni-variate case, that the transformation matrix needs to be determined only once (globally) and can then be pre-multiplied to the system matrices after conducting the entire assembly procedure. This allows to use B-splines (or NURBS) as both shape functions and trial functions within the initial computations, which reduces both the implementational and computational effort considerably.
For NURBS basis functions, the vector of bi-variate approximate dual basis functions~${}^R\tilde{\boldsymbol{\lambda}}^{p,q}_{r,s}\left(\xi,\eta\right)$ is computed from the vector of bi-variate NURBS basis functions~$\boldsymbol{R}^{p,q}\left(\xi,\eta\right)$ by
\begin{equation}
\label{eq:2d_NURBS_transformation_matrix}
    {}^R\tilde{\boldsymbol{\lambda}}^{p,q}_{r,s}\left(\xi,\eta\right)^T
		={}^R\boldsymbol{S}^{p,q}_{r,s}\cdot\boldsymbol{R}^{p,q}\left(\xi,\eta\right)^T
        \qquad\text{with}\qquad        {}^R\boldsymbol{S}^{p,q}_{r,s}=\boldsymbol{\Lambda}^{-1}\cdot\boldsymbol{S}^{p,q}_{r,s}\left(\boldsymbol{\Xi},\boldsymbol{H}\right)\cdot\boldsymbol{\Lambda}^{-1}\,,
\end{equation}
where the diagonal matrix of weights is now defined by $\boldsymbol{\Lambda}=\operatornamewithlimits{diag}\left(w_{ij}\right)$ for $i=1,\dots,n$ and $j=1,\dots,m$.

\subsection{Comparison}
\begin{figure}[th]
	\centering
	\begin{subfigure}{0.49\textwidth}
		\includegraphics[width=\textwidth]{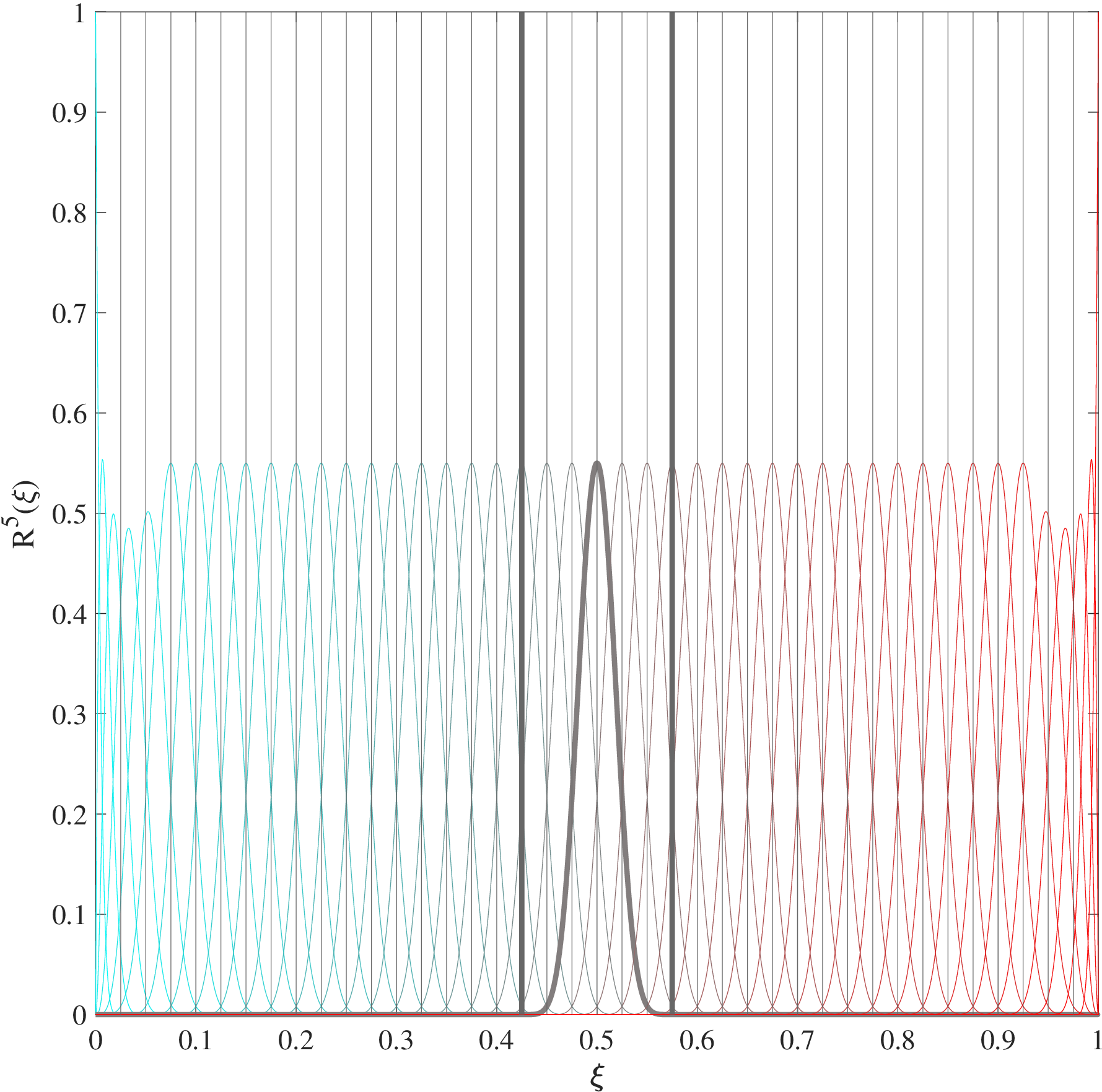}
		\caption{uni-variate NURBS}
	\end{subfigure}
	\hfill
	\begin{subfigure}{0.49\textwidth}
		\includegraphics[width=\textwidth]{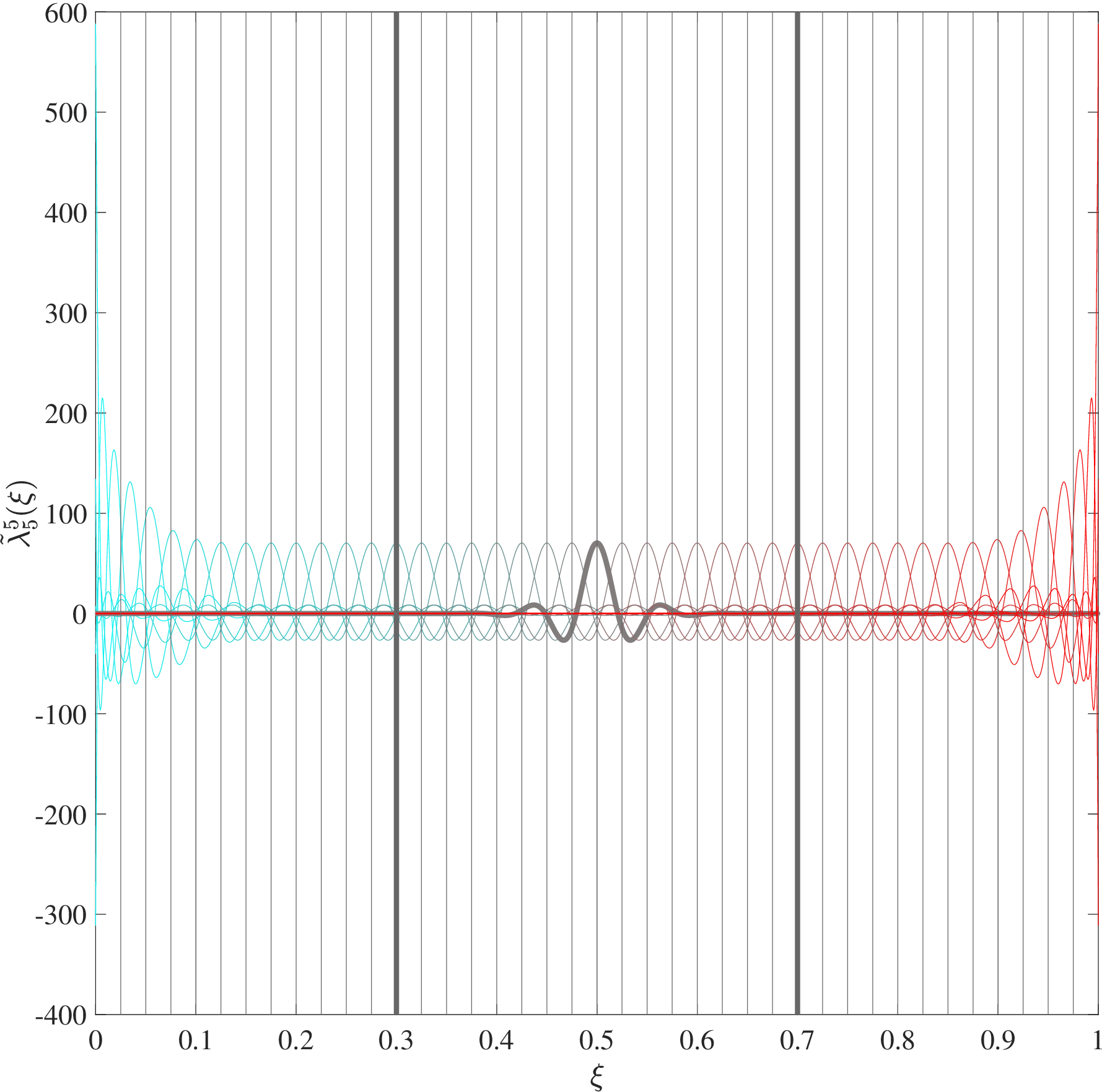}
		\caption{corresponding uni-variate approximate dual basis functions}
	\end{subfigure}
	\\
	\begin{subfigure}{0.49\textwidth}
		\includegraphics[width=\textwidth]{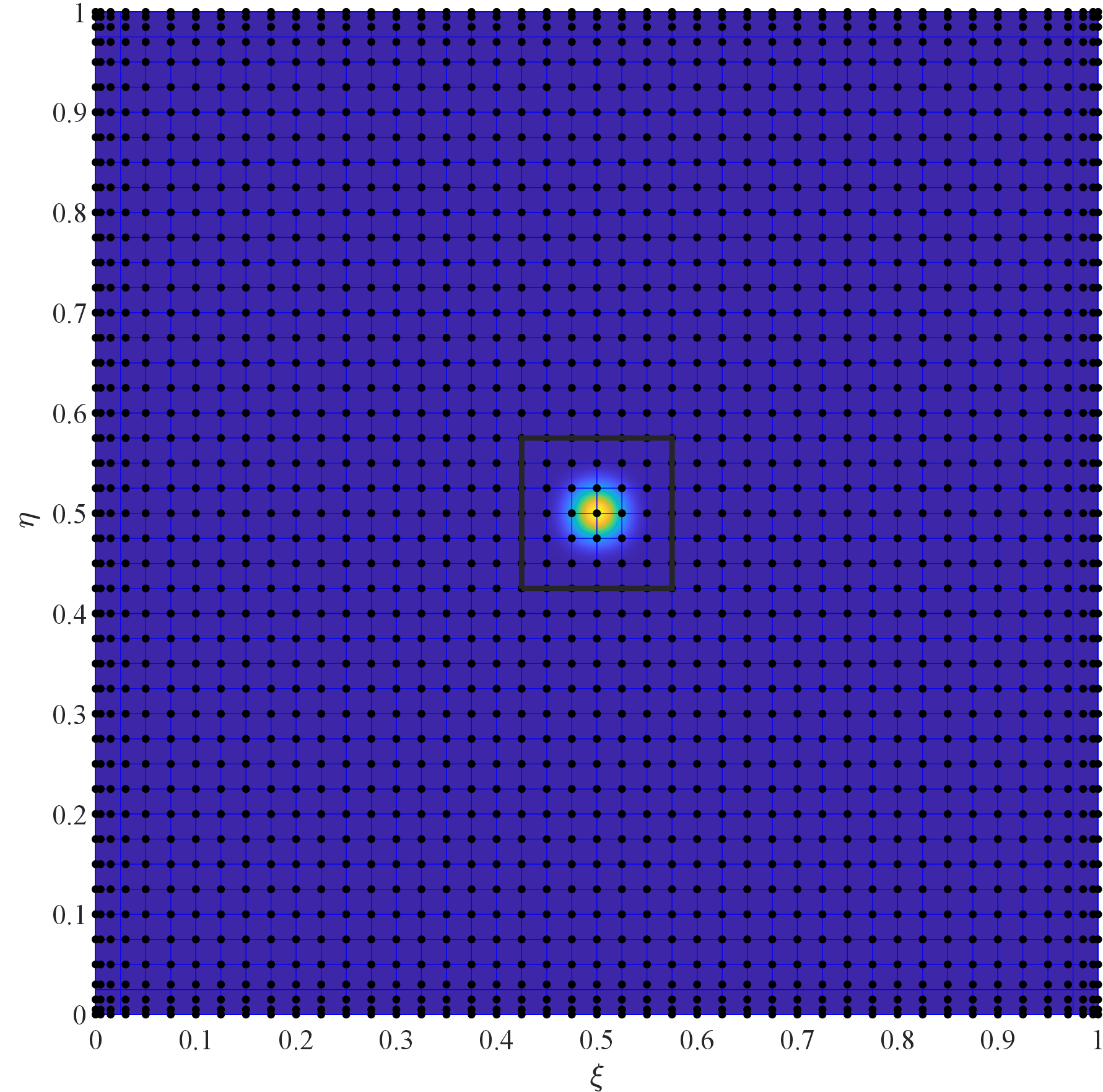}
		\caption{bi-variate B-splines}
	\end{subfigure}
	\hfill
	\begin{subfigure}{0.49\textwidth}
		\includegraphics[width=\textwidth]{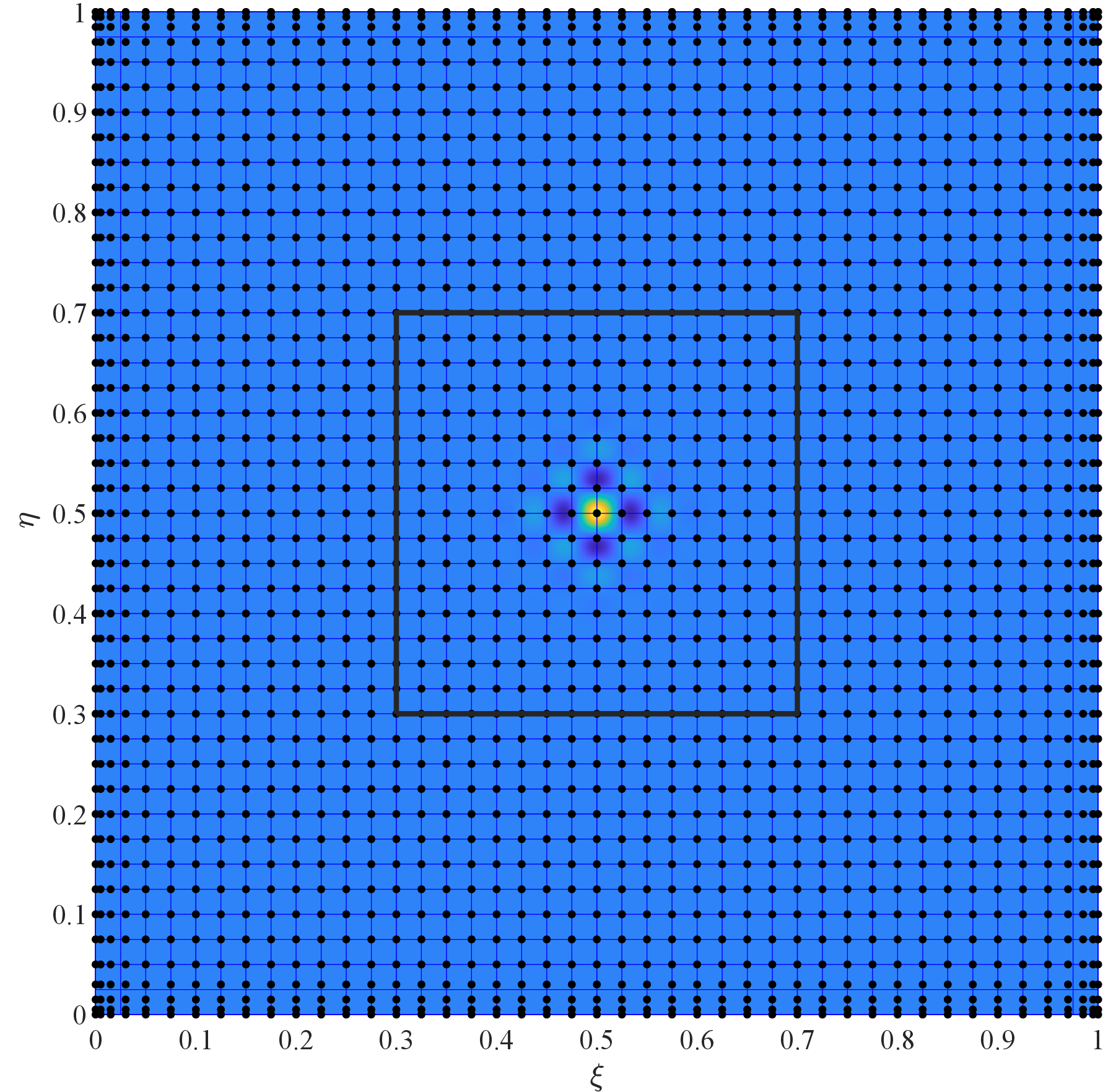}
		\caption{corresponding bi-variate approximate dual basis functions}
	\end{subfigure}
	\caption{Comparison of NURBS and their approximate dual basis functions for $p=q=r=s=5$ using the Quadratic plate with complex loading from Sec.~\ref{sec:numerical_exmples}}
	\label{fig:basisFunctions}
\end{figure}

In Fig.~\ref{fig:basisFunctions}, B-Splines and their approximate dual basis functions are compared for $p=r=5$ using the Quadratic plate with complex loading examined in Sec.~\ref{sec:numerical_exmples}. In Fig.~\ref{fig:basisFunctions}~(a) and (b) the uni-variate shape functions are compared. It becomes obvious that approximate dual basis functions take on negative as well as much higher numeric values compared to B-spline basis functions that have to comply with the partition-of-unity-criterion. Furthermore, approximate dual basis functions with maximum reproduction $r=p$ possess an increased support of $3\cdot p+1$ compared to B-splines which have support in $p+1$ elements. For the highlighted basis function, the support of the resulting bi-variate basis function is indicated in Fig.~\ref{fig:basisFunctions}~(c) and (d), assuming that the same basis functions are employed for the second parametric direction as well.
Employing approximate dual basis functions as test functions while using B-splines as shape functions, a system matrix with higher bandwidth, but diagonal dominance can be achieved. This allows to develop an efficient static condensation procedure for mixed formulations as proposed in Sec.~\ref{sec:procedure}.
Enhanced approximate dual basis functions, which will be introduced in Sec.~\ref{sec:enhancedAD}, only differ from approximate dual basis functions at knots with limited continuity (knot vectors which have internal knots in the coarsest mesh). In that case, the support of enhanced approximate dual basis functions is slightly extended. The resulting functions look very similar to standard approximate dual basis functions. For the example of Fig.~\ref{fig:basisFunctions}, there is no difference between enhanced and standard approximate dual basis functions. Thus, we omit to display the enhanced approximate dual basis functions.

\section{Treatment of patches with limited internal continuity}
\label{sec:internal_continuity}
In the frame of isogeometric analysis, approximate dual basis functions have been used for the first time in~\cite{DornischEtAl2017} in order to obtain a fast and accurate coupling of NURBS patches with the mortar method. The results of~\cite{DornischEtAl2017} show that, despite the application of row-sum lumping, optimal convergence rates of the global stress error can be obtained for all geometries with unlimited internal continuity, i.e., the knot vectors of the coarsest mesh of order $p$ in the first parametric direction are of the type
\begin{equation}
    \boldsymbol{\Xi}=(\boldsymbol{\xi}_1,\boldsymbol{\xi}_2) \quad\text{with}\quad \boldsymbol{\xi}_i=\left(\xi_i\cdot\boldsymbol{1}_{p+1}\right)\,,\quad i\in\{1,2\}\quad\text{and}\quad \xi_2>\xi_1\,,
\end{equation}
where $\boldsymbol{1}_{p+1}$ is a row vector of length $p+1$ of which all entries are equal to $1$. The condition for the knot vector~$\boldsymbol{H}$ in the second parametric direction is defined analogously. Hence, its explicit notation is omitted. If the condition of unlimited internal continuity is not met, which means that the knot vector of the coarsest mesh contains internal knots, the optimal convergence rates of the global stress error cannot be achieved in some cases in~\cite{DornischEtAl2017}. Knot vectors of this type are defined by
\begin{equation}
\label{eq:knotvectorlimitedcontinuity}
    \boldsymbol{\Xi}=(\boldsymbol{\xi}_1,\hat{\boldsymbol{\xi}},\boldsymbol{\xi}_2) \quad\text{with}\quad \hat{\boldsymbol{\xi}}=\left(\hat{\xi}_1\cdot\boldsymbol{1}_{c_1},\ldots,\hat{\xi}_l\cdot\boldsymbol{1}_{c_l}\right)\,,\quad 1\leq j\leq l\,,\quad\hat{\xi}_{j+1}>\hat{\xi}_{j}\,,
\end{equation}
where obviously $\xi_2>\hat\xi_l$ and $\hat\xi_1>\xi_1$ has to hold. Therein, the multiplicity of each internal knot $j$ is denoted by $c_j$. Accordingly, the continuity at the internal knot $j$ is $C^{p-c_j}$. The deterioration of the convergence rate for knot vectors of the type defined in Eq.~\eqref{eq:knotvectorlimitedcontinuity} also occurs for the proposed displacement-stress mixed plate formulation with lumping-based static condensation if approximate dual basis functions are employed, as we will show for the numerical example investigated in Sec.~\ref{sec:numerical_exmples}.
\begin{remark}
    Since the deterioration of the convergence rate bases on the smoothness of the underlying B-splines, we anticipate, that this effect will occur in any similar constraining method that uses approximate dual basis functions or any other linear combination of B-splines or NURBS as test functions. To the authors' knowledge, the numerical examples of all publications on similar methods do not contain knot vectors with limited continuity. Thus, the deterioration of the convergence rate cannot be discovered in those investigations. However, in industrial applications, knot vectors with limited continuity occur quite frequently.
\end{remark}In the following, we provide a short guideline on how to remove this problematic behavior, based on the previous works~\cite{DornischStoeckler2021} and~\cite{Stoeckler2025}. These simple cures can be applied straightforwardly to other works based on approximate duals, such as~\cite{HeldEtAl2024,NguyenEtAl2024, NguyenEtAl2023, NguyenEtAl2025}.
\subsection{Reduction of continuity}
\label{sec:reductioncontinuity}
The deterioration of the convergence rate of the global stress error for knot vectors with limited internal continuities has been studied in~\cite{DornischStoeckler2021} in the context of the isogeometric mortar method for NURBS as test functions. However, since approximate dual basis functions are linear combinations of NURBS (see Eq.~\eqref{eq:2d_NURBS_transformation_matrix}), the findings hold for approximate duals in the same manner. In the following, we will repeat the very simple cure and explain why the problem also occurs for the formulation proposed in this contribution. In~\cite{DornischStoeckler2021}, knot vectors of the type defined in Eq.~\eqref{eq:knotvectorlimitedcontinuity} are studied. The case $c_j<p\quad\forall\quad j\in\{1,\ldots,l\}$, which means that at each knot $j$ at least $C^1$-continuity prevails, only requires increasing the multiplicity to $c_j+1$, i.e., reducing the continuity by one. The resulting knot vector
\begin{equation}
\label{eq:knotvectorenhanced}
    \boldsymbol{\Xi}^+=(\boldsymbol{\xi}_1,\hat{\boldsymbol{\xi}}^+,\boldsymbol{\xi}_2) \quad\text{with}\quad \hat{\boldsymbol{\xi}}^+=\left(\hat{\xi}_1\cdot\boldsymbol{1}_{c_1+1},\ldots,\hat{\xi}_l\cdot\boldsymbol{1}_{c_l+1}\right)\,,
\end{equation}
where the inequalities from Eq.~\eqref{eq:knotvectorlimitedcontinuity} still hold,
has to be used for the discretization of geometry and all unknown fields instead of the original coarsest knot vector. Thus, if mesh refinement or order elevation is used, it has to based on $\boldsymbol{\Xi}^+$ instead of $\boldsymbol{\Xi}$. This just constitutes some kind of a priori mesh refinement and adds only a few additional degrees of freedom. The mathematical proof in~\cite{DornischStoeckler2021} shows that this very simple alteration is sufficient to recover optimal convergence rates for NURBS. For approximate dual basis functions, this improves the situation significantly, but as will be seen in Sec.~\ref{sec:C1_internal}, in some extreme cases, the convergence rate still deteriorates for finer meshes. The effect is mathematically explained  in~\cite{Stoeckler2025}, where a corresponding cure is also proposed. In Sec.~\ref{sec:enhancedAD} of this contribution, the main effects and their remedy will be revisited.
For the case of $c_j=p$, which means that $C^0$-continuity prevails at knot $j$, the interface has to be split into two parts for the mortar method. Within the mortar method, the split of the interface is not problematic, since it does not propagate into the patch. However, for the method proposed in this paper, splitting a patch of the shear parameter discretization into several patches while the deformation mesh remains unaltered would entail a complicated implementation. A much simpler solution for the case $c_j=p$ will be presented in the next subsection.
\begin{remark}
    In the mortar method, the reason for the need to decrease continuity is that the normal vector to the interface has to be interpolated. The normal vector is computed from the first spatial derivatives of the NURBS-defined geometry field. Consider a point $\boldsymbol{x}$ with $C^1$-continuity along the mortar interface. The normal vector will have $C^0$-continuity at the very same point. If the normal vector is interpolated with the same basis functions as the geometry, this interpolation will be of $C^1$-continuity at the point $\boldsymbol{x}$. Thus, the interpolation introduces a considerable error. If the enhanced knot vector of Eq.~\eqref{eq:knotvectorenhanced} is used instead, the normal vector does not change (it is still of $C^0$-continuity at the point $\boldsymbol{x}$ since mesh refinement does not alter the geometry), but the interpolation will now also possess $C^0$-continuity at this point. Thus, the interpolation error is entirely removed for the case of NURBS. Since approximate dual basis functions are linear combinations of NURBS, the interpolation is 'smeared' around point $\boldsymbol{x}$ again, which reintroduces an interpolation error (cf.~Sec.~\ref{sec:enhancedAD}).
\end{remark}
    For the displacement-stress mixed plate formulation with static condensation proposed in this contribution, independent shear forces $S_\alpha$ for $\alpha=1,2$ are interpolated independently and forced to be equivalent to the shear forces derived from the interpolated deformations (cf. third line of Eq.~\eqref{eq:strongform}). Thus, at a $C^1$-point, the interpolated shear forces are interpolated by $C^1$-basis functions, while the derived shear forces depend on the first derivative of the interpolated deformations, which makes them $C^0$ at the very point. Hence, the cure of Eq.~\eqref{eq:knotvectorenhanced} is also required for our proposed formulation.

\begin{remark}
    It is to be noted that, according to the mathematical proof in~\cite{DornischStoeckler2021}, the multiplicity has to be increased at every internal knot as indicated in Eq.~\eqref{eq:knotvectorenhanced}. For the present numerical studies, it proved sufficient to increase the multiplicity only at points with initial $C^1$-continuity. Although the reduction at points with $C^2$-continuity appeared not to be required to obtain the optimal convergence rate, employing Eq.~\eqref{eq:knotvectorenhanced} is recommended to comply with the theoretical requirement.
\end{remark}
\subsection{Enhanced approximate dual basis functions}
\label{sec:enhancedAD}
Summing up the considerations of Sec.~\ref{sec:reductioncontinuity}, two drawbacks remain. At knots with initial $C^1$-continuity, the convergence rate of the global stress error is deteriorated for very fine meshes. Furthermore, at knots with $C^0$-continuity, patches have to be split in order to use approximate dual basis functions. These two drawbacks can be overcome by enhanced approximate dual basis functions, which have been recently proposed in \cite{Stoeckler2025}. These functions are able to approximate functions from bent Sobolev spaces with optimal approximation order, whereas the standard approximate dual basis functions only approximate functions from smooth Sobolev spaces with optimal order. The improved approximation properties yield proper convergence rates also for knot vectors containing knots with $C^0$-continuity. For details we refer to~\cite{Stoeckler2025}. The numerical results in Sec.~\ref{sec:numerical_exmples} will assess the improvement of accuracy induced by using enhanced approximate duals.

The enhanced approximate dual basis functions are also computed from B-splines using an adapted transformation matrix ${}^{enh}\boldsymbol{S}^p_r(\boldsymbol{\Xi}^+)$, that is computed similarly to that for approximate dual basis functions in Eq.~\eqref{eq:dual_transformation_matrix}, but involving an additional term consisting of a matrix product of similar structure~(cf.~Eq.~(19) of~\cite{Stoeckler2025}):
\begin{equation}
    \label{eq:enhanced_dual_transformation_matrix}	{}^{enh}\boldsymbol{S}_r^p\left(\boldsymbol{\Xi^+}\right)=\boldsymbol{S}_r^p\left(\boldsymbol{\Xi^+}\right)+\left[\left(\prod_{\hat{k}=1}^{p}\boldsymbol{D}_{\boldsymbol{\Xi^+},p+\hat{k}}\right)\cdot
	\boldsymbol{U}_{\boldsymbol{\Xi^+},p+1}\cdot
	\left(\prod_{\hat{k}=1}^{p}\boldsymbol{D}_{\boldsymbol{\Xi^+},p+\hat{k}}\right)^T\right]
\end{equation}
For the computation of $\boldsymbol{U}_{\boldsymbol{\Xi^+},p+1}$, it is required to solve a system of equations. See Eq.~(27) of~\cite{Stoeckler2025} for details. The knot vector $\boldsymbol{\Xi}^+$ has to be computed according to Sec.~\ref{sec:reductioncontinuity}, where the continuity of every internal knot with at least $C^1$-continuity has to be lowered by one. Internal knots with $C^0$-continuity do not have to be altered. The resulting matrix has a slightly increased bandwidth.
Once constructed for uni-variate B-splines, the required adaptations regarding their computation for NURBS and bi-variate functions correspond to those conducted for approximate dual basis functions in Sec.~\ref{sec:approximate_duals}. Finally, enhanced approximate dual basis functions for NURBS can be computed as follows:
\begin{equation}
\label{eq:2d_NURBS_transformation_matrix_enh}
    {}^{enh,R}\tilde{\boldsymbol{\lambda}}^{p,q}_{r,s}\left(\xi,\eta\right)^T
		={}^{enh,R}\boldsymbol{S}^{p,q}_{r,s}\cdot\boldsymbol{R}^{p,q}\left(\xi,\eta\right)^T
        \quad\text{with}\quad        {}^{enh,R}\boldsymbol{S}^{p,q}_{r,s}=\boldsymbol{\Lambda}^{-1}\cdot{}^{enh,R}\boldsymbol{S}^{p,q}_{r,s}\left(\boldsymbol{\Xi},\boldsymbol{H}\right)\cdot\boldsymbol{\Lambda}^{-1}
\end{equation}

\section{Derivation of mixed plate formulation with adapted interpolation orders}
\label{sec:derivation}

In this section, an isogeometric displacement-stress-mixed plate formulation with adapted approximation orders as proposed in \cite{StammenDornisch2023} is derived based on a standard mixed plate formulation \cite{Zienkiewicz.2006b}. This will serve as a model problem for the static condensation procedure proposed in the subsequent section. A Reissner--Mindlin plate formulation is chosen, because analytical solutions of sufficiently complex load and stress states exist, which is a prerequisite for a proper assessment of global error convergence rates.

Introducing the differential operators
\begin{equation}
		\boldsymbol{L}^T=
		\left[
		\begin{array}{ccc}
			\frac{\partial}{\partial x} & 0 & \frac{\partial}{\partial y} \\
			0 & \frac{\partial}{\partial y} & \frac{\partial}{\partial x}
		\end{array}
		\right]
		\quad\text{and}\quad
		\boldsymbol{\nabla}
		=
		\left(
		\begin{array}{c}
			\frac{\partial}{\partial x} \\ \frac{\partial}{\partial y}
		\end{array}
		\right)
		\text{,}
\end{equation}
the following general governing equations can be derived for a linear plate subjected to a  distributed transversial load $f^{\textrm{ext}}$:
\begin{equation}
	\begin{split}
		\boldsymbol{\nabla}^T\cdot\boldsymbol{S}+f^{\textrm{ext}}=0\\
		\boldsymbol{L}^T\cdot\boldsymbol{M}+\boldsymbol{S}=\boldsymbol{0}
	\end{split}
\end{equation}
Please note that we use the symbol '$\cdot$' to denote a standard matrix multiplication.
This set of equations is transformed into a mixed form by eliminating the bending moments $\boldsymbol{M}=\left(
M_{x},M_{y},M_{xy}
\right)^{{T}}$ using the corresponding constitutional relation $\boldsymbol{M}=\hat{\boldsymbol{D}}_M\cdot\boldsymbol{L}\cdot\boldsymbol{\Theta}$. Furthermore, the constitutive equation for shear $\boldsymbol{S}=\hat{\boldsymbol{D}}_S\left(\boldsymbol{\nabla}\cdot w-\boldsymbol{\Theta}\right)$ is added in an inverted version as an additional equation. This procedure introduces the shear forces $\boldsymbol{S}=\left(
S_{1},S_{2}
\right)^{{T}}$ as independent parameters, additionally to the standard deformation unknowns (displacement $w$ and rotations $\boldsymbol{\Theta}=\left(
    \Theta_1,\Theta_2\right)^T$) of the irreducible form:
\begin{equation}
	\begin{split}
		\boldsymbol{\nabla}^T\cdot\boldsymbol{S}+f^{\textrm{ext}}=0\\
		\boldsymbol{L}^T\cdot\left(\hat{\boldsymbol{D}}_M\cdot\boldsymbol{L}\cdot\boldsymbol{\Theta}\right)+\boldsymbol{S}=\boldsymbol{0}\\
		\left(\boldsymbol{\nabla}\cdot w-\boldsymbol{\Theta}\right)-\hat{\boldsymbol{D}}_S^{-1}\cdot\boldsymbol{S}=\boldsymbol{0}
	\end{split}
\end{equation}
Therein, the material matrices $\hat{\boldsymbol{D}}_M$ and $\hat{\boldsymbol{D}}_S$ are defined based on the standard material parameters (Young's modulus $E$, Poisson's ratio $\nu$, shear correction factor $\kappa=\frac{5}{6}$ and shear modulus $G$) as well as the plate thickness $t$:
\begin{equation}
	\begin{split}
		\hat{\boldsymbol{D}}_M
		&=
		\frac{E\cdot t^3}{12\cdot\left(1-\nu^2\right)}\cdot
		\left[
		\begin{array}{ccc}
			1 & \nu & 0\\
			\nu & 1 & 0\\
			0 & 0 & \frac{1-\nu}{2}
		\end{array}
		\right]
		\\
		\hat{\boldsymbol{D}}_S
		&=
		\quad
		\kappa\cdot G\cdot t\cdot
		\left[
		\begin{array}{cc}
			1 & 0\\
			0 & 1
		\end{array}
		\right]
	\end{split}
\end{equation}
The corresponding weak form, defined within the plate domain $\Omega$, can be derived by introducing the variations of the unknowns and employing integration by parts if appropriate:
\begin{equation}
\label{eq:strongform}
	\begin{split}
		&\int_\Omega\left(\boldsymbol{\nabla}\cdot\delta w\right)^T\cdot\boldsymbol{S}~\text{d}A=\int_\Omega\delta w^T\cdot f^{\textrm{ext}}~\text{d}A\\&\int_\Omega\left(\boldsymbol{L}\cdot\delta\boldsymbol{\Theta}\right)^T\cdot\hat{\boldsymbol{D}}_M\cdot\left(\boldsymbol{L}\cdot\boldsymbol{\Theta}\right)-\delta\boldsymbol{\Theta}^T\cdot\boldsymbol{S}~\text{d}A={0}\\
		&\int_\Omega\delta \boldsymbol{S}^T\cdot\left(\left(\boldsymbol{\nabla}\cdot w-\boldsymbol{\Theta}\right)-\hat{\boldsymbol{D}}_S^{-1}\cdot\boldsymbol{S}\right)~\text{d}A={0}
	\end{split}
\end{equation}
The discretization is performed by independent interpolations of the unknown quantities
\begin{equation}
\label{eq:interpolation_general}
		w=\boldsymbol{R}_w\cdot\tilde{\boldsymbol{w}}
		\text{,}\quad
		\boldsymbol{\Theta}=\boldsymbol{R}_\Theta\cdot\tilde{\boldsymbol{\Theta}}
		\quad\text{and}\quad
		\boldsymbol{S}=\boldsymbol{R}_S\cdot\tilde{\boldsymbol{S}}=\begin{bmatrix}
		    \boldsymbol{R}_{S_1}&\boldsymbol{0}\\\boldsymbol{0}&\boldsymbol{R}_{S_2}
		\end{bmatrix}\cdot\tilde{\boldsymbol{S}}
		\text{,}
\end{equation}
where the matrices $\boldsymbol{R}_\chi$ contain the NURBS basis functions for each field, whereby each basis function is multiplied by the identity matrix of the appropriate dimension (cf.~Eq.~\eqref{eq:discretization}), and $\tilde{\boldsymbol{\chi}}$ is the associated vector of nodal values for each type of parameter ${\chi\in\left\{w, \Theta, S_1, S_2\right\}}$.
The associated vector of unknowns and the corresponding load-vector can be subdivided into their irreducible and mixed part, respectively:
\begin{alignat}{5}
\label{eq:unknowns}
		\tilde{\boldsymbol{v}}
		&=
		\left(
		\begin{array}{c}
			\tilde{\boldsymbol{d}}\\
			\tilde{\boldsymbol{S}}\\
		\end{array}
		\right)
		&&\quad\text{with}\quad
		\tilde{\boldsymbol{d}}
		&&=
		\left(
		\begin{array}{c}
			\tilde{\boldsymbol{w}}\\
			\tilde{\boldsymbol{\Theta}}\\
		\end{array}
		\right)
		&&\quad\text{and}\quad
		\tilde{\boldsymbol{S}}
		&&=
		\left(
		\begin{array}{c}
			\tilde{\boldsymbol{S}}_1\\
			\tilde{\boldsymbol{S}}_2\\
		\end{array}
		\right)
		\\
		\boldsymbol{f}
		&=
		\left(
		\begin{array}{c}
			\boldsymbol{f}_{d}\\
			\boldsymbol{0}\\
		\end{array}
		\right)
		&&\quad\text{with}\quad
		\boldsymbol{f}_{d}
		&&=
		\left(
		\begin{array}{c}
			\boldsymbol{f}_w\\
			\boldsymbol{0}
		\end{array}
		\right)
		&&\quad\text{and}\quad
		\boldsymbol{f}_{{w}}
		&&=
		\int_\Omega\boldsymbol{R}_w^T\cdot f^{\textrm{ext}}~\text{d}A
\end{alignat}
Thus, the system matrix for a linear mixed Reissner-Mindlin plate results as follows:
\begin{equation}
\label{eq:system_matrix_general}
	\boldsymbol{K}
	=
	\int_\Omega
	\left[
	\begin{array}{ccc}
		\boldsymbol{0}&
		\boldsymbol{0}&
		\left(\boldsymbol{\nabla}\cdot\boldsymbol{R}_w\right)^T\cdot\boldsymbol{R}_{S}
		\\
		\boldsymbol{0}&
		\left(\boldsymbol{L}\cdot\boldsymbol{R}_\Theta\right)^T\cdot\hat{\boldsymbol{D}}_M\cdot\left(\boldsymbol{L}\cdot\boldsymbol{R}_\Theta\right)&
		-\boldsymbol{R}_\Theta^T\cdot	\boldsymbol{R}_S
		\\
		 \boldsymbol{R}_S^T\cdot\boldsymbol{\nabla}\cdot\boldsymbol{R}_w&
		-\boldsymbol{R}_S^T\cdot	\boldsymbol{R}_\Theta&
		-\boldsymbol{R}_S^T\cdot\hat{\boldsymbol{D}}_S^{-1}\cdot\boldsymbol{R}_S
	\end{array}
	\right]
	~\text{d}A
\end{equation}
In \citep{StammenDornisch2023}, we examined that the interpolation of the shear forces should be adapted using selectively reduced approximation orders, while the irreducible part should not be modified, in order to counteract locking and obtain the optimal convergence rate. Thus, we redefine the discretizations of Eq.~\eqref{eq:interpolation_general} by allowing for different basis functions for the two shear forces $\boldsymbol{S}=\left(
S_{1},S_{2}
\right)^{{T}}$ and specify all employed orders for the basis functions:
\begin{equation}
\label{eq:discretization}
	\begin{aligned}
        {w}&=\boldsymbol{R}_{w}\cdot\tilde{\boldsymbol{w}}~\qquad&&\text{with }\qquad &&\boldsymbol{R}_w&=&\boldsymbol{R}_{w~}^{p,p}&=&\begin{bmatrix}
		    R^{p,p}_{ij},\dots,R^{p,p}_{nm}
		\end{bmatrix}\\
		\boldsymbol{\Theta}&=\boldsymbol{R}_{\Theta}\cdot\tilde{\boldsymbol{\Theta}}~
        \qquad&&\text{with }\qquad &&\boldsymbol{R}_\Theta&=&\boldsymbol{R}_\Theta^{p,p}&=&\begin{bmatrix}
		    R^{p,p}_{ij}\boldsymbol{I}_2,\dots,R^{p,p}_{nm}\boldsymbol{I}_2
		\end{bmatrix}\\
		{S}_1&=\boldsymbol{R}_{S_1}\cdot\tilde{\boldsymbol{S}}_1
        \qquad&&\text{with }\qquad &&\boldsymbol{R}_{S_1}&=&\boldsymbol{R}_{S_1}^{p-1,p}&=&\begin{bmatrix}
		    R^{p-1,p}_{ij},\dots,R^{p-1,p}_{n_S m}
		\end{bmatrix}\\
		{S}_2&=\boldsymbol{R}_{S_2}\cdot\tilde{\boldsymbol{S}}_2
        \qquad&&\text{with }\qquad &&\boldsymbol{R}_{S_2}&=&\boldsymbol{R}_{S_2}^{p,p-1}&=&\begin{bmatrix}
		    R^{p,p-1}_{ij},\dots,R^{p,p-1}_{n m_S}
		\end{bmatrix}
	\end{aligned}
\end{equation}
In this equation, $\boldsymbol{I}_2$ represents the identity matrix of dimension 2, as this is the number of rotation parameters that are interpolated using the same basis functions. Due to the required split in the discretization of the shear forces, the dimension of each shear force space is $1$, as for the displacement parameter. It should be noted that the number of nodes in the direction of the reduced order is lower, which is denoted by $n_S$ and $m_S$. The exact amount depends on the knot vectors employed. In a multi-patch setting, the deformations and rotations are chosen to be $C^0$-continuous over the patch intersections, while the shear forces are interpolated discontinuously between the individual patches.

\section{Efficient static condensation procedure}
\label{sec:procedure}
One of the main parts of this work is to propose an efficient static condensation procedure based on enhanced approximate dual basis functions. It is derived based on the system matrices defined in Sec.~\ref{sec:derivation}. However, the general procedure can be applied to any formulation of similar structure. Independently of our research, ~\cite{NguyenEtAl2024} proposed a static condensation approach for a thin shell formulation that is based on approximate dual basis functions. However, in contrast to those investigations, we include the proper treatment of general knot vectors, see Sec.~\ref{sec:internal_continuity}, which requires a reduction of continuity and employing enhanced approximate dual basis functions of NURBS, at the example of a Reissner--Mindlin plate formulation.

The first step is to aggregate the irreducible and the independent mixed parts of $\boldsymbol{K}$ according to Eq. \eqref{eq:unknowns}.
Since the shear constitutive matrix $\hat{\boldsymbol{D}}_S$ is diagonal and the shear forces $S_1$ and $S_2$ are interpolated independently as indicated in Eq.~\eqref{eq:interpolation_general}, the part of the system matrix coupling the virtual and physical shear parameters can be written as
\begin{equation}
\boldsymbol{K}_{S-S}=
    \begin{bmatrix}
        \begin{array}{cc}
    		\boldsymbol{K}_{S_1-S_1} & \boldsymbol{0} \\
    		  \boldsymbol{0}  & \boldsymbol{K}_{S_2-S_2}
    	\end{array}
    \end{bmatrix}\,.
    \label{eq:shear_part_total}
\end{equation}
Accordingly, we define
\begin{equation}
\boldsymbol{K}_{d-S}=
\left[
	\begin{array}{rr}
		\boldsymbol{K}_{d-S_1}  & \boldsymbol{K}_{d-S_2}
	\end{array}
	\right]\qquad\text{and}\qquad
\boldsymbol{K}_{S-d}=
	\left[
	\begin{array}{r}
		\boldsymbol{K}_{S_1-d}\\
		\boldsymbol{K}_{S_2-d}
	\end{array}
	\right]\,.
\end{equation}
Thus, the system matrix of Eq.~\eqref{eq:system_matrix_general} can be subdivided as follows:
\begin{equation}
	\boldsymbol{K}=
	\left[
	\begin{array}{ccc}
		\boldsymbol{K}_{d-d} & \boldsymbol{K}_{d-S_1} & \boldsymbol{K}_{d-S_2} \\
		\boldsymbol{K}_{S_1-d} & \boldsymbol{K}_{S_1-S_1} & \boldsymbol{0} \\
		\boldsymbol{K}_{S_2-d} & \boldsymbol{0}  & \boldsymbol{K}_{S_2-S_2}
	\end{array}
	\right]
    \label{eq:split_system}
\end{equation}
\begin{figure}[t]
	\centering
	\includegraphics[width=\textwidth]{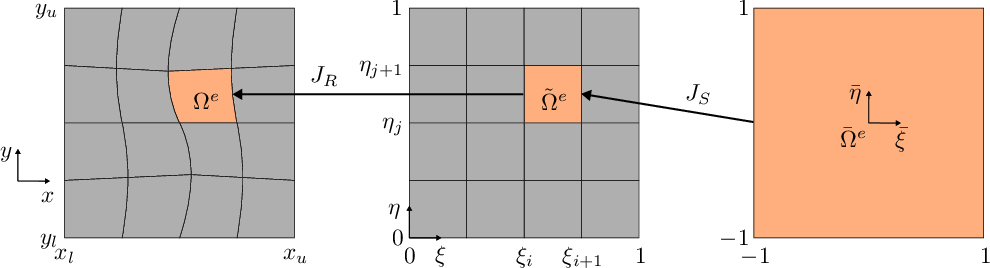}
	\caption{{Mapping between element in physical domain $\Omega^e$, element in parametric domain $\tilde{\Omega}^e$ and parent element $\bar{\Omega}^e$}}
	\label{fig:mapping}
\end{figure}
The resulting sub-matrices are computed element-wise and assembled globally/on patch-level later on. Due to the fact that integrations are performed for the parent element $\bar{\Omega}^e$, the Jacobian $J=J_S\cdot J_R$ is used for the equations related to the virtual displacements and rotations. As the approximate duality only holds for B-splines in the parametric domain, we employ a weighted integration with $\frac{W^2}{J_R}$ for the equations of Eq.~\eqref{eq:split_system} that are associated with the virtual (dual) shear parameters. Later on, this adapted weighting in conjunction with a pre-multiplication of the relevant system matrix parts with the transformation matrices of Sec.~\ref{sec:approximate_duals} ensures the diagonal dominance of the sub-matrices $\boldsymbol{K}_{S_\alpha-S_\alpha}$  for~$\alpha\in \left\{1,2\right\}$, required to minimize the error introduced by the proposed lumping procedure.
Due to the diagonality of $\hat{\boldsymbol{D}}_S$ and $\boldsymbol{K}_{S-S}$, only the factor $\frac{1}{\kappa\cdot G\cdot t}$ of $\hat{\boldsymbol{D}}_S^{-1}$ remains. Hence, we multiply all rows associated with the virtual (dual) shear parameters with ${\kappa\cdot G\cdot t}$, to achieve that the matrices $\boldsymbol{K}_{S_\alpha-S_\alpha}^e$, which only depend on the NURBS basis functions and the weighting function $w$ that have been defined in Eq.~\eqref{eq:NURBS_basis}. Finally, the sub-matrices are computed employing the equations

\begin{align}
\label{element_matrix_factor}
	\begin{split}
		\boldsymbol{K}_{d-d}^e
		&=
		\int_{\bar{\Omega}^e}
		\left[
		\begin{array}{cc}
			\boldsymbol{0}&
			\boldsymbol{0}
			\\
			\boldsymbol{0}&
			\left(\boldsymbol{L}\cdot\boldsymbol{R}_\Theta\right)^T\cdot\hat{\boldsymbol{D}}_M\cdot\left(\boldsymbol{L}\cdot\boldsymbol{R}_\Theta\right)
		\end{array}
		\right]
		\cdot J
		~\text{d}A
	\\
		\boldsymbol{K}_{d-S_\alpha}^e
		&=
		\int_{\bar{\Omega}^e}
		\left(
		\begin{array}{c}
			\left(\boldsymbol{\nabla}\cdot\boldsymbol{R}_w\right)^T\cdot\boldsymbol{R}_{S_\alpha}
			\\
			-\boldsymbol{R}_\Theta^T\cdot\boldsymbol{R}_{S_\alpha}
		\end{array}
		\right)
		\cdot J
		~\text{d}A
	\\
		\boldsymbol{K}_{S_\alpha-d}^e
		&=
		\int_{\bar{\Omega}^e}
		\left(
		\begin{array}{cc}
			\boldsymbol{R}_{S_\alpha}^T\cdot\boldsymbol{\nabla}\cdot\boldsymbol{R}_w
			-\boldsymbol{R}_{S_\alpha}^T\cdot\boldsymbol{R}_\Theta
		\end{array}
		\right)
		\cdot\kappa\cdot G\cdot t\cdot W^2\cdot J_S
		~\text{d}A
	\\
		\boldsymbol{K}_{S_\alpha-S_\alpha}^e
		&=
		\int_{\bar{\Omega}^e}
		-\boldsymbol{R}_{S_\alpha}^T\cdot\boldsymbol{R}_{S_\alpha}
		\cdot W^2 \cdot J_S
		~\text{d}A
		\end{split}
\end{align}
for~$\alpha\in \left\{1,2\right\}$.
Therein, the Jacobian determinant $J_S$ for the mapping from parent element to parameter space is calculated by the partial derivative of the parametric coordinates ${\boldsymbol{\xi}=\left(
\xi,\eta\right)^T}$ with respect to the parent element coordinates $\bar{\boldsymbol{\xi}}=\left(\begin{array}{cc}
\bar{\xi},\bar{\eta}\end{array}\right)^T$ and that of the mapping from parameter space to physical space $J_R$ is determined by the partial derivative of the physical coordinates $\boldsymbol{x}=\left(
x,y\right)^T$ with respect to the parametric coordinates $\boldsymbol{\xi}$, i.e.,
\begin{equation}
	J_S=\frac{\partial \boldsymbol{\xi}}{\partial \bar{\boldsymbol{\xi}}}
	\quad\text{and}\quad
	J_R=\frac{\partial \boldsymbol{x}}{\partial \boldsymbol{\xi}}\,.
\end{equation}
The parts of the system matrix are assembled from the corresponding element matrices according to the assembly of the domain $\Omega=\overset{e}{\cup}\Omega_e$ from its elements $\Omega_e$:
\begin{equation}
\label{eq:elementwiseAssembly}
	\boldsymbol{K}_{X}=\overset{e}{\cup}\boldsymbol{K}_{X}^e
	\text{,}\quad
	X\in\left\{d-d, d-S_\alpha, S_\alpha-d, S_\alpha-S_\alpha\right\}
	\text{,}\quad
    \text{for}~\alpha\in \left\{1,2\right\}
\end{equation}

\begin{figure}[t]
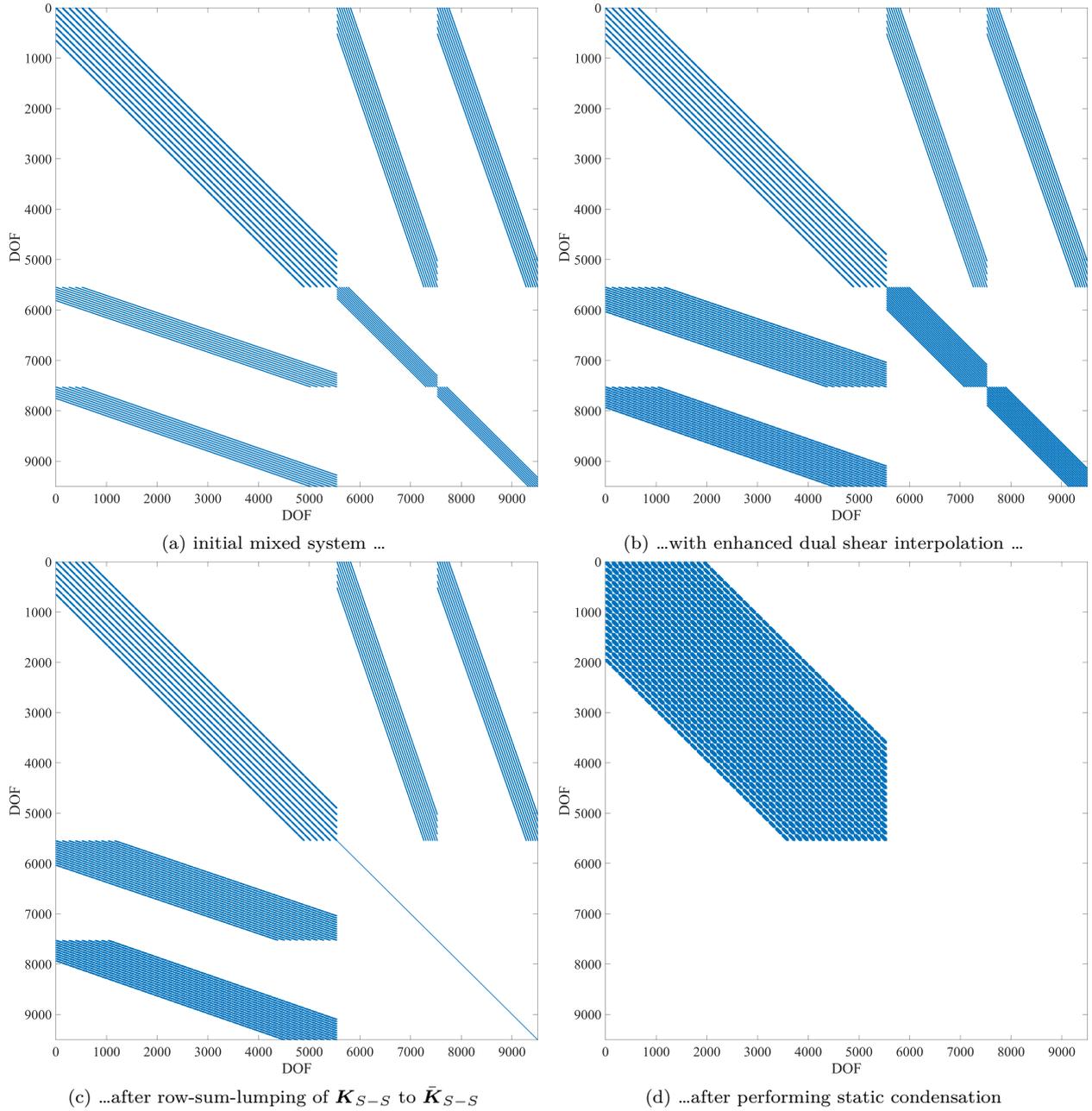

	\centering
	\begin{subfigure}{0.49\textwidth}
		\includegraphics[width=\textwidth]{figures/mixed.png}
		\caption{initial mixed system \dots}
	\end{subfigure}
	\hfill
	\begin{subfigure}{0.49\textwidth}
		\includegraphics[width=\textwidth]{figures/duals.png}
		\caption{\dots with enhanced dual shear interpolation \dots}
	\end{subfigure}
	\\
	\begin{subfigure}{0.49\textwidth}
		\includegraphics[width=\textwidth]{figures/lumped.png}
		\caption{\dots after row-sum-lumping of $\boldsymbol{K}_{S-S}$ to $\bar{\boldsymbol{K}}_{S-S}$}
	\end{subfigure}
	\hfill
	\begin{subfigure}{0.49\textwidth}
		\includegraphics[width=\textwidth]{figures/cond.png}
		\caption{\dots after performing static condensation}
	\end{subfigure}
	\caption{Sparsity of mixed system matrix for benchmark example \nameref{sec:numerical_exmples} examined in Sec.~\ref{sec:quadraticPlate_ud} in case $t=0.1~m$, $r=p=5$ and $40\times40$ elements for all sub-steps of the proposed efficient static condensation procedure.}
	\label{fig:procedure}
\end{figure}

The second step is to eliminate the additionally introduced shear parameters $\tilde{\boldsymbol{S}}$ from the mixed formulation, which reduces the computational effort of the analysis but preserves its benefits with regard to accuracy and alleviation of locking phenomena. This procedure is commonly known as static condensation:
\begin{equation}
	\left(\boldsymbol{K}_{d-d}-\boldsymbol{K}_{d-S}\cdot\boldsymbol{K}_{S-S}^{-1}\cdot\boldsymbol{K}_{S-d}\right)\cdot\tilde{\boldsymbol{d}}
	=
	\boldsymbol{f}_d
\end{equation}
Due to the fact that the inverse of the shear part
$\boldsymbol{K}_{S-S}$ (see Eq.~\eqref{eq:shear_part_total})
of the system matrix has to be determined, static condensation can get computationally expensive, since $\boldsymbol{K}_{S-S}$ consists of two banded matrices (cf. Fig.~\ref{fig:procedure}~(a)), where the bandwidth depends on the order~$p$.
In order to reduce the computational costs associated with the static condensation to a minimum, lumping the shear part $\boldsymbol{K}_{S-S}$ to a diagonal matrix, whose inverse can be determined directly, is intended. As this potentially causes a huge error, depending on the initial bandwidth of this matrix part and its diagonal dominance prior to lumping, employing approximate dual basis functions for the interpolation of the test functions for the shear parameters is considered. Due to their nearly-bi-orthogonality, this ensures diagonal dominant sub-matrices. This minimizes the error introduced by the proposed lumping procedure.

Thus, the third step is to employ (enhanced) approximate dual basis functions for the virtual shear parameters, which is equivalent to pre-multiplying the relevant system sub-matrices with the corresponding transformation matrices
\begin{equation}
\label{eq:transformationMatrices}
   \begin{aligned}
    \boldsymbol{T}_1=
    \left\{
    \begin{array}{rr}
    {}^R\boldsymbol{S}_{r_1,s_1}^{p-1,p}\text{,} & \text{approximate duals}\\
    {}^{enh,R}\boldsymbol{S}_{r_1,s_1}^{p-1,p}\text{,} & \text{enhanced approximate duals}
    \end{array}
    \right.
    \qquad&\text{with}\qquad  r_1\leq p-1,~ s_1\leq p\\
     \boldsymbol{T}_2=
     \left\{
    \begin{array}{rr}
    {}^R\boldsymbol{S}_{r_2,s_2}^{p,p-1}\text{,} & \text{approximate duals}\\
    {}^{enh,R}\boldsymbol{S}_{r_2,s_2}^{p,p-1}\text{,} & \text{enhanced approximate duals}
    \end{array}
    \right.
    \qquad&\text{with}\qquad  r_2\leq p,~ s_2\leq p-1 \,,
     \end{aligned}
\end{equation}
which are defined in Eq.~\eqref{eq:2d_NURBS_transformation_matrix} or Eq.~\eqref{eq:2d_NURBS_transformation_matrix_enh}, respectively, depending on whether approximate dual basis functions or enhanced approximate dual basis functions are employed within the proposed procedure. This
yields the adapted system matrix
\begin{equation}
	\label{eq:dualSystemBsplines}
	{}^{PG}\boldsymbol{K}=
	\left[
	\begin{array}{rrr}
		\boldsymbol{K}_{d-d} & \boldsymbol{K}_{d-S_1}  & \boldsymbol{K}_{d-S_2} \\
		\boldsymbol{T}_1\cdot\boldsymbol{K}_{S_1-d} & \boldsymbol{T}_1\cdot\boldsymbol{K}_{S_1-S_1} & \boldsymbol{0}\\		\boldsymbol{T}_2\cdot\boldsymbol{K}_{S_2-d} & \boldsymbol{0} & \boldsymbol{T}_2\cdot\boldsymbol{K}_{S_2-S_2}
	\end{array}
	\right]\,,
\end{equation}
where the index of ${}^{PG}\boldsymbol{K}$ denotes that the formulation becomes of Petrov-Galerkin type, since the employed test and trial functions for the shear parameters are not identical. This also leads to a non-symmetric system matrix. However, it is to be noted that the transformations applied in Eq.~\eqref{eq:dualSystemBsplines} do not change the underlying function space. Thus, the eigenvalues of the system matrix remain positive.
Due to the transformation, both matrix parts that couple the respective virtual and physical shear parameters
\begin{equation}
    {}^{PG}\boldsymbol{K}_{S_\alpha-S_\alpha}=\boldsymbol{T}_\alpha\cdot\boldsymbol{K}_{S_\alpha-S_\alpha}
\end{equation}
that occur in Eq.~\eqref{eq:dualSystemBsplines}, and need to be inverted within static condensation, possess a slightly increased bandwidth compared to the initial matrix part $\boldsymbol{K}_{S_\alpha-S_\alpha}$~(cf.~Fig.~\ref{fig:procedure}~(a) and (b)), but are diagonal dominant. Both the increase of bandwidth and the diagonal dominance depend of the chosen degree of reproduction. In Fig.~\ref{fig:procedure} and throughout the numerical examples, full reproduction is employed.

The fourth step is to apply row-sum lumping to both matrix parts ${}^{PG}\boldsymbol{K}_{S_\alpha-S_\alpha}$, which couple the virtual and physical shear parameters. This operation can be expressed by
\begin{equation}
{}^{PG}\tilde{\boldsymbol{K}}_{S_\alpha-S_\alpha}=\operatorname{diag}({}^{PG}\boldsymbol{K}_{S_\alpha-S_\alpha}\cdot\boldsymbol{1})
\,,
\end{equation}
where $\boldsymbol{1}$ is a vector of ones with adequate length.
Since the matrices ${}^{PG}\boldsymbol{K}_{S_\alpha-S_\alpha}$ are diagonal dominant, the error induced by the lumping is very small.
\begin{remark}
    For the similar case of mortar-based patch coupling, it was shown in~\cite{DornischEtAl2017} that the error decay rate of the lumping error is mostly the same as for the solution approximation error. Considering the additional measures proposed in~\cite{DornischStoeckler2021} and~\cite{Stoeckler2025}, which are explained in Sec.~\ref{sec:internal_continuity} in detail, the lumping error always possesses the same decay rate as the solution approximation error. Proving this mathematically is ongoing work within the framework of mortar methods. Within the present paper, this fact is only shown numerically.
\end{remark}
Due to the underlying integration space and the additional factors introduced in Eq.~\eqref{element_matrix_factor} for the shear part of the element matrix, according to~\cite{ChuiEtAl2004}, row-sum lumping of ${}^{PG}\boldsymbol{K}_{S_\alpha-S_\alpha}$ always yields the identity matrix~$\boldsymbol{I}$, i.e., ${}^{PG}\tilde{\boldsymbol{K}}_{S_\alpha-S_\alpha}=\boldsymbol{I}$. This yields the partially-lumped system matrix
\begin{equation}
	\label{eq:dualSystemBsplineslumped}
	{}^{PG}\tilde{\boldsymbol{K}}=
	\left[
	\begin{array}{rrr}
		\boldsymbol{K}_{d-d} & \boldsymbol{K}_{d-S_1}  & \boldsymbol{K}_{d-S_2} \\
		\boldsymbol{T}_1\cdot\boldsymbol{K}_{S_1-d} & \boldsymbol{I} & \boldsymbol{0}\\
		\boldsymbol{T}_2\cdot\boldsymbol{K}_{S_2-d} & \boldsymbol{0} & \boldsymbol{I}
	\end{array}
	\right]\,.
\end{equation}
The static condensation of the resulting system of equations ${}^{PG}\tilde{\boldsymbol{K}}\tilde{\boldsymbol{v}}=\boldsymbol{f}$ is now trivial, since no matrix inversions are required. This enables a significant speedup in comparison to the standard approach of static condensation. The final system to be solved is given by
\begin{equation}
\label{eq:condensedSystem}
	\left(\boldsymbol{K}_{d-d}-\sum_{\alpha=1}^2{\boldsymbol{K}_{d-S_\alpha}\cdot\boldsymbol{T}_\alpha\cdot\boldsymbol{K}_{S_\alpha-d}}\right)\cdot\tilde{\boldsymbol{d}}
	=
	\boldsymbol{f}_d\,.
\end{equation}
This results in a banded system matrix with increased bandwidth, but significantly less degrees of freedom. This is depicted exemplarily in Fig.~\ref{fig:procedure}~(d).

Due to the inter-element continuity of the employed spline functions, static condensation cannot be performed on element-level. Hence, for single-patch problems, it needs to be conducted globally. For multi-patch geometries, static condensation can be performed either globally or on patch-level, since the shear forces are interpolated discontinuously in-between patches~(cf.~Sec.~\ref{sec:derivation}). For the sake of efficiency, it is conducted on patch-level for the proposed formulation. The transformation required for approximate dual basis functions or enhanced approximate dual basis functions can be performed at any level. To gain a reasonable efficiency while facilitating the implementation, it is conducted globally for single-patch geometries as indicated in Eq.~\eqref{eq:condensedSystem} and performed on element-level if the geometry to be analyzed consists of several patches. Following this approach, the deformation parts of the global system matrix $\boldsymbol{K}_{d-d}$ and $\boldsymbol{K}_{d-S_\alpha}$ for $\alpha\in\{1,2\}$ can still be assembled as indicated in Eq.~\eqref{eq:elementwiseAssembly}. The transformation of the global system matrix parts ${}^{PG}\boldsymbol{K}_{S_\alpha-d}$ would require to assemble a global transformation matrix from the patch-wise transformation matrices $\boldsymbol{T}_\alpha$ for $\alpha\in\{1,2\}$ in a suitable manner. Since this would be a huge block-diagonal matrix, it is anticipated to result in an increased overall computational effort compared to the element-wise transformation pursued in the present approach. 
Hence, the global or patch-wise pre-multiplication of the transformation matrices $\boldsymbol{T}_\alpha$ is replaced by the element-wise assembly $^{PG}\boldsymbol{K}_{S_\alpha-d}=\overset{e}{\cup}\,^{PG}\boldsymbol{K}_{S_\alpha-d}^e$
of the already transformed element matrices~$^{PG}\boldsymbol{K}_{S_\alpha-d}^e=\boldsymbol{T}^e_\alpha\cdot \boldsymbol{K}_{S_\alpha-d}^e$ for $\alpha\in \left\{1,2\right\}$.
Therefore, the local transformation matrices $\boldsymbol{T}_\alpha^e$ have to be computed by extracting the relevant parts of the patch-wise transformation matrices $\boldsymbol{T}_\alpha$ of Eq.~\eqref{eq:transformationMatrices} for $\alpha\in\{1,2\}$. The independent shear part requires no further special attention, because it becomes the identity matrix as derived in Eq.~\eqref{eq:dualSystemBsplineslumped} if (enhanced) approximate dual basis functions are employed following the procedure proposed in this contribution.\\
\begin{remark}
    The proposed formulation provides a unified methodology for both NURBS and B-splines. In the case of B-splines, the weighting occurring within the element matrices in Eq.~\eqref{element_matrix_factor} simplifies to $W=1$, the diagonal matrices~$\boldsymbol{\Lambda}$ in Eq.~\eqref{eq:2d_NURBS_transformation_matrix} become the identity matrix and the NURBS basis functions in Eq.~\eqref{eq:discretization} can be replaced by B-spline basis functions. Besides the costs for the computation of the basis functions itself, there is no significant difference in computational costs for the static condensation procedure by considering the more general case of NURBS in comparison to the limited case of B-splines as used in~\cite{NguyenEtAl2024}.
\end{remark}

The whole procedure is visualized in
Fig.~\ref{fig:procedure}, where the resulting changes in the matrix structure are depicted for each step, based on the numerical example examined in Sec.~\ref{sec:quadraticPlate_ud}. In these plots, the order of degrees of freedom first comprises the deflection and rotations per control points and subsequently all shear force parameters separated per direction. Fig.~\ref{fig:procedure}~(a) depicts the initial structure of the mixed system matrix. In Fig.~\ref{fig:procedure}~(b), the increase of bandwidth of the lower sub-matrices due to the pre-multiplication of the equations associated with the virtual shear parameters with the transformation matrices $\boldsymbol{T}_\alpha$ is visualized for the case of full reproduction. In comparison, employing dual shape functions of a lower reproduction degree, the enlargement of bandwidth is smaller. The effect of row-sum-lumping of the shear parts ${}^{PG}\boldsymbol{K}_{S_\alpha-S_\alpha}$ of the system matrix to ${}^{PG}\tilde{\boldsymbol{K}}_{S_\alpha-S_\alpha}$ is shown in Fig.~\ref{fig:procedure}~(c). Fig.~\ref{fig:procedure}~(d) depicts the final matrix structure after performing the whole proposed static condensation procedure.

\section{Numerical Examples}
\label{sec:numerical_exmples}

In this section, the robustness of the proposed dual lumping procedure is examined with respect to discretizations that include different types of limited internal continuities.
In particular, the results of the following formulations are compared:
\begin{itemize}
    \item standard purely displacement-based isogeometric formulation (std.)
    \item isogeometric displacement-stress mixed formulation without any lumping or condensation procedure (mxd.)
    \item entirely NURBS-based isogeometric displacement-stress mixed formulation with row-sum-lumping and subsequent static condensation (lmp.)
    \item displacement-stress mixed formulation with row-sum-lumping and subsequent static condensation employing approximate dual basis functions for the interpolation of the virtual stress parameters~(AD).
    \item displacement-stress mixed formulation with row-sum-lumping and subsequent static condensation employing enhanced approximate dual basis functions for the interpolation of the virtual stress parameters (eAD)
\end{itemize}
Additionally, the necessity of performing knot insertion as well as the possibility of omitting the weights of the NURBS basis functions employed for the discretization of the mixed fields is studied. Throughout these entire investigations, full reproduction is employed for both types of dual basis functions.

\begin{figure}[t]
	\centering
	\begin{subfigure}{0.24\textwidth}
		\includegraphics[width=\textwidth]{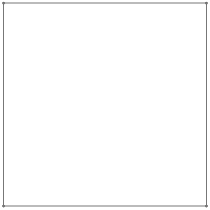}
		\caption{Undistorted geometry}
        \label{fig:example:sp:undistorted}
	\end{subfigure}
	\hfill
	\begin{subfigure}{0.24\textwidth}
		\includegraphics[width=\textwidth]{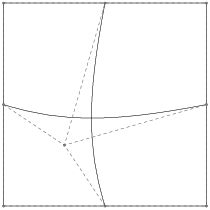}
		\caption{NURBS geometry}
                \label{fig:example:sp:NURBS}
	\end{subfigure}
	\hfill
	\begin{subfigure}{0.24\textwidth}
		\includegraphics[width=\textwidth]{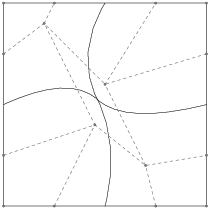}
		\caption{$C^1$- geometry}
                \label{fig:example:sp:C1}
	\end{subfigure}
	\hfill
	\begin{subfigure}{0.24\textwidth}
		\includegraphics[width=\textwidth]{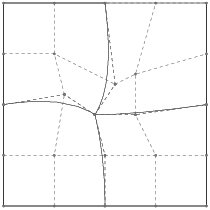}
		\caption{$C^0$- geometry}
                \label{fig:example:sp:C0}
	\end{subfigure}
	\caption{Initial meshes with isolines (solid), control points (circles) and control polygon (dashed)}
	\label{fig:msh}
\end{figure}

\begin{figure}[t]
	\centering
	\begin{subfigure}{0.24\textwidth}
		\includegraphics[width=\textwidth]{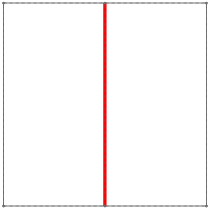}
		\caption{Multi-patch-geometry with linear interface}
                \label{fig:example:mp:linear}
	\end{subfigure}
	\hfill
	\begin{subfigure}{0.24\textwidth}
		\includegraphics[width=\textwidth]{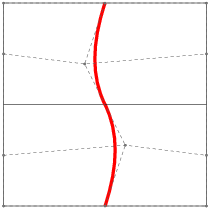}
		\caption{Multi-patch-geometry with $C^1$-continuity}
        \label{fig:example:mp:C1}
	\end{subfigure}
	\hfill
	\begin{subfigure}{0.24\textwidth}
		\includegraphics[width=\textwidth]{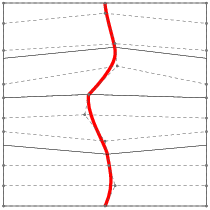}
		\caption{Multi-patch-geometry with various continuities}
        \label{fig:example:mp:Cvar}
	\end{subfigure}
	\caption{Initial meshes with patch interface (bold red solid), isolines (solid), control points (circles) and control polygon (dashed)}
	\label{fig:msh_mp}
\end{figure}

Due to the availability of a corresponding analytical solution, which is provided in Eq.~\eqref{eq:deformation_solution}, a fully clamped quadratic plate with complex loading is employed as numerical example. As the inherent complexity of the resulting displacement field allows to conduct meaningful error convergence studies, it is widely used for the verification of plate elements~\cite{ChinosiLovadina1995,BeiraoDaVeigaEtAl2012a}. The aim of this contribution is to study the influence of discretizations with limited internal continuity and distortions. Furthermore, we also study multi-patch settings, which yield discontinuity in the shear parameters. Thus, we use a multitude of discretization types but always use the same mechanical example. Plots of the coarsest mesh for the single-patch discretizations are depicted in Fig.~\ref{fig:msh}. The multi-patch discretizations are plotted in Fig.~\ref{fig:msh_mp}.
The initial surface orders, knot vectors, and control points of the examined geometries are included in \ref{apndx}. Furthermore, in \ref{supp_data} we provide the corresponding IGES files as supplementary material to this article, to allow direct access to the exact geometry representations that are employed for the analysis of the underlying benchmark problem.

\begin{figure}[t]
	\centering
	\begin{subfigure}{0.32\textwidth}
		\includegraphics[width=\textwidth]{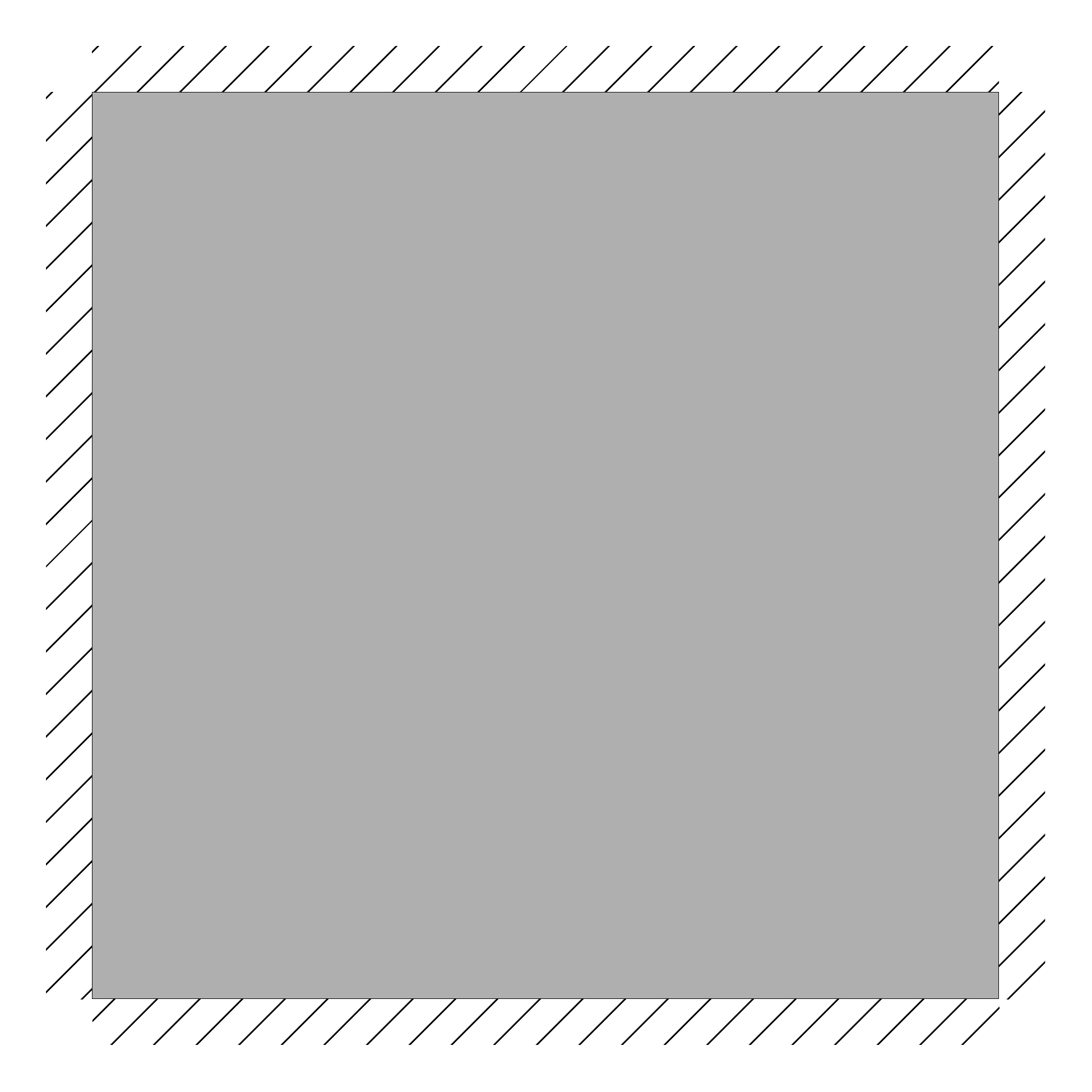}
		\caption{System geometry}
	\end{subfigure}
	\hfill
	\begin{subfigure}{0.32\textwidth}
		\includegraphics[width=\textwidth]{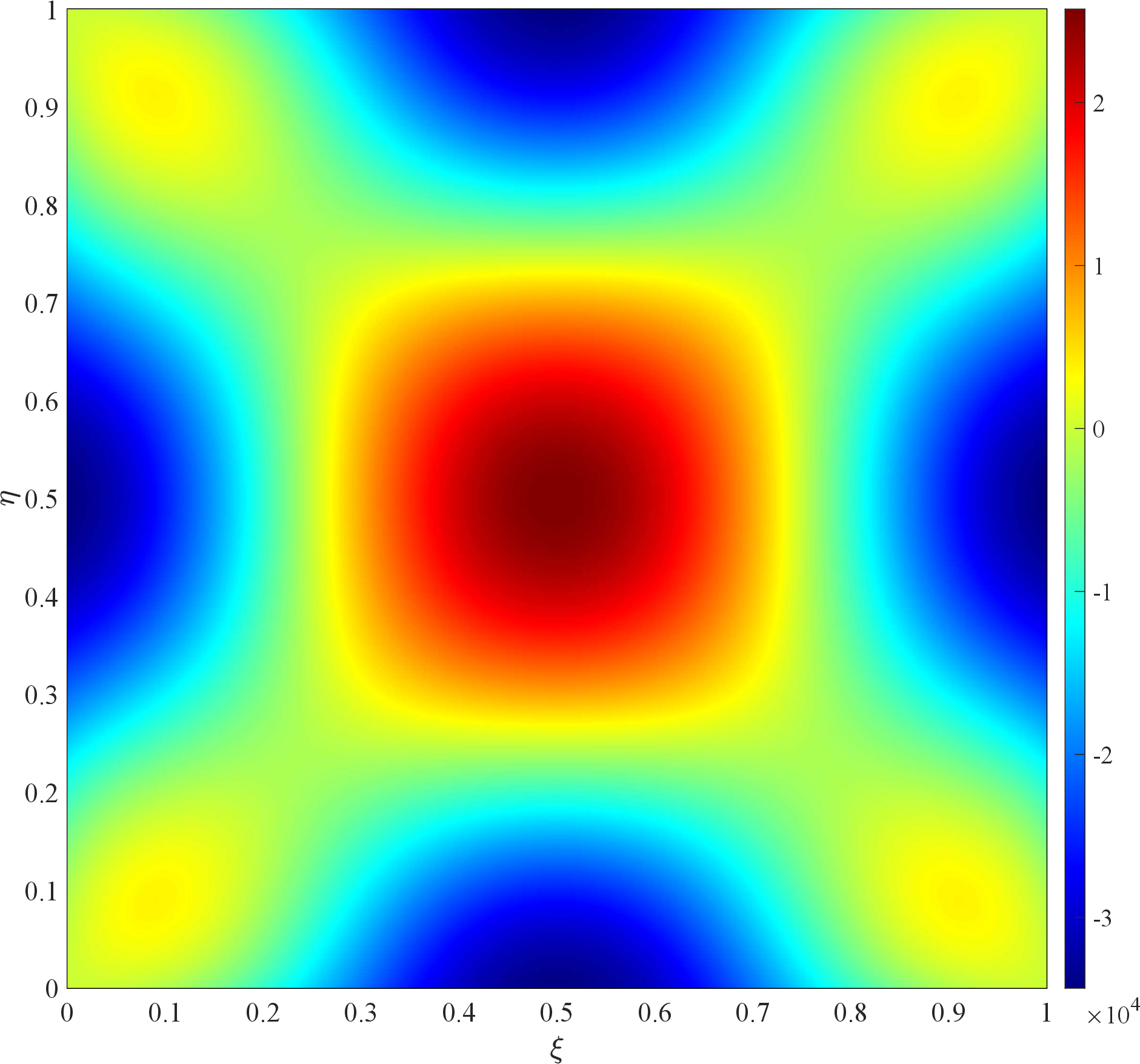}
		\caption{Load function}
	\end{subfigure}
	\hfill
	\begin{subfigure}{0.32\textwidth}
		\includegraphics[width=\textwidth]{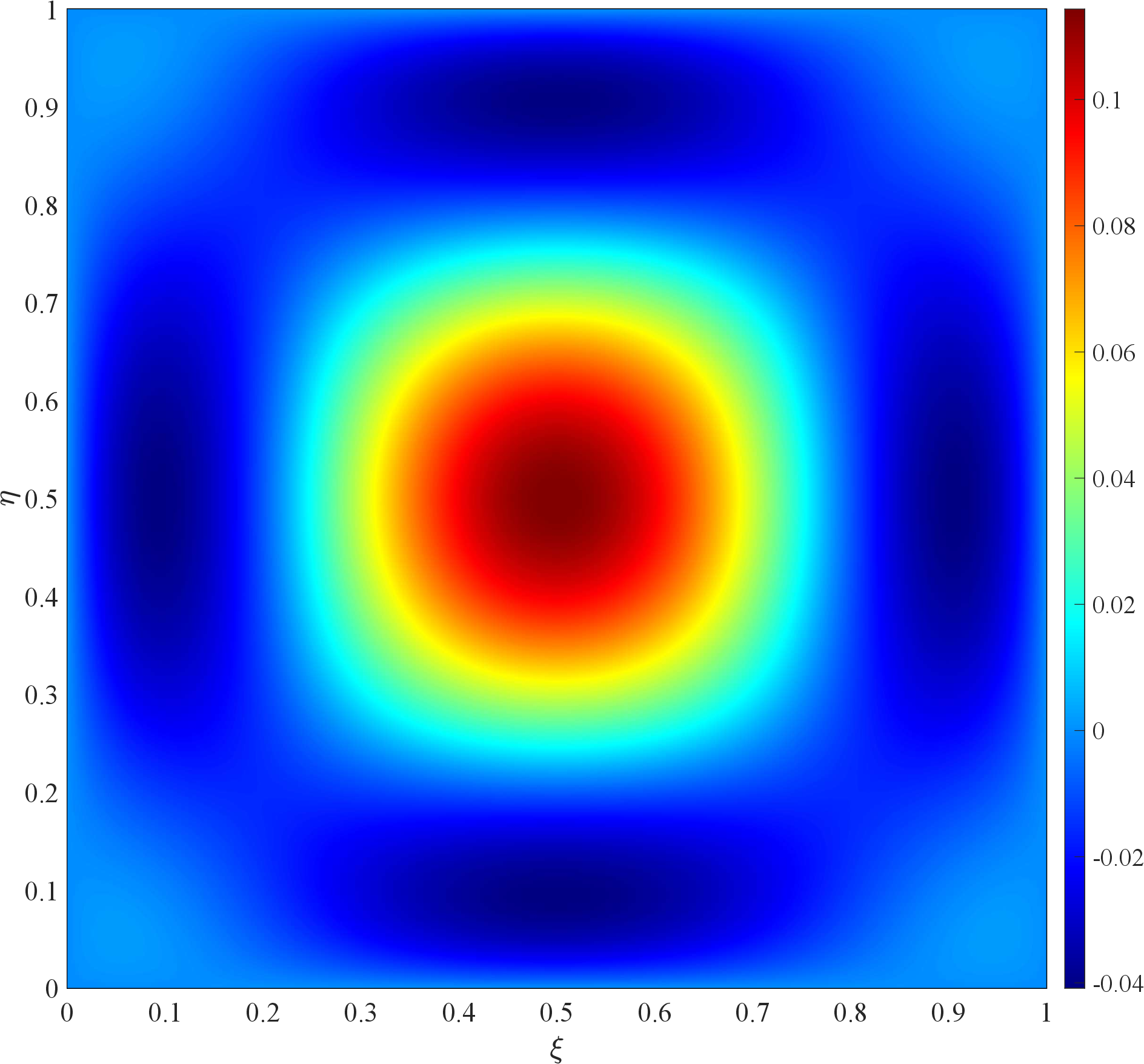}
		\caption{Reference solution of displacement}
	\end{subfigure}
	\caption{System geometry, load function and reference solution of displacement for the investigated benchmark problem}
	\label{fig:q+w}
\end{figure}

The transversal loading was adapted from \citep{ChinosiLovadina1995} as follows:
\begin{equation}
	\begin{split}
	f^{\textrm{ext}}\left(\xi\text{,}\eta\right)&=f_0^{\textrm{ext}}\cdot\frac{E\cdot\left(\frac{t}{1~\text{m}}\right)^3}{12\cdot(1-\nu^2)}\cdot\bigl(f^{\textrm{ext}}_1\left(\xi\text{,}~\eta\right)+f^{\textrm{ext}}_2\left(\xi\text{,}~\eta\right)\bigr)
	\text{,}\quad\text{where} \\
	f_0^{\textrm{ext}}&=100
	\text{,} \\
	f^{\textrm{ext}}_1\left(\xi\text{,}~\eta\right)&=12\cdot \hat{f}_2\left(\xi\text{,}~\eta\right)\cdot\bigl(2\cdot\eta^2\cdot\left(\eta-1\right)^2+\hat{f}_1\left(\xi\text{,}~\eta\right)\bigr)
	\quad\text{and} \\
	f^{\textrm{ext}}_2\left(\xi\text{,}~\eta\right)&=12\cdot \hat{f}_1\left(\xi\text{,}~\eta\right)\cdot\bigl(2\cdot\xi^2\cdot\left(\xi-1\right)^2+\hat{f}_2\left(\xi\text{,}~\eta\right)\bigr)
	\text{,}\quad\text{where}\\
	\hat{f}_1\left(\xi\text{,}~\eta\right)&=\xi\cdot\left(\xi-1\right)\cdot\left(5\cdot\eta^2-5\cdot\eta+1\right)
	\quad\text{and} \\
	\hat{f}_2\left(\xi\text{,}~\eta\right)&=\eta\cdot\left(\eta-1\right)\cdot\left(5\cdot\xi^2-5\cdot\xi+1\right)
	\end{split}
\end{equation}
Therein, the dimensionless length-parameters $\xi_1=x/L$ and $\xi_2=y/L$ are employed.
The corresponding reference solution for the deformation over the whole plate domain is adapted from \citep{ChinosiLovadina1995} as
\begin{equation}
\label{eq:deformation_solution}
	\begin{split}
	w^{\textrm{ex}}\left(\xi\text{,}~\eta\right)&=w_0\left(\xi\text{,}~\eta\right)-\frac{2\cdot t^2}{5\cdot(1-\nu)}\cdot\bigl(w_1\left(\xi\text{,}~\eta\right)+w_2\left(\xi\text{,}~\eta\right)\bigr)~\text{m}
	\text{,}\quad\text{where}\\
	w_0\left(\xi\text{,}~\eta\right)&=\frac{1}{3}\cdot\xi^3\cdot\left(\xi-1\right)^3\cdot\eta^3\cdot\left(\eta-1\right)^3
	\text{,}\\
	w_1\left(\xi\text{,}~\eta\right)&=\eta^2\cdot\left(\eta-1\right)^2\cdot\xi\cdot\left(\xi-1\right)\cdot \hat{f}_2\left(\xi\text{,}~\eta\right)
	\quad\text{and}\\
	w_2\left(\xi\text{,}~\eta\right)&=\xi^2\cdot\left(\xi-1\right)^2\cdot\eta\cdot\left(\eta-1\right)\cdot \hat{f}_1\left(\xi\text{,}~\eta\right)
	\quad\text{.}
	\end{split}
\end{equation}
The system geometry, load function as well as the displacement reference solution are depicted in Fig.~\ref{fig:q+w}.
The results for the displacement $w$ computed by the previously mentioned investigated formulations
are compared by means of the global $L_2$- error norm $||w||_{L_2}=\sqrt{\int_\Omega (w-w^{\textrm{ex}})^2\,\mathrm{d}A}$ to the exact solution $w^{\textrm{ex}}$ for $L=1~\text{m}$, $E=10000~\frac{\text{kN}}{\text{m}^2}$, $\nu=0.3$ and $t\in \left\{1\text{,}~0.1\text{,}~0.01\text{,}~0.001\text{,}~0.0001\right\}~\text{m}$, for both single-patch and multi-patch meshes involving various continuities. The consideration of a broad range of slenderness ratios and discretization types allows a proper judgement of the ability of the proposed method to alleviate transverse shear locking for complicated patch layouts, which frequently occur in industrial applications. However, we exclude the study of trimmed NURBS patches, which will be part of future studies.


\subsection{Undistorted mesh}
\label{sec:quadraticPlate_ud}
\FloatBarrier
\begin{figure}[thp]
	\centering
	\begin{subfigure}{0.49\textwidth}
		\includegraphics[width=\textwidth]{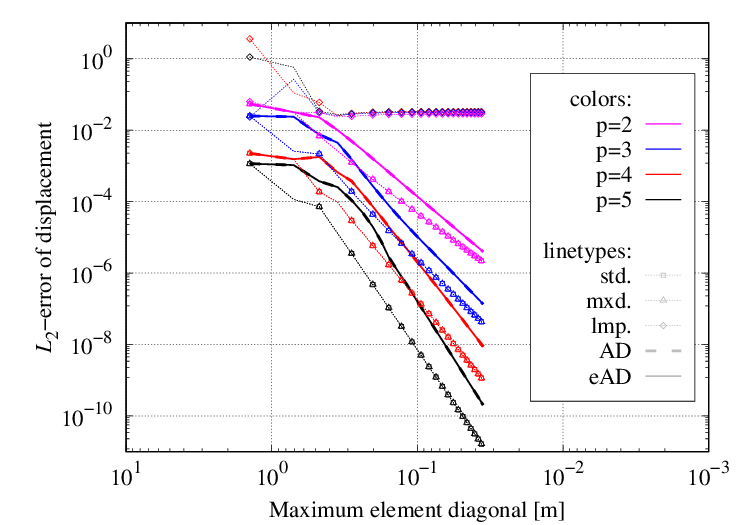}
		\caption{$t=1$~m}
	\end{subfigure}
	\hfill
	\begin{subfigure}{0.49\textwidth}
		\includegraphics[width=\textwidth]{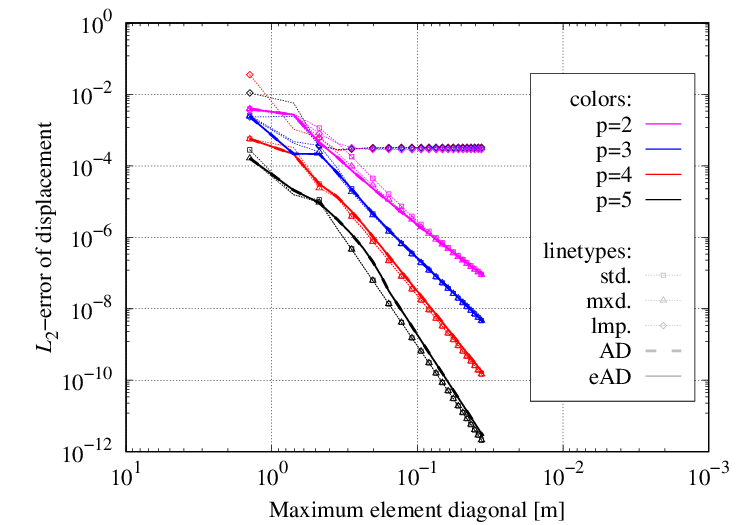}
		\caption{$t=0.1$~m}
	\end{subfigure}
	\\
	\begin{subfigure}{0.49\textwidth}
		\includegraphics[width=\textwidth]{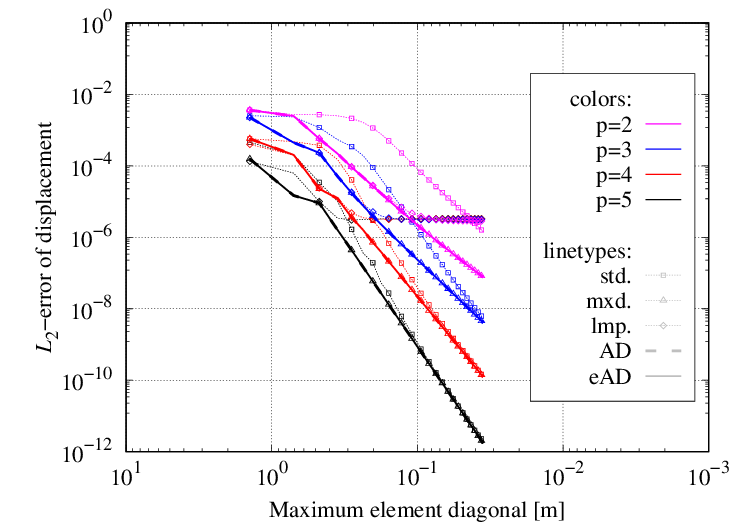}
		\caption{$t=0.01$~m}
	\end{subfigure}
	\hfill
	\begin{subfigure}{0.49\textwidth}
		\includegraphics[width=\textwidth]{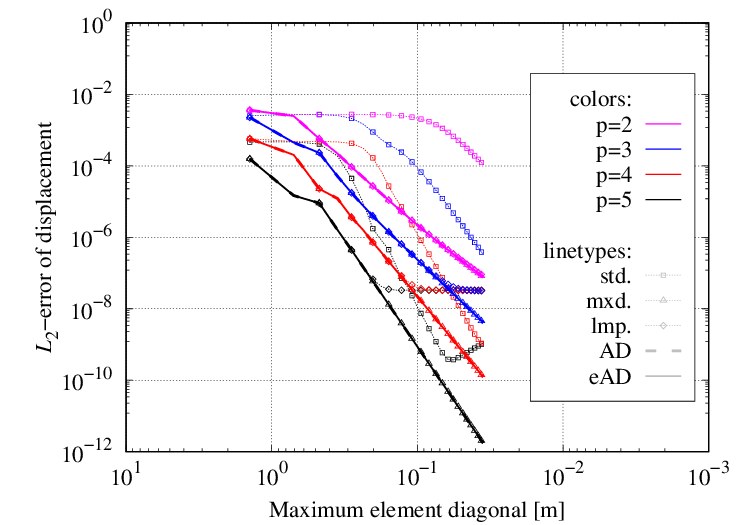}
		\caption{$t=0.001$~m}
	\end{subfigure}
	\\
	\begin{subfigure}{0.49\textwidth}
		\includegraphics[width=\textwidth]{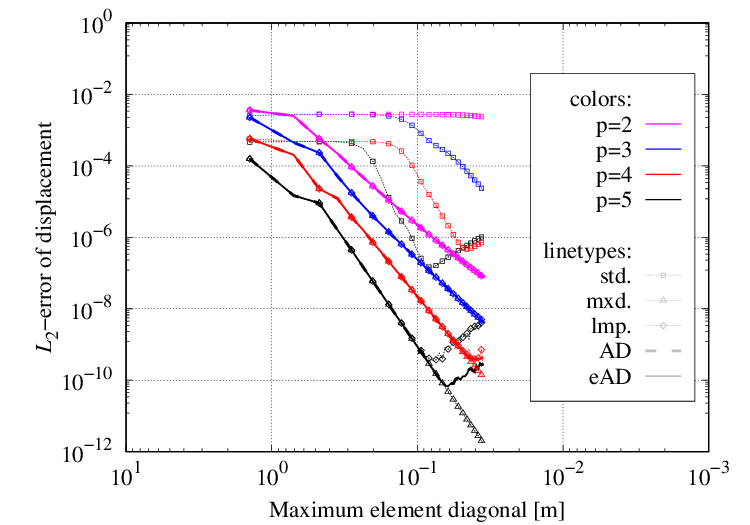}
		\caption{$t=0.0001$~m}
	\end{subfigure}
	\caption{$L_2$-error of displacement for numerical example with undistorted mesh}
	\label{fig:ud}
\end{figure}
Fig.~\ref{fig:ud} depicts the results for the undistorted mesh, whose initial configuration is plotted in Fig.~\ref{fig:msh}~(a).
With decreasing thickness, locking effects evolve employing the standard formulation~(std.). The lower the approximation degree, the more severe the occurring locking. For the present investigation, this effect sets in for a slenderness ratio $\frac{L}{t}=100$. For the thinnest case, considerable locking occurs, even for $p=5$. The results of the standard formulation do not converge to the correct solution for $p=2$  but develop a plateau of an almost constant error. Beginning from $p=3$, convergence is obtained for fine meshes. However, the error level is at least two orders of magnitude larger than for all considered mixed formulations.
The standard mixed formulation not involving any lumping procedure~(mxd.) delivers the most accurate and stable results, which possess the expected convergence rates. Specifically, for $p=3$, the accuracy of the standard formulation~(std.) is diminished by almost four orders of magnitude compared to the corresponding mixed variants~(mxd.).
Already for this simplest mesh, the results of the formulation variant involving NURBS-lumping~(lmp.) reveal severe error plateaus. The error due to the employed NURBS-lumping procedure decreases with decreasing thickness, but is a sufficient reason not to recommend this basic lumping procedure due to a lack of reliability and the fact that it only provides sufficiently accurate results for very thin cases, independently of the employed approximation order. This raises the need for an adapted lumping procedure that offers reliable results while maintaining the advantage of increased efficiency.
For this simplest mesh with unlimited internal continuity, the variants with lumping based on approximate duals~(AD) and enhanced approximate duals~(eAD) yield identical transformation matrices~$\boldsymbol{T }_\alpha$ and thus provide identical results. In contrast to the variant involving NURBS-lumping~(lmp.), both variants involving dual-based lumping~(AD, eAD) yield optimal convergence rates. For the two thickest cases, the error level is slightly higher compared to the mixed formulation without lumping~(mxd.), but, for the two medium-thin cases, the results are identical to the computationally more expensive standard mixed formulation~(mxd.). For the thinnest case with a corresponding slenderness ratio of $\frac{L}{t}=10,000$, the quality of the solution deteriorates slightly after an error-threshold of approximately $5\cdot10^{-10}$ or $1\cdot10^{-10}$ is reached for $p=4$ or $p=5$, respectively. We attribute this to a minor deterioration in the condition number.  Considering that this effect occurs only for very fine meshes of the extremely thin case, and that the magnitude of error is just slightly higher for the variants employing approximate dual basis functions~(AD, eAD) compared to the mixed variant for very thick cases, the proposed variants employing approximate dual basis functions~(AD, eAD) can be recommended for all thicknesses and refinement ranges, as the involved lumping increases the efficiency due to the reduced solution time.

\FloatBarrier
\subsection{Distorted NURBS mesh}
\label{sec:quadraticPlate_N}
\FloatBarrier
\begin{figure}[thp]
	\centering
	\begin{subfigure}{0.49\textwidth}
		\includegraphics[width=\textwidth]{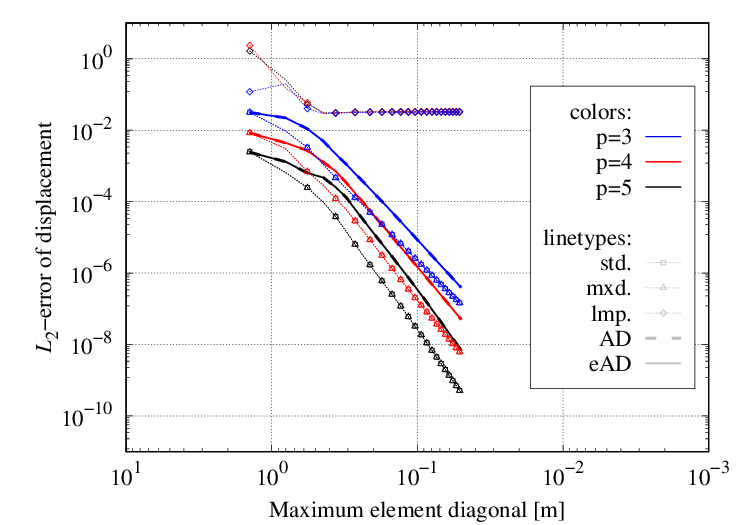}
		\caption{$t=1$~m}
	\end{subfigure}
	\hfill
	\begin{subfigure}{0.49\textwidth}
		\includegraphics[width=\textwidth]{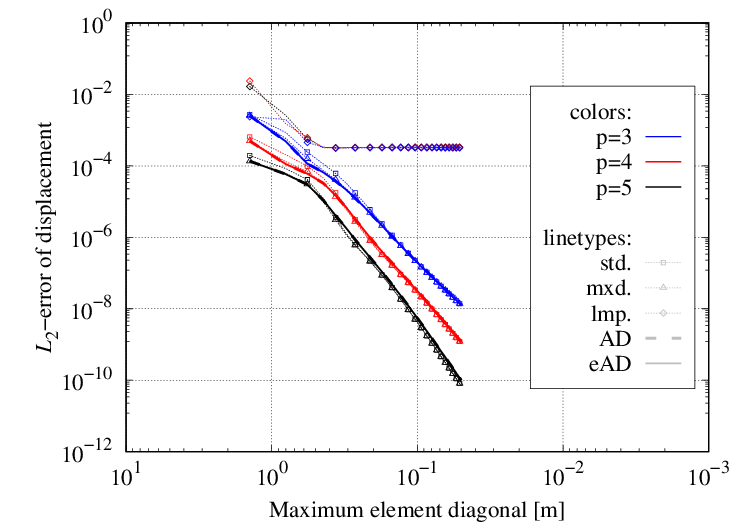}
		\caption{$t=0.1$~m}
	\end{subfigure}
	\\
	\begin{subfigure}{0.49\textwidth}
		\includegraphics[width=\textwidth]{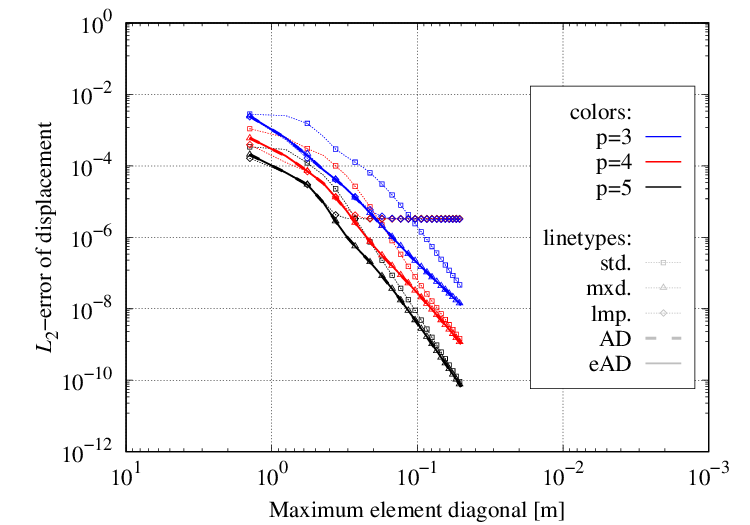}
		\caption{$t=0.01$~m}
	\end{subfigure}
	\hfill
	\begin{subfigure}{0.49\textwidth}
		\includegraphics[width=\textwidth]{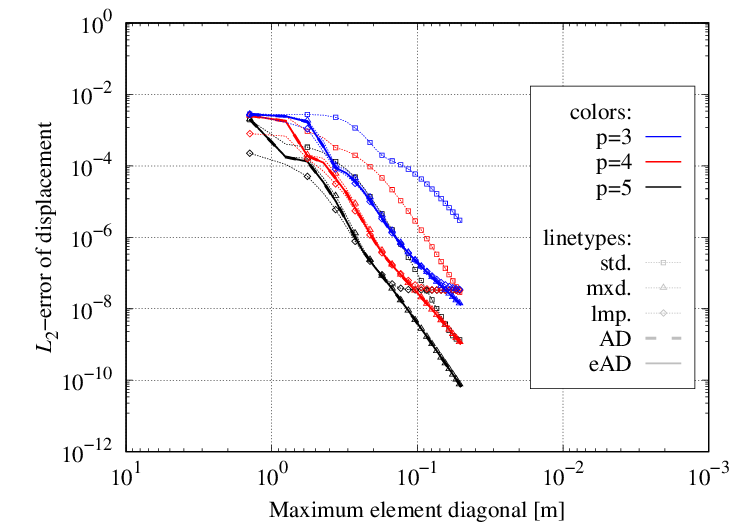}
		\caption{$t=0.001$~m}
	\end{subfigure}
	\\
	\begin{subfigure}{0.49\textwidth}
		\includegraphics[width=\textwidth]{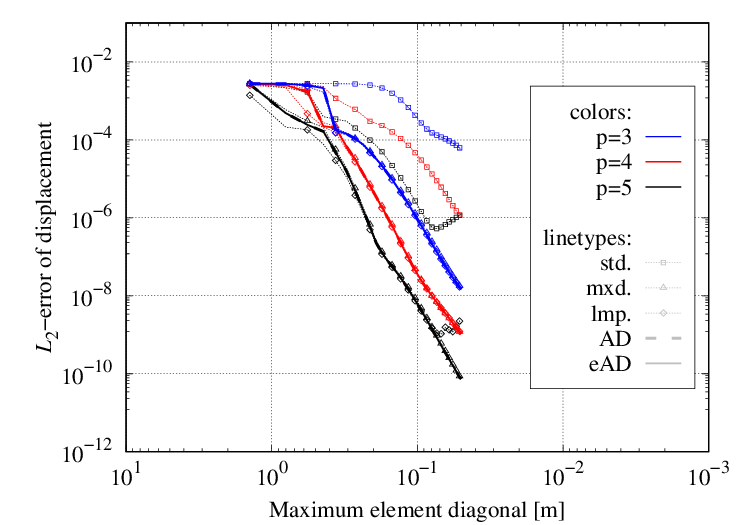}
		\caption{$t=0.0001$~m}
	\end{subfigure}
	\caption{$L_2$-error of displacement for numerical example with NURBS mesh}
	\label{fig:w}
\end{figure}

\begin{figure}[thp]
	\centering
	\begin{subfigure}{0.49\textwidth}
		\includegraphics[width=\textwidth]{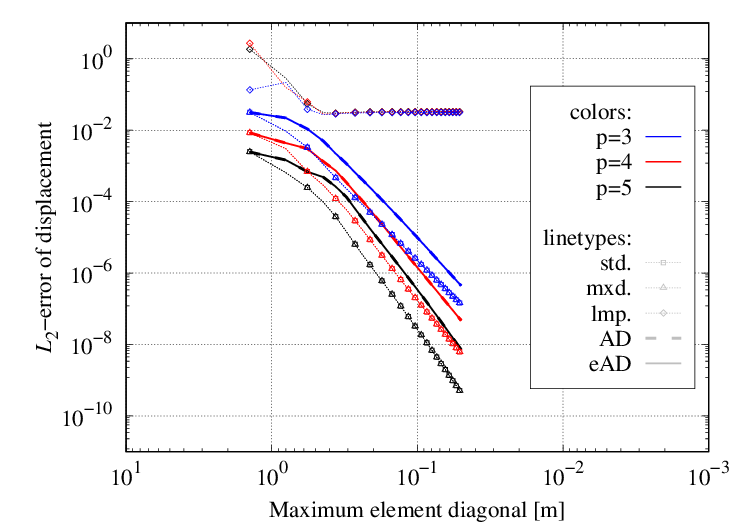}
		\caption{$t=1$~m}
	\end{subfigure}
	\hfill
	\begin{subfigure}{0.49\textwidth}
		\includegraphics[width=\textwidth]{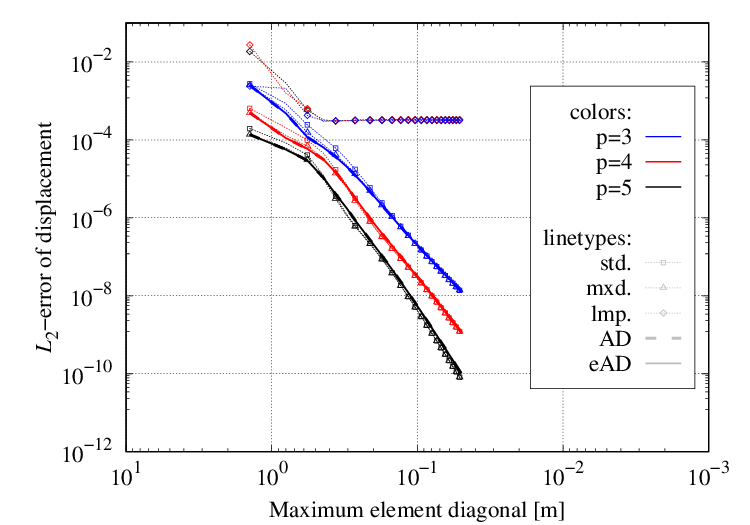}
		\caption{$t=0.1$~m}
	\end{subfigure}
	\\
	\begin{subfigure}{0.49\textwidth}
		\includegraphics[width=\textwidth]{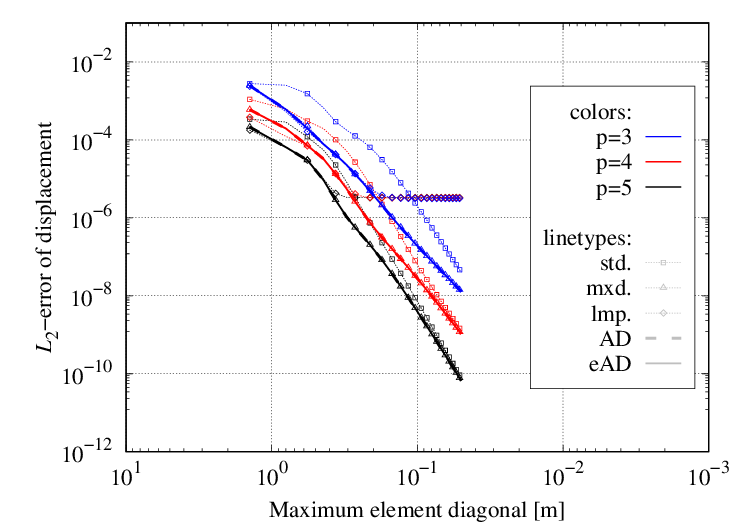}
		\caption{$t=0.01$~m}
	\end{subfigure}
	\hfill
	\begin{subfigure}{0.49\textwidth}
		\includegraphics[width=\textwidth]{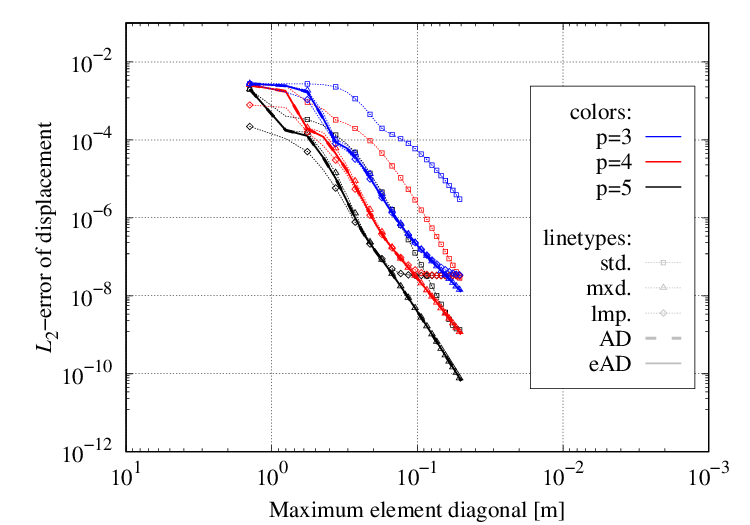}
		\caption{$t=0.001$~m}
	\end{subfigure}
	\\
	\begin{subfigure}{0.49\textwidth}
		\includegraphics[width=\textwidth]{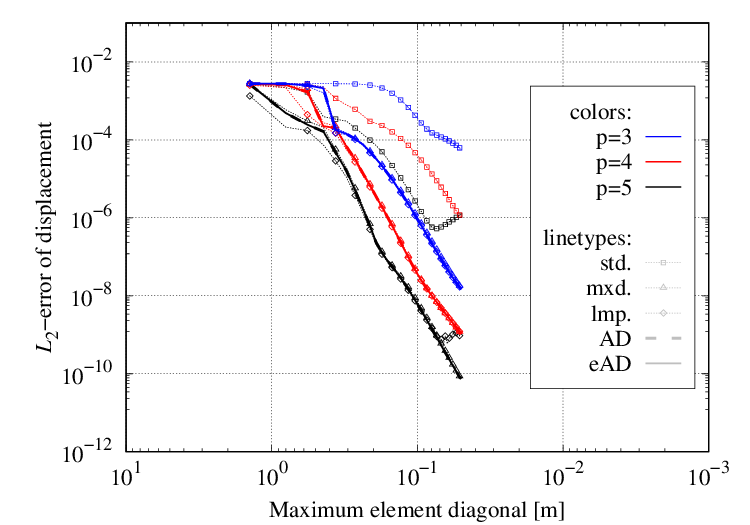}
		\caption{$t=0.0001$~m}
	\end{subfigure}
	\caption{$L_2$-error of displacement for numerical example with NURBS mesh: Study on neglection of weights for stress basis functions}
	\label{fig:w_nw}
\end{figure}

Fig.~\ref{fig:w} depicts the results for a distorted NURBS mesh with full internal continuity, of which the coarsest instantiation is plotted in Fig.~\ref{fig:msh}~(b). This discretization both demonstrates the effect of mesh distortion and verifies the implementation of NURBS and their dual basis functions.
Like for the undistorted mesh, investigated in Sec.~\ref{sec:quadraticPlate_ud}, locking effects evolve with decreasing thickness, employing the standard formulation~(std.). Again, the observed locking effect is more severe, the lower the approximation degree is.
Also for the distorted mesh, the standard mixed formulation not involving any lumping procedure~(mxd.) provides the most accurate and stable results, which possess the expected convergence rates.
As for the undistorted mesh, the results of the formulation variant involving NURBS-lumping~(lmp.) reveal severe error plateaus whose error magnitudes decrease with decreasing thickness. This confirms the previous finding that this unadapted lumping procedure is not recommended and an adapted lumping procedure is required.
Again, as expected, the variants with lumping based on approximate duals~(AD) and enhanced approximate duals~(eAD) provide identical results with optimal convergence rates. Like for the undistorted mesh, for the two thickest cases, the error level is slightly higher compared to the mixed formulation without lumping~(mxd.). In contrast, for all remaining thinner cases, the results are almost identical to the computationally more expensive standard mixed formulation~(mxd.). For the two thinnest variants, the initial error of the dual variants~(AD, eAD) is even slightly lower than that of the standard mixed formulation~(mxd.). Furthermore, no numerical instabilities occur for the distorted mesh, which makes the lumping-variants employing approximate dual basis functions (AD, eAD) even more recommendable for a general application to increase the efficiency of mixed formulations.

Fig.~\ref{fig:w_nw} provides the results for a simplified formulation, where B-splines (instead of NURBS) are used as basis functions for the interpolation of the stress parameters $\boldsymbol{S}_\alpha$. The results are visually identical, and hence no deteriorating effects can be observed in comparison to the results in Fig.~\ref{fig:w}. Thus, in order to simplify the implementation of the proposed method slightly, it is possible to omit the weighting factors of the NURBS basis functions for the additional shear parameters. However, this study will also be conducted in Sec.~\ref{sec:example:mp:multicon} for the most severe NURBS discretization of the investigated benchmark example.

\subsection{Internal $C^1$-continuity}
\begin{figure}[t]
	\centering
	\begin{subfigure}{0.49\textwidth}
		\includegraphics[width=\textwidth]{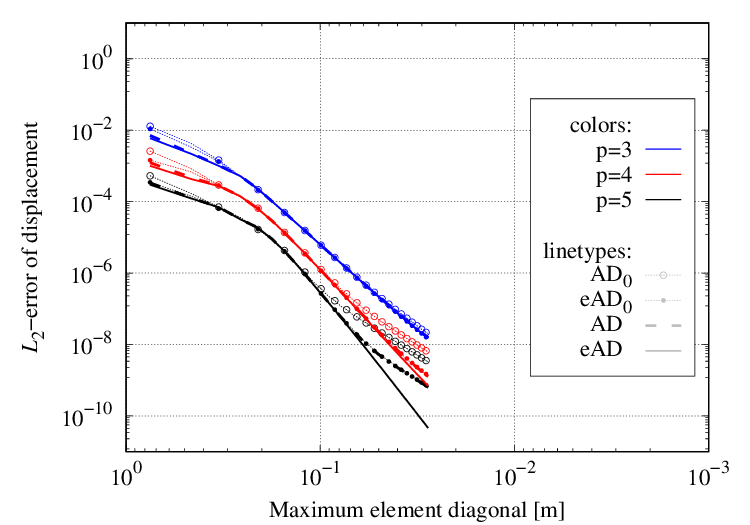}
		\caption{$t=1$~m}
	\end{subfigure}
	\hfill
	\begin{subfigure}{0.49\textwidth}
		\includegraphics[width=\textwidth]{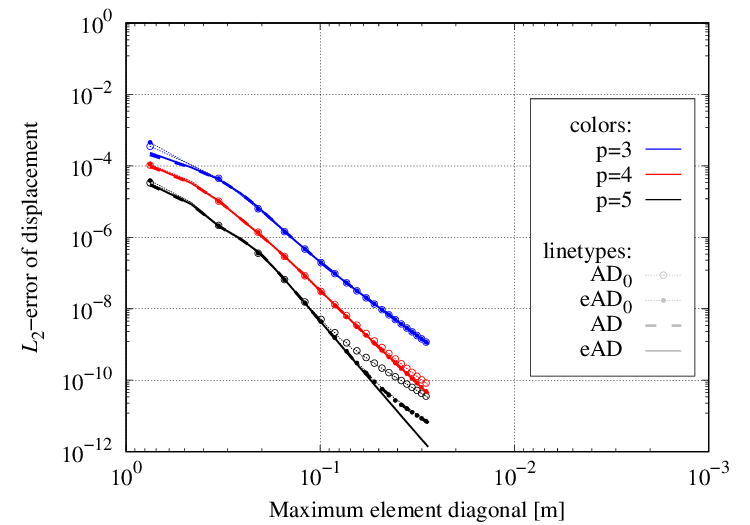}
		\caption{$t=0.1$~m}
	\end{subfigure}
	\\
	\begin{subfigure}{0.49\textwidth}
		\includegraphics[width=\textwidth]{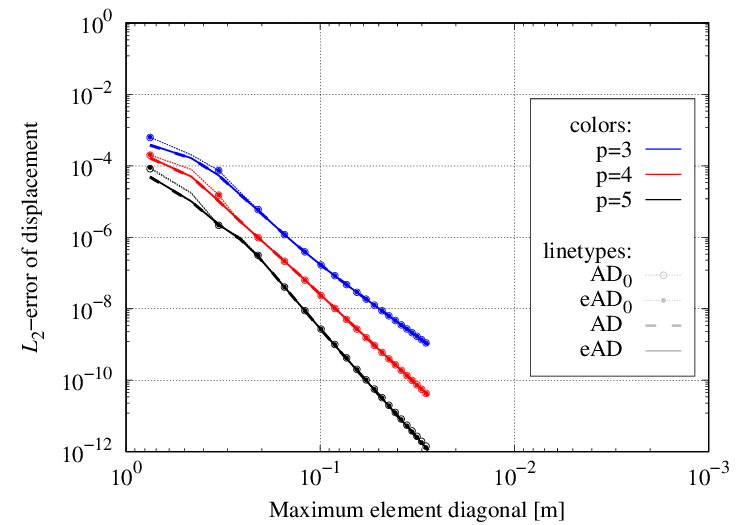}
		\caption{$t=0.01$~m}
	\end{subfigure}
	\hfill
	\begin{subfigure}{0.49\textwidth}
		\includegraphics[width=\textwidth]{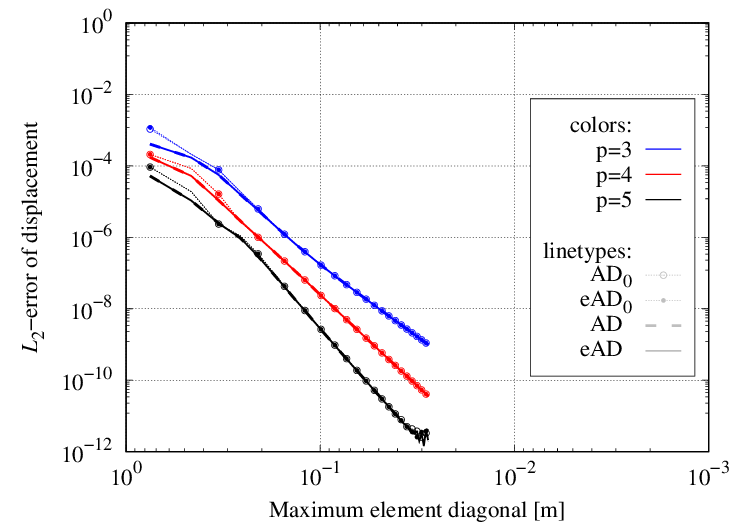}
		\caption{$t=0.001$~m}
	\end{subfigure}
	\\
	\begin{subfigure}{0.49\textwidth}
		\includegraphics[width=\textwidth]{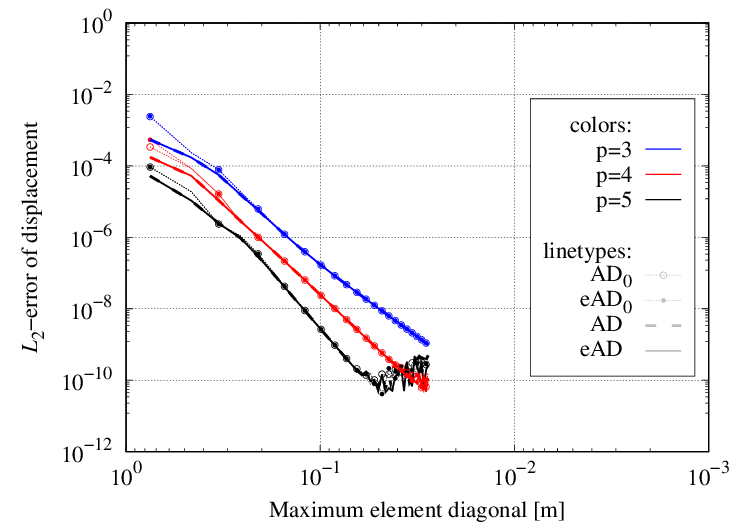}
		\caption{$t=0.0001$~m}
	\end{subfigure}
	\caption{$L_2$-error of displacement for numerical example with $C^1$-continuity: Study on requirement of knot insertion}
	\label{fig:C1_C}
\end{figure}

In this discretization with internal $C^1$-continuity, of which the coarsest mesh is depicted in Fig.~\ref{fig:msh}~(c), the requirement for the reduction of continuity at internal knots with limited continuity, as defined in Sec.~\ref{sec:reductioncontinuity}, is shown. The results in Fig.~\ref{fig:C1_C} clearly confirm that without this reduction of continuity,  the AD and eAD lumping variants yield suboptimal convergence rates in the fine limit. This effect increases for higher approximation orders. For the graphs denoted by an index '0' (AD$_0$ and eAD$_0$) in Fig.~\ref{fig:C1_C}, the knot vector is not altered, i.e.\ the initial continuity is maintained. The results for the two thickest variants in Fig.~\ref{fig:C1_C}~(a)-(b) show a clear deterioration of the convergence rate. For the graphs denoted by AD and eAD, the continuity is reduced as proposed in Sec.~\ref{sec:reductioncontinuity}, as is done in the remainder of this paper. The results are clearly improved compared to the AD$_0$ and eAD$_0$ results for the two thickest variants, respectively.  The change from AD$_0$ to AD reduces the error level but the convergence rate seems to be rather unchanged. Only if we use the eAD variant, which is proposed in Sec.~\ref{sec:enhancedAD} in combination with the continuity reduction of Sec.~\ref{sec:reductioncontinuity}, the lumping schemes achieve optimal convergence rates. For the thinner versions in Fig.~\ref{fig:C1_C}~(c)-(e), only minimal differences are visible. The slightly lower error level for coarse meshes for the variants with reduced continuity is attributed to the additional degrees of freedom introduced by this procedure.
Slight numerical instabilities occur for $p=4$ and $p=5$ for very fine meshes and the thinnest variants, but considering the obtainable precision of the employed double precision numerics under consideration of the condition number, these oscillations are not surprising. Thus, this study clearly confirms that only the enhanced formulation with reduced internal continuity~(eAD) retains the optimal convergence rates for all considered thickness ratios.

\label{sec:C1_internal}

\begin{figure}[t]
	\centering
	\begin{subfigure}{0.49\textwidth}
		\includegraphics[width=\textwidth]{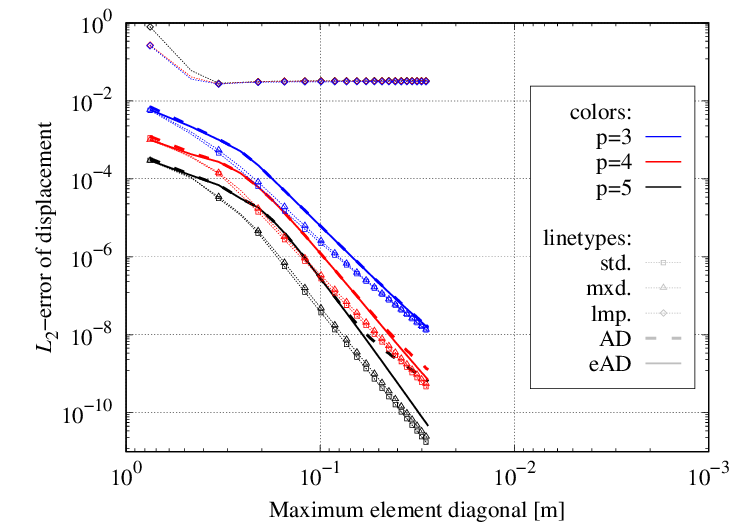}
		\caption{$t=1$~m}
	\end{subfigure}
	\hfill
	\begin{subfigure}{0.49\textwidth}
		\includegraphics[width=\textwidth]{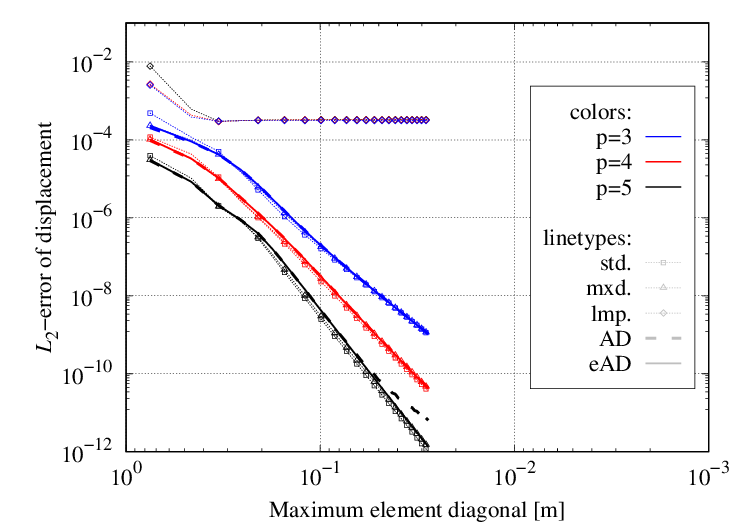}
		\caption{$t=0.1$~m}
	\end{subfigure}
	\\
	\begin{subfigure}{0.49\textwidth}
		\includegraphics[width=\textwidth]{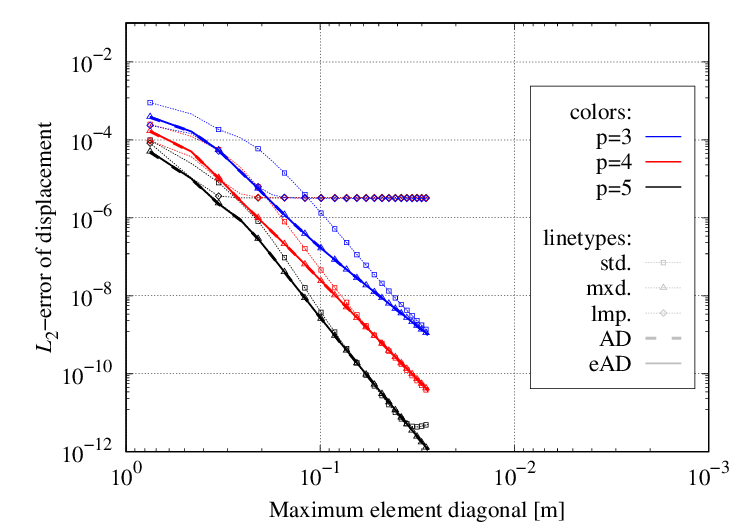}
		\caption{$t=0.01$~m}
	\end{subfigure}
	\hfill
	\begin{subfigure}{0.49\textwidth}
		\includegraphics[width=\textwidth]{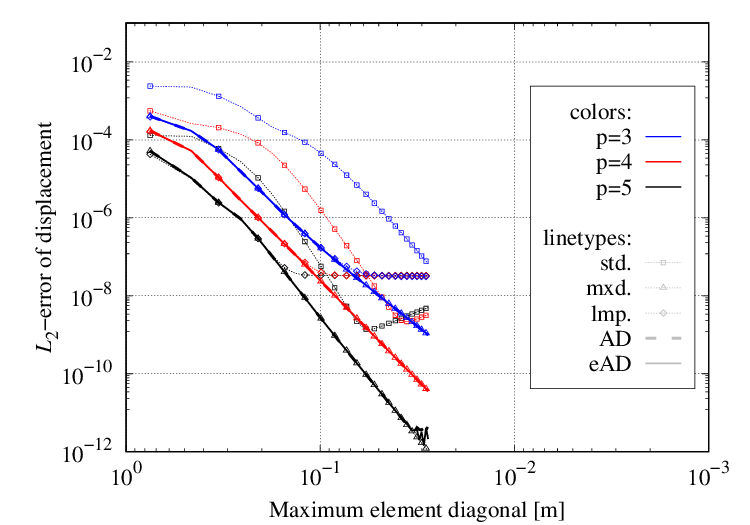}
		\caption{$t=0.001$~m}
	\end{subfigure}
	\\
	\begin{subfigure}{0.49\textwidth}
		\includegraphics[width=\textwidth]{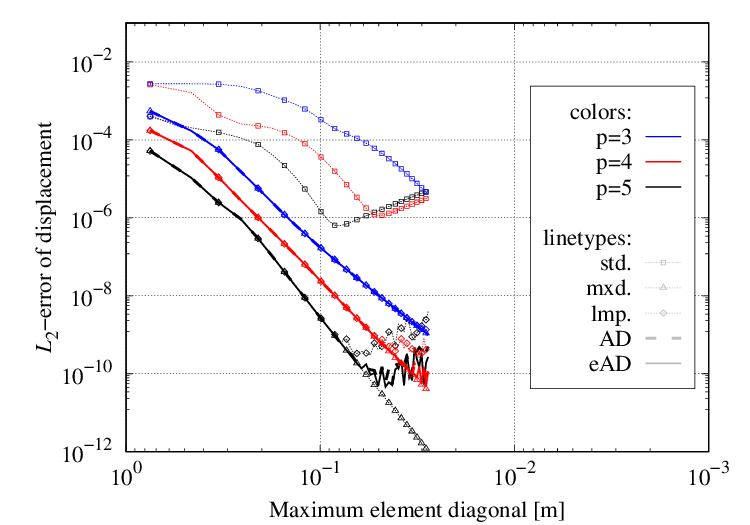}
		\caption{$t=0.0001$~m}
	\end{subfigure}
	\caption{$L_2$-error of displacement for numerical example with $C^1$-continuity}
	\label{fig:C1}
\end{figure}
\FloatBarrier
Fig.~\ref{fig:C1} compares the accuracy of the dual variants with required knot insertion~(AD, eAD) to that of the other examined variants. This clearly shows that enhanced approximate duals are required to properly treat geometries with internal $C^1$-continuity, which are very common in real-world CAD files.
As for the previous discretization types, the purely displacement-based method~(std.) leads to intensifying locking the lower the thickness and the lower the order of the basis functions is chosen. Especially for the thinnest case, severe locking occurs, yielding an error level that is several orders of magnitude higher compared to the other methods.
The standard mixed formulation not involving any lumping procedure~(mxd.) delivers the most accurate and stable results, which possess the expected convergence rates.
Like for the previously investigated discretizations, severe error plateaus occur for the formulation variant involving NURBS-lumping~(lmp.). Even though the error introduced decreases with decreasing thickness, this again confirms a sufficient reason not to employ this simplest lumping variant and emphasizes the need for an enhanced lumping procedure.
Unlike for the previously investigated meshes, the mesh with internal $C^1$-continuity yields different results for both dual variants~(AD, eAD). Notably, the variant with lumping based on approximate duals (AD) results in a loss in convergence rate for the two thickest cases given in Fig.~\ref{fig:C1}~(a)-(b). In contrast to that, the variant with lumping based on enhanced approximate duals~(eAD) manages to preserve the expected convergence rate. This effect even increases with an increasing approximation degree. Thus, the enhanced variant~(eAD) achieves almost equivalent results compared to the standard mixed formulation, for all but the thickest variant, where its resulting error level is slightly increased. Furthermore, for the thinnest variant, numerical instabilities occur for $p=4$ and $p,5$ after approaching an error magnitude of approximately $10^{-10}$. Conclusively, only the proposed variant involving enhanced approximate dual basis functions~(eAD) and the standard mixed formulation are capable of obtaining optimal convergence rates for all slenderness ratios.

\FloatBarrier
\subsection{Internal $C^0$-continuity}
\label{sec:C0_internal}
\FloatBarrier

\begin{figure}[thp]
\centering
\begin{subfigure}{0.49\textwidth}
\includegraphics[width=\textwidth]{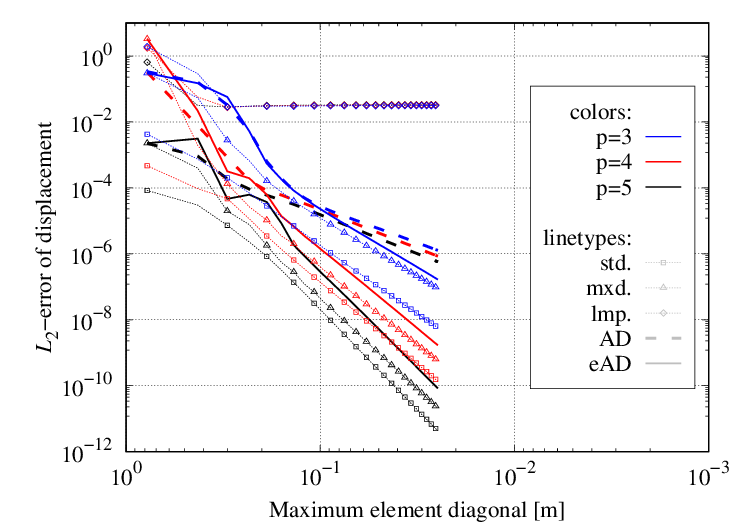}
\caption{$t=1$~m}
\end{subfigure}
\hfill
\begin{subfigure}{0.49\textwidth}
\includegraphics[width=\textwidth]{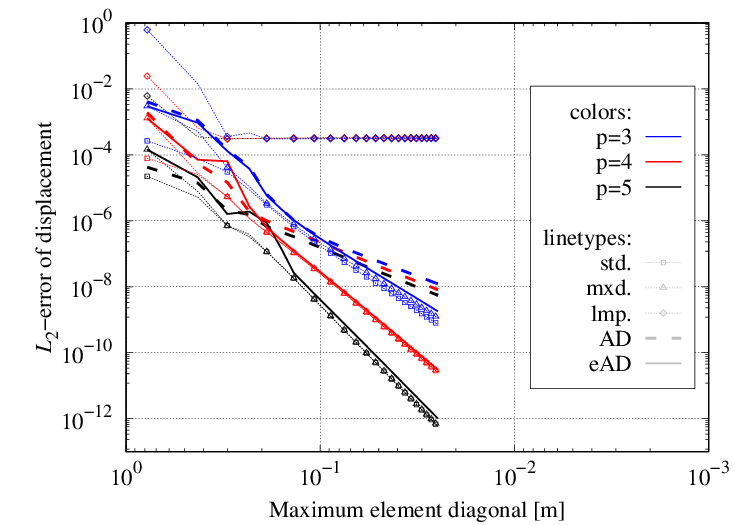}
\caption{$t=0.1$~m}
\end{subfigure}
\\
\begin{subfigure}{0.49\textwidth}
\includegraphics[width=\textwidth]{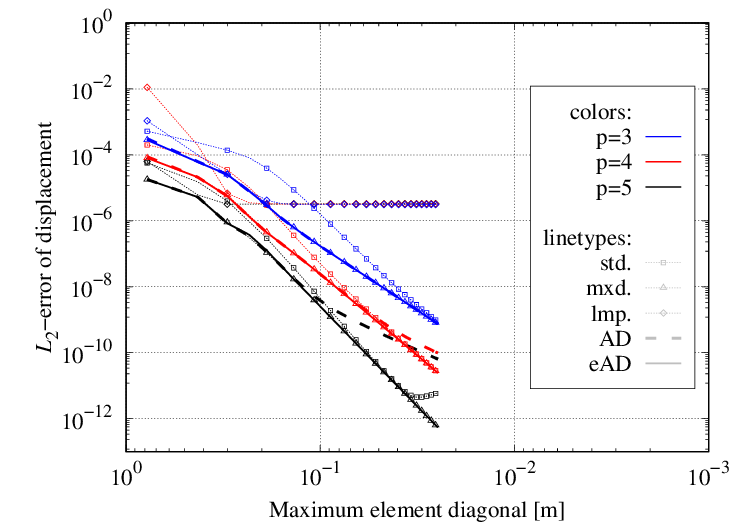}
\caption{$t=0.01$~m}
\end{subfigure}
\hfill
\begin{subfigure}{0.49\textwidth}
\includegraphics[width=\textwidth]{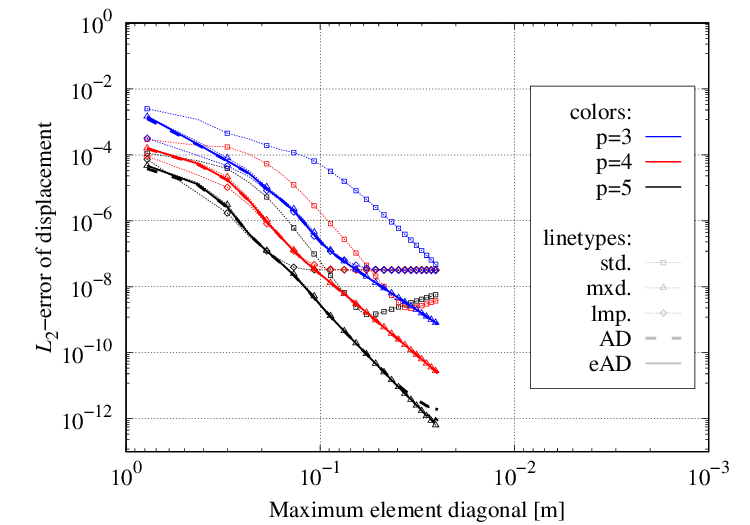}
\caption{$t=0.001$~m}
\end{subfigure}
\\
\begin{subfigure}{0.49\textwidth}
\includegraphics[width=\textwidth]{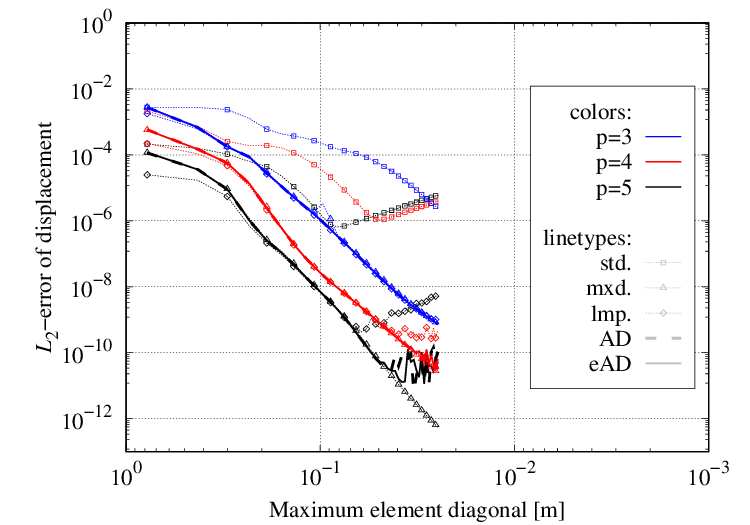}
\caption{$t=0.0001$~m}
\end{subfigure}
\caption{$L_2$-error of displacement for numerical example with $C^0$-continuity}
\label{fig:C0}
\end{figure}

Fig.~\ref{fig:C0} depicts the results for a distorted mesh with internal $C^0$-continuity, of which the coarsest instantiation is plotted in Fig.~\ref{fig:msh}~(d). Even if such geometries are rather uncommon in practical applications, they may nevertheless occur in some cases. Within this contribution, this example is employed to demonstrate the requirement of employing enhanced approximate dual basis functions~(eAD) to analyze discretizations that include such limited continuities reliably.
Again, with decreasing thickness, locking effects of similar extent as for the previously investigated mesh types evolve, if the standard formulation~(std.) is employed.
As investigated for all single-patch meshes, not involving any lumping procedure within a standard mixed formulation~(mxd.), provides the most accurate and stable results of the expected convergence rates.
Similar to the previous findings, a serious error plateau develops if a mixed formulation involving simple NURBS-lumping~(lmp.) is employed. This effect is independent of the approximation order but alleviates with decreasing thickness.
For this discretization with internal $C^0$-continuity, employing approximate dual basis functions~(AD) results in a significantly impaired convergence rate, especially for the thick variants with a corresponding slenderness ratio of $\frac{L}{t}=1$ and $\frac{L}{t}=10$. In particular, the convergence rate reduces to approximately~$1.5$, regardless of approximation order and thickness. This is an expected behavior that has also been observed in ~\cite{DornischEtAl2017,DornischStoeckler2021} within the scope of a mortar method. In contrast, the variant employing enhanced approximate duals~(eAD) is able to obtain optimal convergence rates. The refinement step in which this superiority sets in depends on the slenderness ratio and the approximation order. The thinner the plate and the lower the selected degree of basis functions, the finer the meshes for which the convergence gets impacted by this effect. However, even for the already thin case $t=0.01~\mathrm {m}$ and $p=5$, there is a significant difference of approximately two orders of magnitude in accuracy comparing the results of the formulation involving standard dual lumping~(AD) and enhanced dual lumping~(eAD). This confirms the need for an adapted lumping procedure~(eAD) when analyzing meshes with limited internal continuity, as only the standard mixed formulation~(mxd.) and this adapted variant~(eAD)  are able to obtain optimal convergence rates for all considered thickness ratios, even if numerical instabilities are introduced in the case of the thinnest variant after reaching an error threshold of $10^{-10}$ for $p=4$ and $p=5$.

Regarding the convergence of the shear stresses, comparable but slightly more unstable results can be achieved, where the order of magnitude of the error can also be higher for the dual variants (AD, eAD) compared to the standard entirely NURBS-based mixed formulation without lumping. Since the shear part of the stiffness matrix is lumped, the higher error level is not surprising.

\begin{figure}[tbp]
	\centering
	\begin{subfigure}{0.49\textwidth}
		\includegraphics[width=\textwidth]{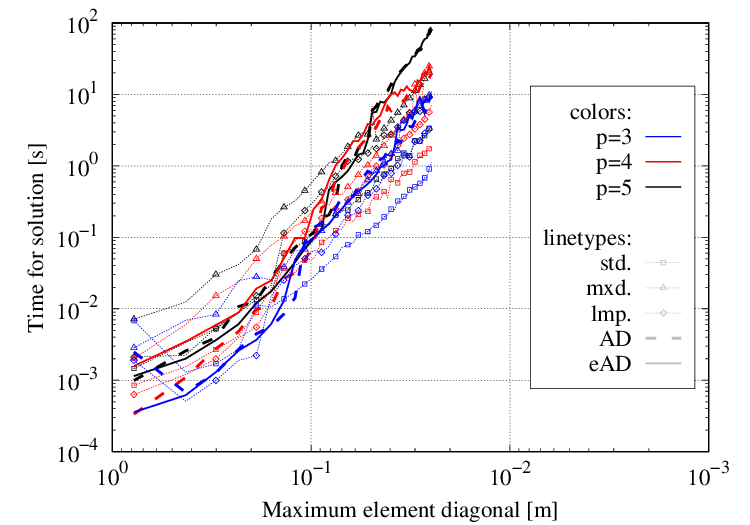}
		\caption{Solution times for $t=1$~m}
	\end{subfigure}
	\hfill
	\begin{subfigure}{0.49\textwidth}
		\includegraphics[width=\textwidth]{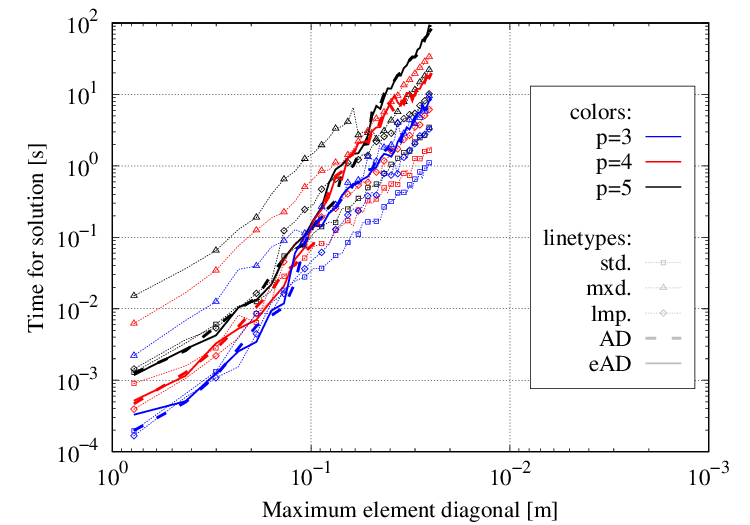}
		\caption{Solution times for $t=0.0001$~m}
	\end{subfigure}
	\\
	\begin{subfigure}{0.49\textwidth}
		\includegraphics[width=\textwidth]{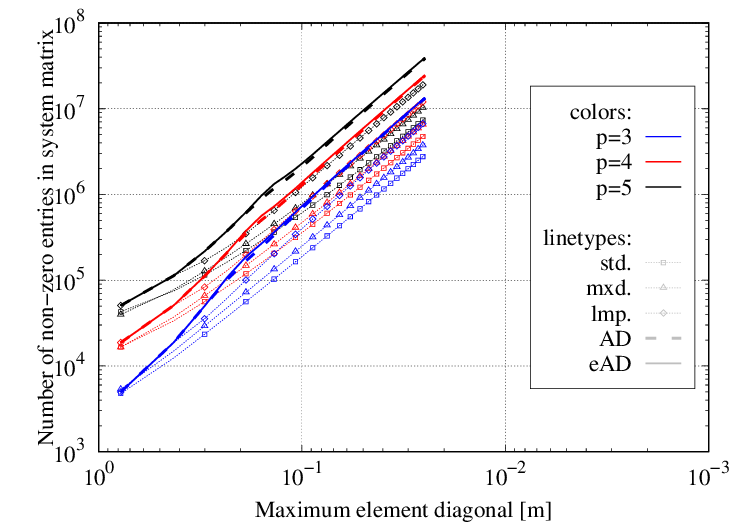}
		\caption{Number of non-zero entries of the system matrix}
	\end{subfigure}
	\hfill
	\begin{subfigure}{0.49\textwidth}
		\includegraphics[width=\textwidth]{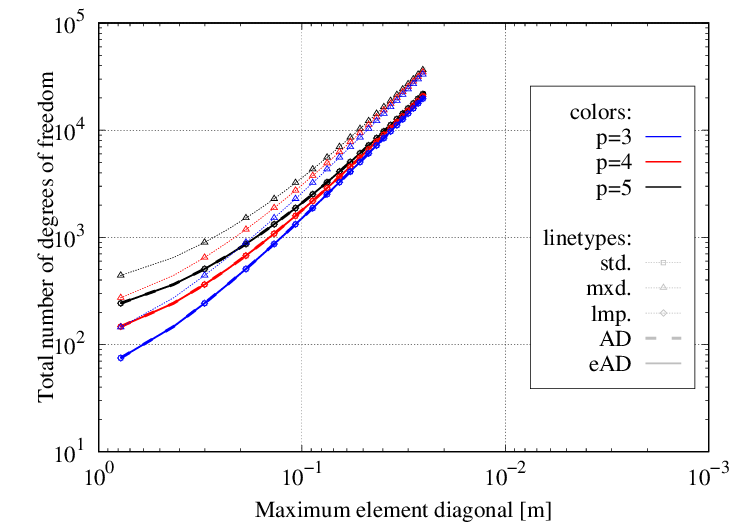}
		\caption{Number of degrees of freedom}
	\end{subfigure}
	\caption{Comparison of computational efficiency for numerical example with $C^0$-continuity}
	\label{fig:C0_eff}
\end{figure}

Furthermore, the effect of the proposed dual condensation procedure on the solution time is studied in Figs.~\ref{fig:C0_eff}~(a) and (b) for the thickest and the thinnest examined case, respectively. There, especially for coarse and moderately fine meshes, a huge advantage in solution time can be observed for the suggested static condensation method, regardless of whether approximate duals or enhanced approximate duals are employed. The aligning effects occurring for very fine meshes are probably attributable to the MATLAB$^\text{\textregistered}$-solver and are expected to reduce within a more efficient framework. To enable a more independent measure for the expected efficiency of the proposed procedure, the number of non-zero entries and the number of degrees of freedom of the system matrix are compared in Figs.~\ref{fig:mp_Cm_nw}~(c) and (d), respectively. These also show the general advantage of the proposed condensation method, which should also pay off regarding computation time.

\subsection{Multi-patch geometry with linear interface}
\label{sec:example:mp:linear}
\begin{figure}[thp]
	\centering
	\begin{subfigure}{0.49\textwidth}
		\includegraphics[width=\textwidth]{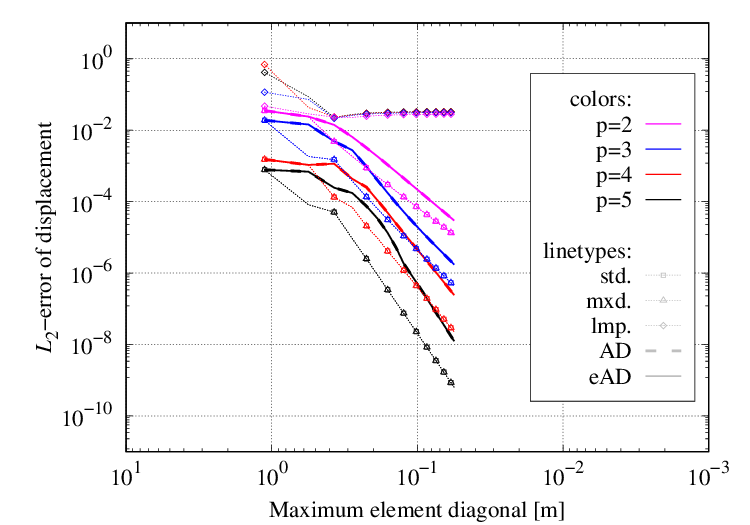}
		\caption{$t=1$~m}
	\end{subfigure}
	\hfill
	\begin{subfigure}{0.49\textwidth}
		\includegraphics[width=\textwidth]{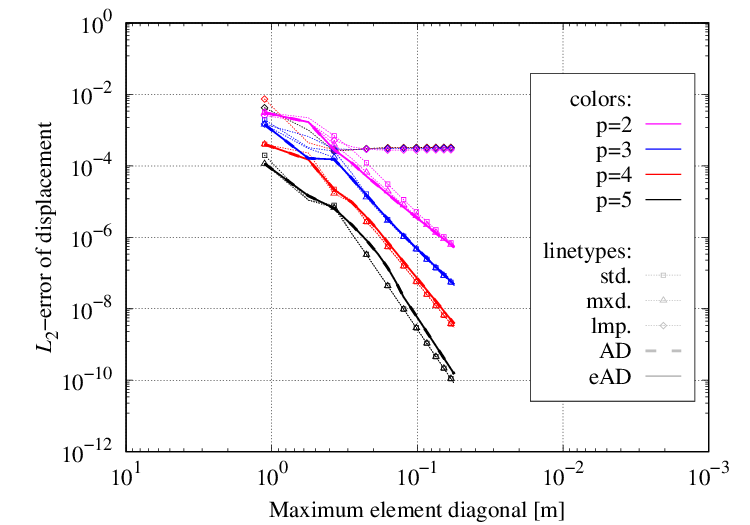}
		\caption{$t=0.1$~m}
	\end{subfigure}
	\\
	\begin{subfigure}{0.49\textwidth}
		\includegraphics[width=\textwidth]{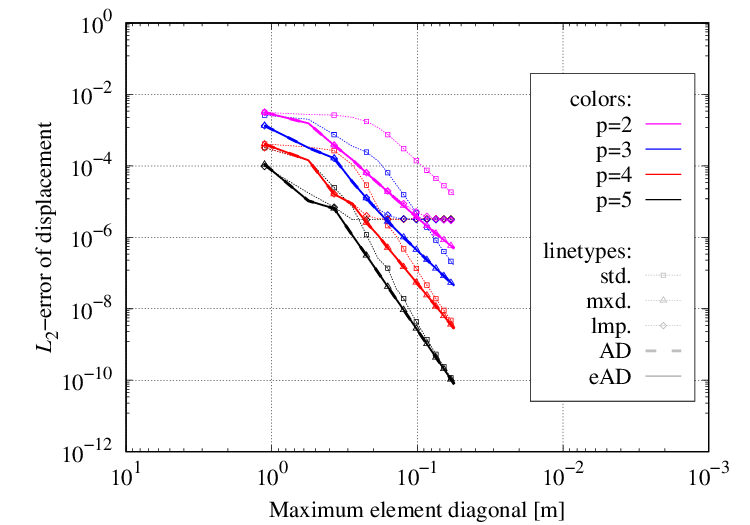}
		\caption{$t=0.01$~m}
	\end{subfigure}
	\hfill
	\begin{subfigure}{0.49\textwidth}
		\includegraphics[width=\textwidth]{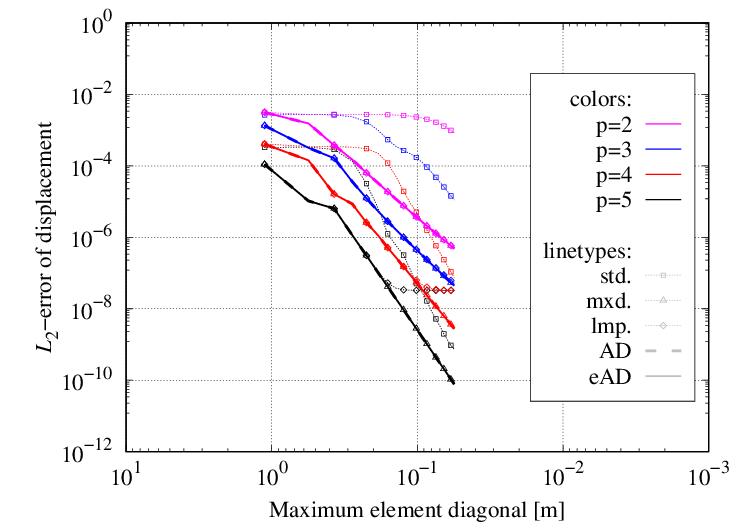}
		\caption{$t=0.001$~m}
	\end{subfigure}
	\\
	\begin{subfigure}{0.49\textwidth}
		\includegraphics[width=\textwidth]{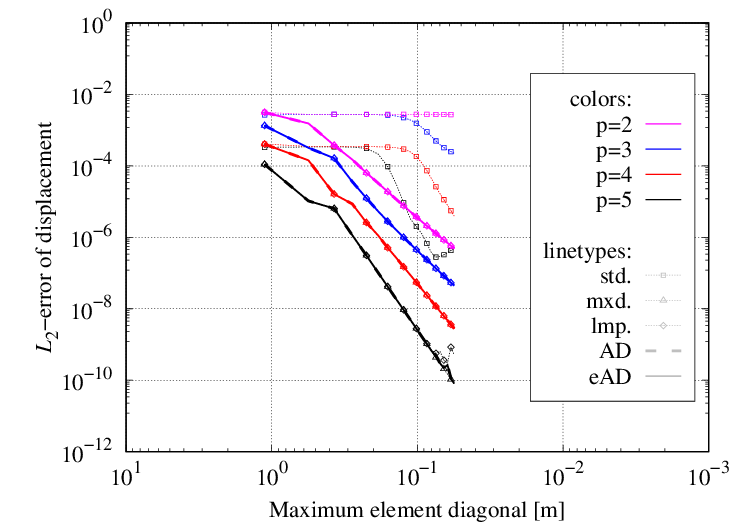}
		\caption{$t=0.0001$~m}
	\end{subfigure}
	\caption{$L_2$-error of displacement for numerical example with multi-patch geometry and linear interface}
	\label{fig:mp_l}
\end{figure}
Fig.~\ref{fig:mp_l} provides the results for a multi-patch mesh with linear interface, of which the coarsest discretization is depicted in Fig.~\ref{fig:msh_mp}~(a).
As investigated for all examined single-patch meshes, the standard formulation~(std.) suffers from severe locking that increases with diminishing thickness and decreasing approximation order for a simple multi-patch mesh with linear interface as well. The intensity of the observed locking effects is similar to that observed for the single-patch discretizations as well.
Furthermore, using a standard mixed formulation~(mxd.) leads to the most accurate and stable results with the required convergence rates, also for the investigated multi-patch geometry.
As observed for the investigated single-patch meshes, using a mixed formulation involving lumping (lmp.) results in a severe error plateau whose error magnitude increases with the thickness of the plate.
Analyzing this simple multi-patch mesh with unlimited internal continuity, as expected, the dual variants (AD, eAD) deliver identical results. Despite the minor numerical instability occurring for the thinnest variant and $p=5$ for one of the finest analyzed meshes, the results imply an applicability of the proposed enhanced dual lumping procedures (AD, eAD) to multi-patch settings as well, as both obtain the optimal convergence rates.

\subsection{Multi-patch geometry with $C^1$-continuity}
\label{sec:example:mp:C1}
\begin{figure}[thp]
\centering
\begin{subfigure}{0.49\textwidth}
	\includegraphics[width=\textwidth]{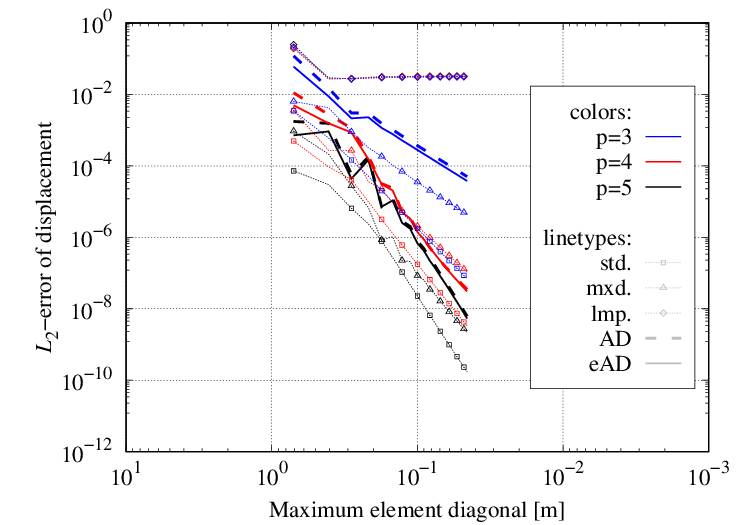}
	\caption{$t=1$~m}
\end{subfigure}
\hfill
\begin{subfigure}{0.49\textwidth}
	\includegraphics[width=\textwidth]{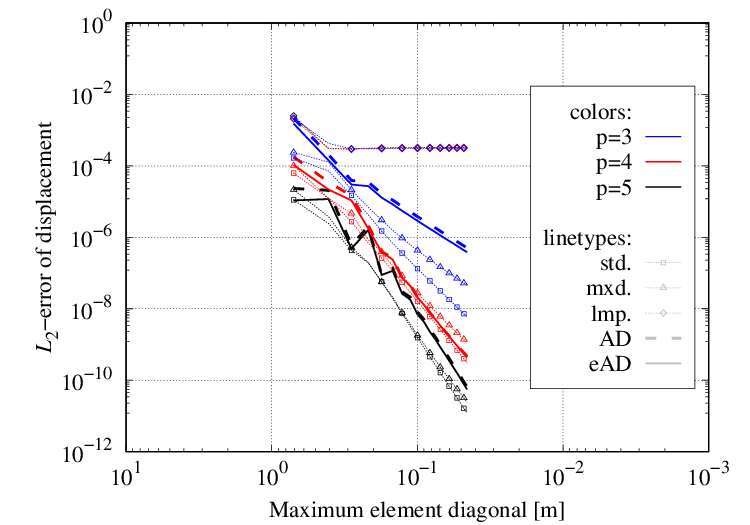}
	\caption{$t=0.1$~m}
\end{subfigure}
\\
\begin{subfigure}{0.49\textwidth}
	\includegraphics[width=\textwidth]{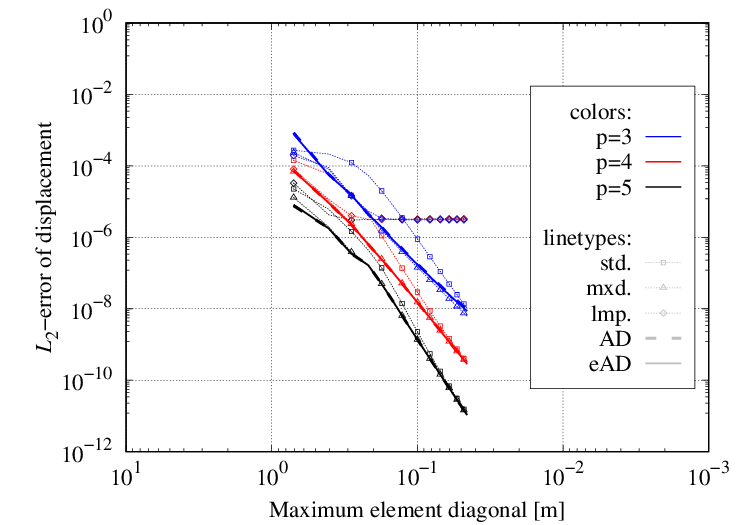}
	\caption{$t=0.01$~m}
\end{subfigure}
\hfill
\begin{subfigure}{0.49\textwidth}
	\includegraphics[width=\textwidth]{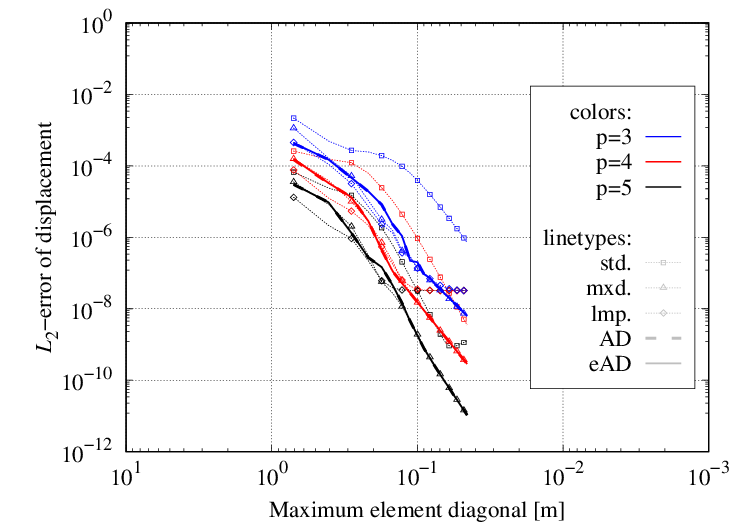}
	\caption{$t=0.001$~m}
\end{subfigure}
\\
\begin{subfigure}{0.49\textwidth}
	\includegraphics[width=\textwidth]{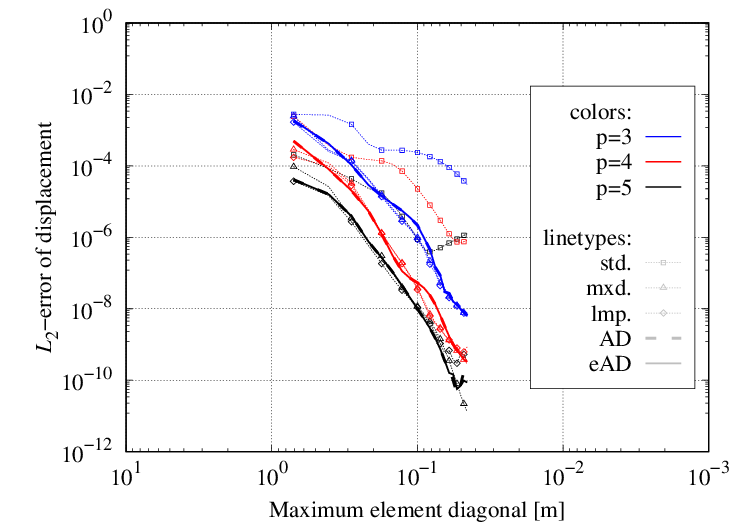}
	\caption{$t=0.0001$~m}
\end{subfigure}
\caption{$L_2$-error of displacement for numerical example with multi-patch geometry with $C^1$-continuity}
\label{fig:mp_C1}
\end{figure}

Fig.~\ref{fig:mp_C1} depicts the results for a distorted multi-patch discretization with $C^1$-continuity, which was also examined in \cite{DornischStoeckler2021}. The coarsest investigated mesh of this type is depicted in Fig.~\ref{fig:msh_mp}~(b).
Regarding the proneness to locking of the standard method~(std.), similar effects, that increase for lower approximation orders and decrease with decreasing thickness, can also be observed for the examined multi-patch geometry with $C^1$-continuity.
While, like for the previous discretization types, the standard mixed method~(mxd.) delivers preferable results for thinner variants (Fig.~\ref{fig:mp_C1}~(c)-(e)), the standard method~(std.) achieves almost identical (Fig.~\ref{fig:mp_C1}~(b)) or even better results (Fig.~\ref{fig:mp_C1}~(a)) for thicker plates.
Also, for this discretization type, the requirement for an adapted lumping procedure is confirmed, as severe error plateaus, whose error magnitudes decrease with decreasing thickness, evolve when employing a mixed formulation involving simple NURBS-lumping~(lmp.) as well.

In contrast to the observations made in Sec.~\ref{sec:C1_internal} for single-patch meshes with $C^1$-continuity, for this multi-patch discretization, the results of both dual variants (AD, eAD) correspond very well. Nevertheless, the enhanced formulation (eAD) yields a slightly smaller error for higher thicknesses, regardless of the interpolation order. To explain why in this example the use of approximate duals (AD) is sufficient to yield optimal convergence rates, we point to the discretizations displayed for this section and Sec.~\ref{sec:C1_internal} in Figs.~\ref{fig:msh_mp}~(b) and~\ref{fig:msh}~(c), respectively. The former is only slightly distorted, while the latter is highly distorted with a large jump in the second derivatives at the central knot.
Despite some oscillations for coarse meshes, both dual variants~(AD and eAD) provide the expected convergence rates for fine meshes, independently of the thickness, although their error levels are slightly higher compared to the standard mixed formulation (mxd.) for the thick cases in Fig.~\ref{fig:mp_C1}~(a)-(b). In thinner cases, AD, eAD and mxd.\ alternate in performing the most accurate.

\subsection{Multi-patch geometry with various continuities}
\label{sec:example:mp:multicon}

\begin{figure}[thp]
\centering
\begin{subfigure}{0.49\textwidth}
	\includegraphics[width=\textwidth]{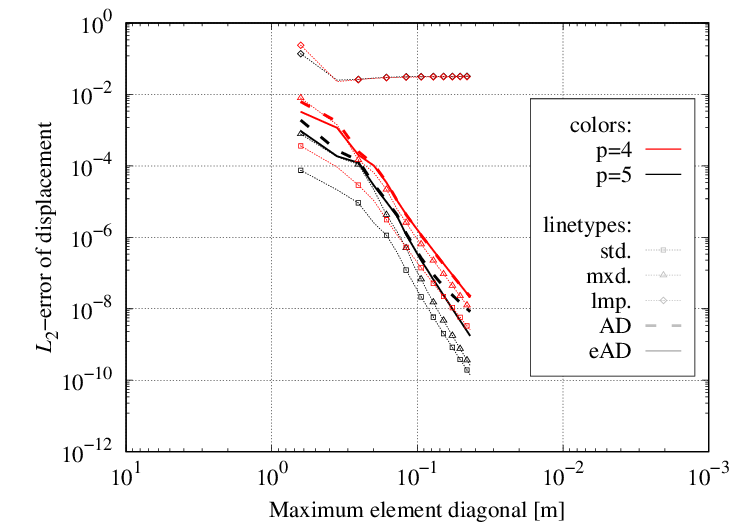}
	\caption{$t=1$~m}
\end{subfigure}
\hfill
\begin{subfigure}{0.49\textwidth}
	\includegraphics[width=\textwidth]{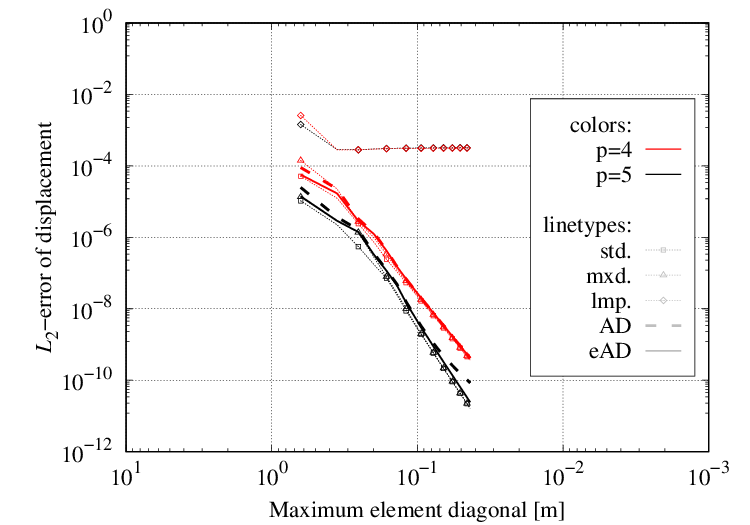}
	\caption{$t=0.1$~m}
\end{subfigure}
\\
\begin{subfigure}{0.49\textwidth}
	\includegraphics[width=\textwidth]{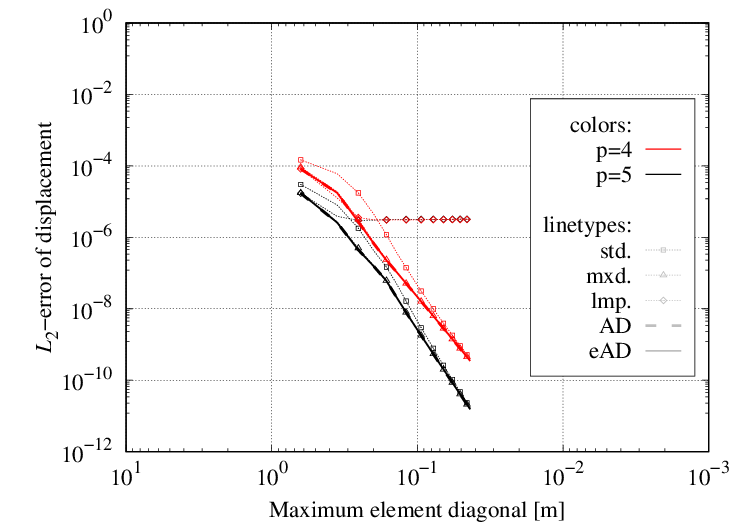}
	\caption{$t=0.01$~m}
\end{subfigure}
\hfill
\begin{subfigure}{0.49\textwidth}
	\includegraphics[width=\textwidth]{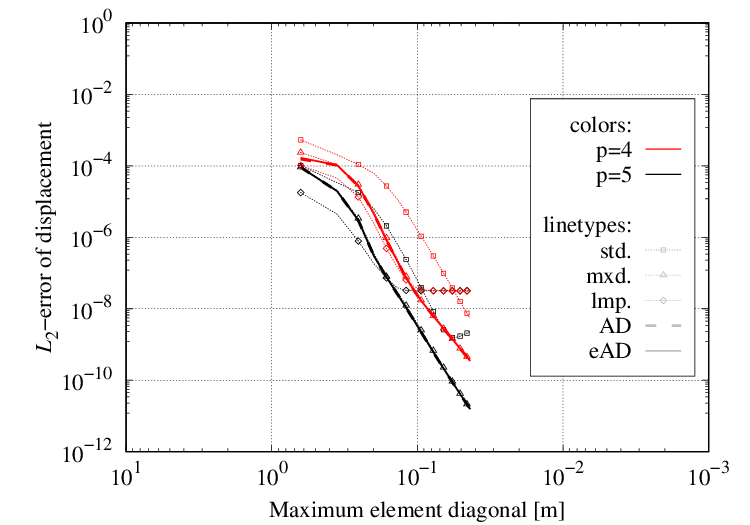}
	\caption{$t=0.001$~m}
\end{subfigure}
\\
\begin{subfigure}{0.49\textwidth}
	\includegraphics[width=\textwidth]{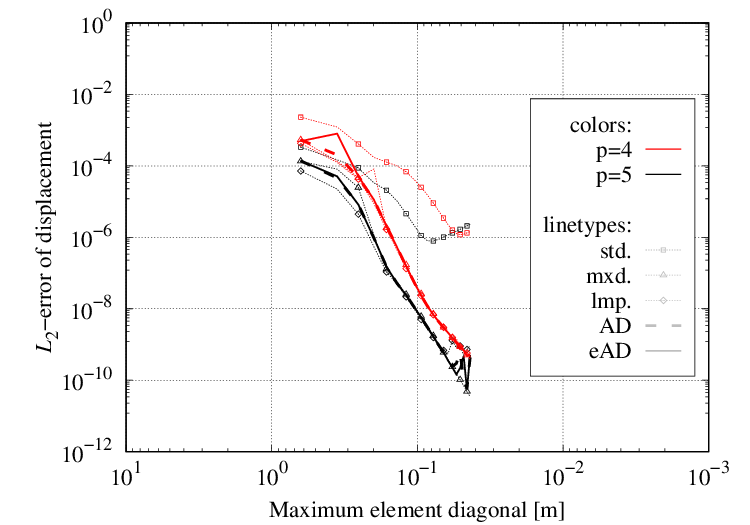}
	\caption{$t=0.0001$~m}
\end{subfigure}
\caption{$L_2$-error of displacement for numerical example with multi-patch geometry with various continuities}
\label{fig:mp_Cm}
\end{figure}

\begin{figure}[thp]
	\centering
	\begin{subfigure}{0.49\textwidth}
		\includegraphics[width=\textwidth]{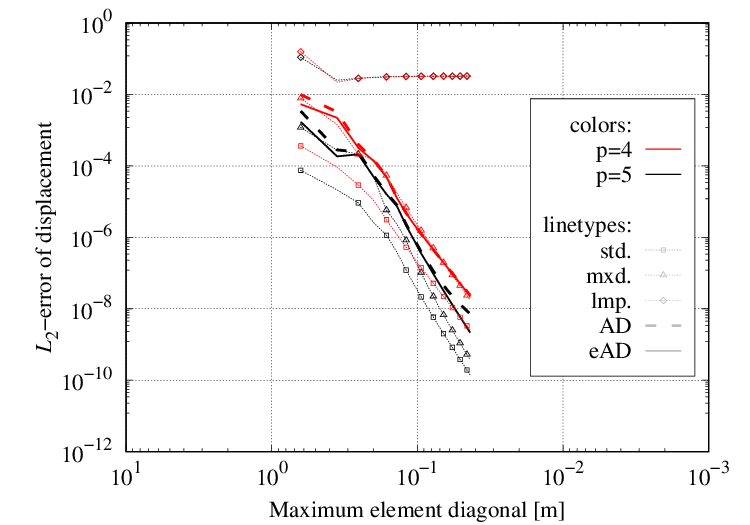}
		\caption{$t=1$~m}
	\end{subfigure}
	\hfill
	\begin{subfigure}{0.49\textwidth}
		\includegraphics[width=\textwidth]{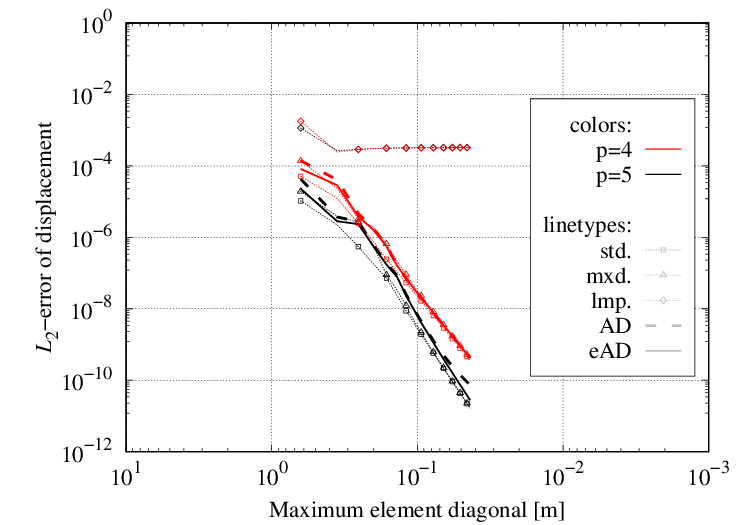}
		\caption{$t=0.1$~m}
	\end{subfigure}
	\\
	\begin{subfigure}{0.49\textwidth}
		\includegraphics[width=\textwidth]{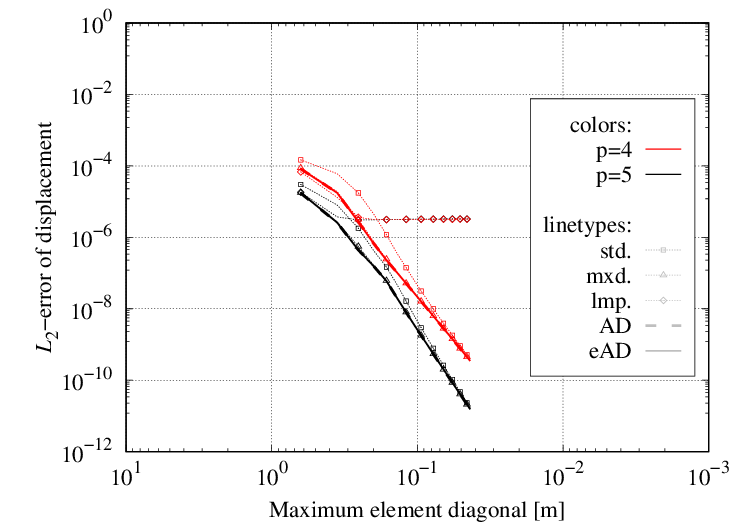}
		\caption{$t=0.01$~m}
	\end{subfigure}
	\hfill
	\begin{subfigure}{0.49\textwidth}
		\includegraphics[width=\textwidth]{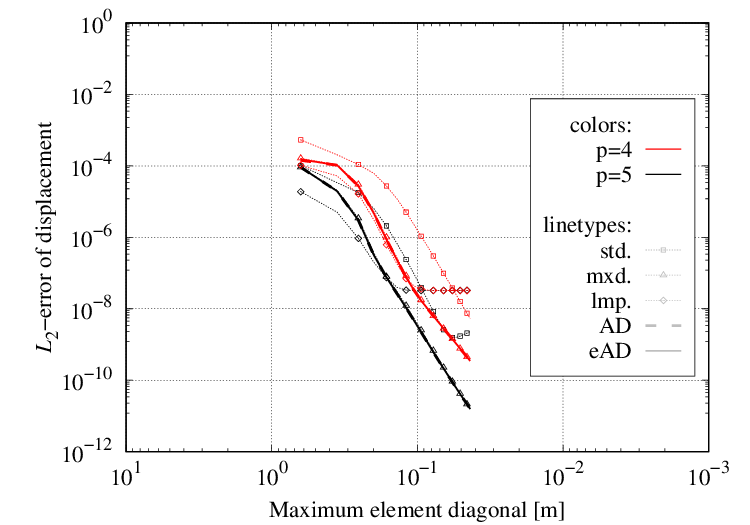}
		\caption{$t=0.001$~m}
	\end{subfigure}
	\\
	\begin{subfigure}{0.49\textwidth}
		\includegraphics[width=\textwidth]{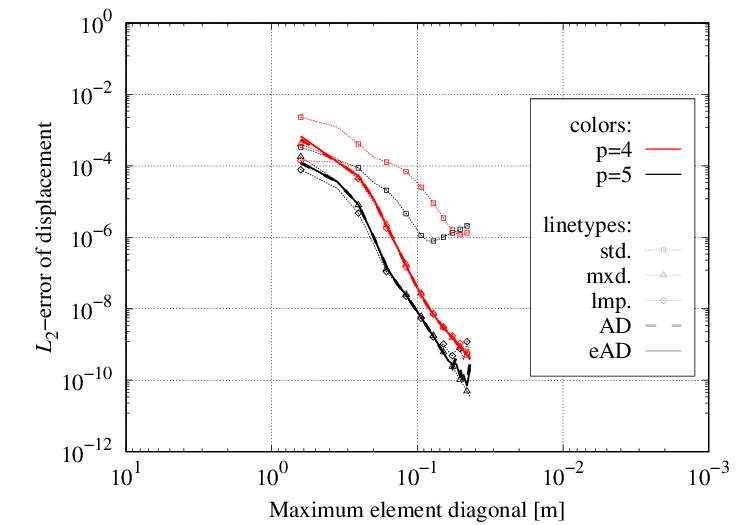}
		\caption{$t=0.0001$~m}
	\end{subfigure}
	\caption{$L_2$-error of displacement for numerical example with multi-patch geometry with various continuities: Study on neglection of weights for stress basis functions}
	\label{fig:mp_Cm_nw}
\end{figure}

Fig.~\ref{fig:mp_Cm} depicts the results for a distorted multi-patch mesh with multiple internal continuities, of which the definition was taken from \cite{DornischStoeckler2021}. Its coarsest discretization is depicted in Fig.~\ref{fig:msh_mp}~(c).
Like for all investigated meshes, the standard formulation~(std.) suffers from increasingly severe locking, the thinner the plate.
The standard mixed method without any involved lumping procedure~(mxd.) delivers the most reliable results for the thin cases~(Fig.~\ref{fig:mp_Cm}~(c)-(e)), while the standard primal method~(std.) achieves a slightly lower error for the thicker plates~(Fig.~\ref{fig:mp_Cm}~(a)-(b)).
As observed for all previously examined discretization types, the mixed formulation conducting lumping for an entirely NURBS-based system matrix~(lmp.) proves to be not recommendable due to the occurring serious error plateau, which is a consequence of the introduced lumping error, even though it decreases with decreasing slenderness.
When comparing the dual variants (AD, eAD), it becomes obvious that for the two thickest cases provided in Fig.~\ref{fig:mp_Cm}~(a)-(b), only the enhanced formulation~(eAD) retains the optimal convergence rate for $p=5$, while it deteriorates when using approximate dual basis functions (AD). Furthermore, for the initial refinement steps, the error of the enhanced variant~(eAD) is slightly lower than that of the unenhanced dual variant (AD). Although minor numerical instabilities occur for fine meshes of the thinnest case and $p=5$, the proposed mixed formulation that involves lumping based on enhanced approximate dual basis functions~(eAD) proves to be generally applicable for the investigated single-patch and multi-patch discretizations, even if those possess points of limited internal continuity.

In Fig.~\ref{fig:mp_Cm_nw}, the stress parameters are interpolated by B-spline basis functions. Thus, it is assessed if the use of B-splines instead of NURBS, which is equivalent to omitting the weights of the corresponding control points, has an impact on the obtainable accuracy. Similarly to Sec.~\ref{sec:quadraticPlate_N}, no significant negative effects are observed, as the results are almost identical to those in Fig.~\ref{fig:mp_Cm}. In the two thickest cases, the error of all mixed formulations is slightly higher if the weight factors are omitted. This confirms the previous findings from Sec.~\ref{sec:quadraticPlate_N} and indicates that the weights can be neglected in order to allow for a slightly more efficient implementation.

Regarding the convergence of the shear stresses, the proposed dual formulations offer optimal convergence rates, but the order of magnitude of the error can be higher compared to the standard entirely NURBS-based mixed formulation without lumping. Since the shear part of the stiffness matrix is lumped, the higher error level is not surprising.

\begin{figure}[tbp]
	\centering
	\begin{subfigure}{0.49\textwidth}
		\includegraphics[width=\textwidth]{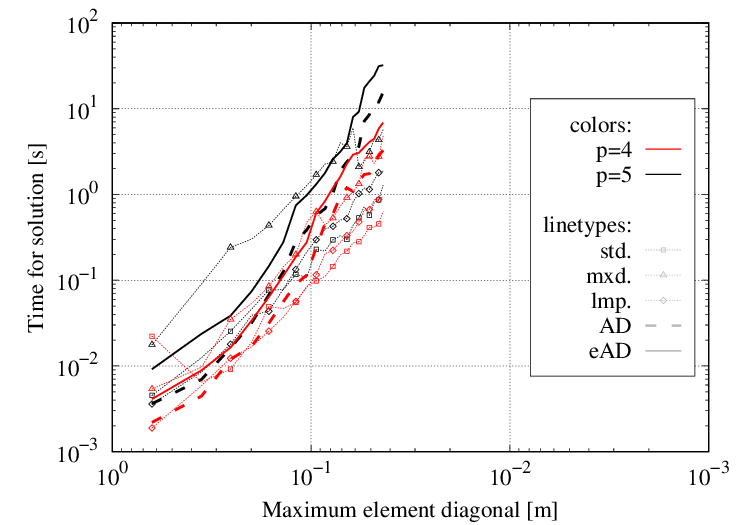}
		\caption{Solution times for $t=1$~m}
	\end{subfigure}
	\hfill
	\begin{subfigure}{0.49\textwidth}
		\includegraphics[width=\textwidth]{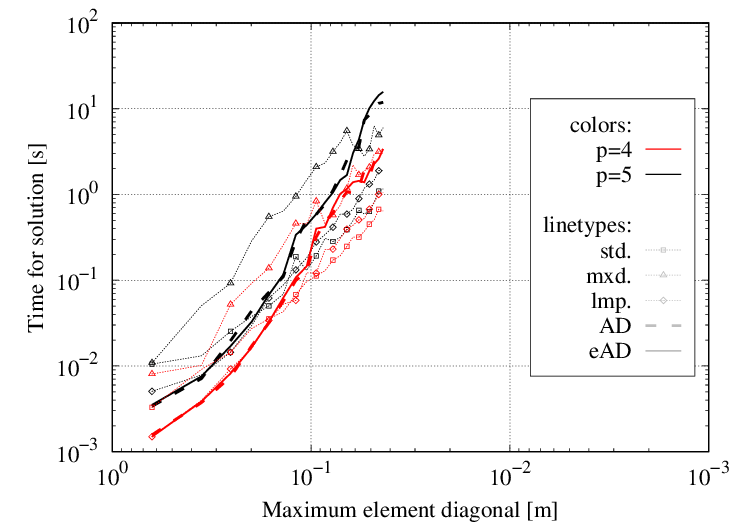}
		\caption{Solution times for $t=0.0001$~m}
	\end{subfigure}
	\\
	\begin{subfigure}{0.49\textwidth}
		\includegraphics[width=\textwidth]{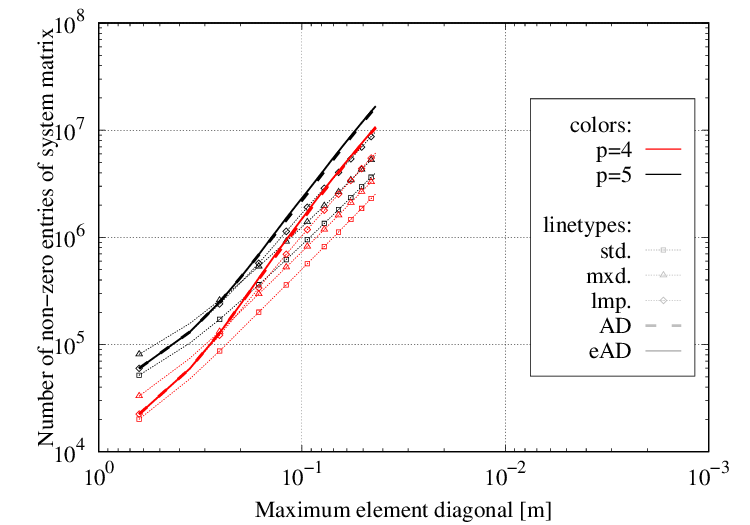}
		\caption{Number of non-zero entries of the system matrix}
	\end{subfigure}
	\hfill
	\begin{subfigure}{0.49\textwidth}
		\includegraphics[width=\textwidth]{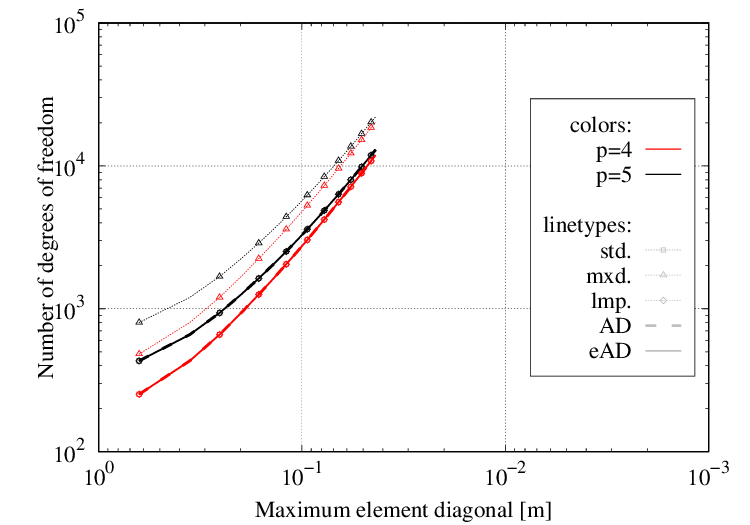}
		\caption{Number of degrees of freedom}
	\end{subfigure}
	\caption{Comparison of computational efficiency for numerical example with multi-patch geometry with various continuities}
	\label{fig:mp_Cm_eff}
\end{figure}

Additionally, for the smallest and the highest thickness studied, the influence of the suggested static condensation method with dual lumping on the computation time required for the solution of the global system of equations is examined in Figs.~\ref{fig:mp_Cm_eff}~(a) and~(b). As for the single-patch example with $C^0$-continuity, a considerable benefit is noticed for both variants of the proposed condensation procedure, particularly for low and moderate degrees of refinement. Probably due to the employed MATLAB$^\text{\textregistered}$-framework, leveling impacts can be observed for the finest meshes of this discretization as well. However, these probably vanish for a more efficient solver. In order to evaluate the efficiency of the proposed method more realistically, Figs.~\ref{fig:mp_Cm_nw}~(c) and (d) compare the number of non-zero entries of the system matrix as well as the number of degrees of freedom. These findings confirm those already made in Sec.~\ref{sec:C0_internal}.

\section{Summary and conclusions}
In this contribution, an efficient static condensation procedure is proposed and examined within the scope of a mixed isogeometric plate formulation. As approximate dual basis functions are employed for the discretization of the virtual parameters of the additional fields, the corresponding matrix part of the system matrix is diagonal-dominant and can thus be lumped to a diagonal matrix at a minimal loss of accuracy. Factorizing the governing equations appropriately, this sub-matrix becomes the unit matrix, which avoids the calculation of an inverse within the proposed static condensation procedure. Hence, the superior convergence and un-locking behavior of a mixed formulation can be achieved with the same number of degrees of freedom as a standard primal formulation, but with minimally increased computation time and at a slight loss of accuracy. This indicates that an even bigger benefit could result for a non-linear formulation that involves solving the global equation in every Newton-Raphson step.
Furthermore, it has been shown that B-Splines can be employed for the interpolation of the mixed fields instead of NURBS basis functions without introducing a significant error. As a central point of this contribution, the requirement of a proper treatment of the initial knot vector for geometries including points of limited internal continuities is shown for both single-patch and multi-patch examples. Additionally, the use of  enhanced approximate dual basis functions is proposed in order to counteract deteriorated convergence rates that are observed when employing approximate dual basis functions.

Further investigations could include an extension towards other types of formulations like non-linear and shell formulations, for instance.
Changing to a more efficient framework would allow to assess the savings in computation time. Additionally, performing static condensation on element-level would be a further desirable enhancement of the present formulation.

\section*{CRediT authorship contribution statement}

\textbf{Lisa Stammen}: Writing -- review and editing, Writing -- original draft, Visualization, Validation, Software, Methodology, Investigation, Formal analysis

\textbf{Wolfgang Dornisch}: Writing -- review and editing, Writing -- original draft, Validation, Supervision, Software, Methodology, Investigation, Funding acquisition, Formal analysis, Conceptualization

\section*{Acknowledgement}
The authors greatly acknowledge the fruitful discussions and cooperation with Prof.~Joachim St\"ockler from TU Dortmund University. The authors thank for providing the scripts for the computation of the extended approximate dual basis functions.

\section*{Funding sources}
This research did not receive any specific grant from funding agencies in the public, commercial, or not-for-profit sectors.

\section*{Declaration of competing interest}
The authors declare that they have no known competing financial interests or personal relationships that could have appeared to influence the work reported in this paper.

\section*{Data availability}
No data was used for the research described in the article.
\newpage
\appendix
\section{Geometry definitions}
\label{apndx}
The discretizations used for the numerical examples can be recovered from the data for knot vectors and control points provided in the following. Furthermore, we provide downloadable IGES files for an easier use of these benchmark discretizations by other researchers.
The initial surface orders and knot vectors of the examined single-patch geometries are as follows:
\begin{equation}
	\footnotesize
	\begin{aligned}
		\text{Sec.~\ref{sec:quadraticPlate_ud}:}&&p&=q=1\:;&&\boldsymbol{\Xi}=\boldsymbol{H}=\left(0,0,1,1\right) \\
		\text{Sec.~\ref{sec:quadraticPlate_N}:}&&p&=q=2\:;&&\boldsymbol{\Xi}=\boldsymbol{H}=\left(0,0,0,1,1,1\right) \\
		\text{Sec.~\ref{sec:C1_internal}:}&&p&=q=2\:;&&\boldsymbol{\Xi}=\boldsymbol{H}=\left(0,0,0,0.5,1,1,1\right) \\
		\text{Sec.~\ref{sec:C0_internal}:}&&p&=q=2\:;&&\boldsymbol{\Xi}=\boldsymbol{H}=\left(0,0,0,0.5,0.5,1,1,1\right)
	\end{aligned}
\end{equation}
The investigated single-patch-meshes including control points and control polygons are depicted in Fig.~\ref{fig:msh}.
The coordinates of control points and their weights are listed in Tab.~\ref{tab:geo}.

Regarding the investigated multi-patch-geometries, surfaces with the following definitions per patch are investigated:
\begin{equation}
	{\scriptstyle \footnotesize
	\begin{aligned}
		 \text{Sec.~\ref{sec:example:mp:linear}:}&&p=1\:\text{,}\:&q=1\:;&&\boldsymbol{\Xi}=\left(0,0,1,1\right)\text{,}\:\boldsymbol{H}=\left(0,0,1,1\right) \\
		 \text{Sec.~\ref{sec:example:mp:C1}:}&&p=1\:\text{,}\:&q=2\:;&&\boldsymbol{\Xi}=\left(0,0,1,1\right)\text{,}\:\boldsymbol{H}=\left(0,0,0,0.5,1,1,1\right) \\
		 \text{Sec.~\ref{sec:example:mp:multicon}:}&&p=1\:\text{,}\:&q=3\:;&&\boldsymbol{\Xi}=\left(0,0,1,1\right)\text{,}\:\boldsymbol{H}=\left(0,0,0,0,0.3,0.3,0.5,0.5,0.5,0.7,1,1,1,1\right)
	\end{aligned}
    }
\end{equation}
Fig.~\ref{fig:msh_mp} depicts the investigated meshes including control points and control polygons.
Tab.~\ref{tab:geo_mp} contains the coordinates of the control points and their weights.
\begin{table}[b]
	\begin{subtable}{0.24\textwidth}
		\centering
		\begin{tabular}{| c| c | c | c |}
			\hline
			x	&	y	&	z	&	w	\\
			\hline
			\hline
			0	&	0	&	0	&	1	\\
			1	&	0	&	0	&	1	\\
			\hline
			0	&	1	&	0	&	1	\\
			1	&	1	&	0	&	1	\\
			\hline
			\hline
		\end{tabular}
		\\~\\~\\~\\~\\~\\~\\~\\~\\~\\~\\~\\~\\~\\~\\~\\~\\~\\~\\~\\~\\~\\~\\
		\caption{Undistorted mesh \\(Sec.~\ref{sec:quadraticPlate_ud})}
	\end{subtable}
	\hfill
	\begin{subtable}{0.24\textwidth}
		\centering
		\begin{tabular}{| c| c | c | c |}
			\hline
			x	&	y	&	z	&	w \\
			\hline
			\hline
			0	&	0	&	0	&	1 	\\
			0.5	&	0	&	0	&	1 	\\
			1	&	0	&	0	&	1 	\\
			\hline
			0	&	0.5	&	0	&	1 	\\
			0.3	&	0.3	&	0	&	1.5 \\
			1	&	0.5	&	0	&	1 	\\
			\hline
			0	&	1	&	0	&	1	\\
			0.5	&	1	&	0	&	1 	\\
			1	&	1	&	0	&	1 	\\
			\hline
			\hline
		\end{tabular}
		\\~\\~\\~\\~\\~\\~\\~\\~\\~\\~\\~\\~\\~\\~\\~\\~\\~\\
		\caption{NURBS mesh \\(Sec.~\ref{sec:quadraticPlate_N})}
	\end{subtable}
	\hfill
	\begin{subtable}{0.24\textwidth}
		\centering
		\begin{tabular}{| c | c | c | c |}
			\hline
			x	&	y	&	z	&	w	\\
			\hline
			\hline
			0		&	0		&	0	&	1	\\
			0.25	&	0		&	0	&	1	\\
			0.75	&	0		&	0	&	1	\\
			1		&	0		&	0	&	1	\\
			\hline
			0		&	0.25	&	0	&	1	\\
			0.45	&	0.4		&	0	&	1	\\
			0.7		&	0.2		&	0	&	1	\\
			1		&	0.25	&	0	&	1	\\
			\hline
			0		&	0.75	&	0	&	1	\\
			0.2		&	0.9		&	0	&	1	\\
			0.5		&	0.6		&	0	&	1	\\
			1		&	0.75	&	0	&	1	\\
			\hline
			0		&	1		&	0	&	1	\\
			0.25	&	1		&	0	&	1	\\
			0.75	&	1		&	0	&	1	\\
			1		&	1		&	0	&	1	\\
			\hline
			\hline
		\end{tabular}
		\\~\\~\\~\\~\\~\\~\\~\\~\\~\\~\\
		\caption{$C^1$ - mesh \\(Sec.~\ref{sec:C1_internal})}
	\end{subtable}
	\hfill
	\begin{subtable}{0.24\textwidth}
		\centering
		\begin{tabular}{| c | c | c | c |}
			\hline
			x	&	y	&	z	&	w	\\
			\hline
			\hline
			0		&	0		&	0		&	1	\\
			0.25	&	0		&	0		&	1	\\
			0.5		&	0		&	0		&	1	\\
			0.75	&	0		&	0		&	1	\\
			1		&	0		&	0		&	1	\\
			\hline
			0		&	0.25	&	0		&	1	\\
			0.25	&	0.25	&	0		&	1	\\
			0.5		&	0.25	&	0		&	1	\\
			0.75	&	0.25	&	0		&	1	\\
			1		&	0.25	&	0		&	1	\\
			\hline
			0		&	0.5		&	0		&	1	\\
			0.3		&	0.55	&	0		&	1	\\
			0.45	&	0.45	&	0		&	1	\\
			0.65	&	0.45	&	0		&	1	\\
			1		&	0.5		&	0		&	1	\\
			\hline
			0		&	0.75	&	0		&	1	\\
			0.25	&	0.75	&	0		&	1	\\
			0.55	&	0.6		&	0		&	1	\\
			0.65	&	0.65	&	0		&	1	\\
			1		&	0.75	&	0		&	1	\\
			\hline
			0		&	1		&	0		&	1	\\
			0.25	&	1		&	0		&	1	\\
			0.5		&	1		&	0		&	1	\\
			0.75	&	1		&	0		&	1	\\
			1		&	1		&	0		&	1	\\
			\hline
			\hline
		\end{tabular}
		\caption{$C^0$ - mesh \\(Sec.~\ref{sec:C0_internal})}
	\end{subtable}
	\caption{Coordinates and weights of the control points of the investigated single-patch meshes}
	\label{tab:geo}
\end{table}
\FloatBarrier
\begin{table}[b]
	\begin{subtable}{0.32\textwidth}
		\centering
		\begin{tabular}{| c| c | c | c |}
			\hline
			x	&	y	&	z	&	w\\
			\hline
			\hline
			0	&	0	&	0	&	1\\
			0.5	&	0	&	0	&	1\\
			1	&	0	&	0	&	1\\
			\hline
			0	&	1	&	0	&	1\\
			0.5	&	1	&	0	&	1\\
			1	&	1	&	0	&	1\\
			\hline
			\hline
		\end{tabular}
		\\~\\~\\~\\~\\~\\~\\~\\~\\~\\~\\~\\~\\~\\~\\~\\~\\~\\~\\~\\~\\~\\~\\~\\~\\~\\
		\caption{Multi-patch-geometry with linear \\interface (Sec.~\ref{sec:example:mp:linear})}
	\end{subtable}
	\begin{subtable}{0.32\textwidth}
		\centering
		\begin{tabular}{| c | c | c | c |}
			\hline
			x	&	y	&	z	&	w\\
			\hline
			\hline
			0	&	0	 	&	0	&	1\\
			0.5	&	0	 	&	0	&	1\\
			1	&	0	 	&	0	&	1\\
			\hline
			0	&	0.25 	&	0	&	1\\
			0.6	&	0.3	 	&	0	&	1\\
			1	&	0.25	&	0	&	1\\
			\hline
			0	&	0.75	&	0	&	1\\
			0.4	&	0.7	 	&	0	&	1\\
			1	&	0.75 	&	0	&	1\\
			\hline
			0	&	1	 	&	0	&	1\\
			0.5	&	1	 	&	0	&	1\\
			1	&	1	 	&	0	&	1\\
			\hline
			\hline
		\end{tabular}
		\\~\\~\\~\\~\\~\\~\\~\\~\\~\\~\\~\\~\\~\\~\\~\\~\\~\\~\\~\\
		\caption{Multi-patch-geometry with \\$C^1$-continuity (Sec.~\ref{sec:example:mp:C1})}
	\end{subtable}
	\begin{subtable}{0.32\textwidth}
		\centering
		\begin{tabular}{| c | c | c | c |}
			\hline
			x	&	y	&	z	&	w\\
			\hline
			\hline
			0		&	0	&	0	&	1\\
			0.5		&	0	&	0	&	1\\
			1		&	0	&	0	&	1\\
			\hline
			0		&	0.1	&	0	&	1\\
			0.55	&	0.1	&	0	&	1.2\\
			1		&	0.1	&	0	&	1\\
			\hline
			0		&	0.2	&	0	&	1\\
			0.52	&	0.2	&	0	&	1.4\\
			1		&	0.2		&	0	&	1\\
			\hline
			0		&	$\frac{11}{30}$	&	0	&	1\\
			0.5		&	0.32		&	0	&	0.8\\
			1		&	$\frac{11}{30}$	&	0	&	1\\
			\hline
			0		&	$\frac{13}{30}$	&	0	&	1\\
			0.4		&	0.45	&	0	&	1\\
			1		&	$\frac{13}{30}$	&	0	&	1\\
			\hline
			0		&	$\frac{8}{15}$	&	0	&	1\\
			0.42	&	0.55	&	0	&	1.3\\
			1		&	$\frac{8}{15}$	&	0	&	1\\
			\hline
			0		&	0.6	&	0	&	1\\
			0.56	&	0.69	&	0	&	1.1\\
			1		&	0.6	&	0	&	1.0\\
			\hline
			0		&	$\frac{23}{30}$	&	0	&	1\\
			0.55	&	0.8	&	0	&	1.5\\
			1		&	$\frac{23}{30}$	&	0	&	1\\
			\hline
			0		&	0.9	&	0	&	1\\
			0.5		&	0.95	&	0	&	0.9\\
			1		&	0.9	&	0	&	1\\
			\hline
			0		&	1	&	0	&	1\\
			0.5		&	1	&	0	&	1\\
			1		&	1	&	0	&	1\\
			\hline
			\hline
		\end{tabular}
		\caption{Multi-patch-geometry with various continuities (Sec.~\ref{sec:example:mp:multicon})}
	\end{subtable}
	\caption{Coordinates and weights of the control points of the investigated multi-patch meshes}
	\label{tab:geo_mp}
\end{table}
\FloatBarrier

\section{Supplementary data}
\label{supp_data}
All IGES files of the employed discretizations are provided within a zip file as supplementary material to the article on the journal's webpage.

  \bibliographystyle{elsarticle-num}
  \bibliography{References.bib}


%
%
%
\end{document}